\documentclass[a4paper,11pt,twoside]{amsart}

\usepackage{latexsym}
\usepackage{amssymb}
\usepackage{amscd}
\usepackage[all]{xy}
\usepackage[mathscr]{euscript}
\usepackage{dsfont}
\usepackage{hyperref}
\usepackage{graphicx}
\usepackage{amsfonts}
\usepackage{amsmath}
\usepackage{amsthm}
\usepackage{aliascnt}

\RequirePackage{xspace}

\newcommand{\thought}[1]{}
\renewcommand{\thought}[1]{ \textbf{[#1]}}

\usepackage{enumerate}        
\newenvironment{roenumerate}{\begin{enumerate}[\upshape (i)]}{\end{enumerate}}

\newcommand\nc {\newcommand}
\newcommand\rnc{\renewcommand}

\newcount\blopone
\newcount\xone
\newcount\xtwo
\newcount\ytwo

\newtheorem{theorem}{Theorem}[subsection]
\newtheorem{thmInt}{Theorem}[section]

\newtheorem{thmApp}{Theorem}[section]
\newaliascnt{prop}{theorem}
\newtheorem{prop}[prop]{Proposition}
\aliascntresetthe{prop}
\newtheorem{com}[theorem]{Comment}
\newtheorem{apl}[theorem]{Application}
\newtheorem{exercise}[theorem]{Exercise}
\newtheorem{redu}[theorem]{Reduction}
\newtheorem{refinement}[theorem]{Refinement}
\newtheorem{summary}[theorem]{Summary}
\newtheorem{importnota}[theorem]{Important Notation}
\newtheorem{prblm}[theorem]{Problem}
\newaliascnt{notation}{theorem}
\newtheorem{notation}[notation]{Notation}
\aliascntresetthe{notation}
\newaliascnt{notationApp}{thmApp}

\aliascntresetthe{notationApp}
\newtheorem{explanation}[theorem]{Explanation}
\newaliascnt{defin}{theorem}
\newtheorem{defin}[defin]{Definition}
\aliascntresetthe{defin}
\newaliascnt{definApp}{thmApp}
\newtheorem{definApp}[definApp]{Definition}
\aliascntresetthe{definApp}
\newtheorem{caution}[theorem]{Caution}
\newaliascnt{remark}{theorem}
\newtheorem{remark}[remark]{Remark}
\aliascntresetthe{remark}
\newaliascnt{remarkApp}{thmApp}
\newtheorem{remarkApp}[remarkApp]{Remark}
\aliascntresetthe{remarkApp}
\newaliascnt{reminder}{theorem}
\newtheorem{reminder}[reminder]{Reminder}
\aliascntresetthe{reminder}
\newtheorem{illustration}[theorem]{Illustration}
\newtheorem{observation}[theorem]{Observation}
\newaliascnt{lemma}{theorem}
\newtheorem{lemma}[lemma]{Lemma}
\aliascntresetthe{lemma}
\newaliascnt{lemmaApp}{thmApp}
\newtheorem{lemmaApp}[lemmaApp]{Lemma}
\aliascntresetthe{lemmaApp}
\newtheorem{construction}[theorem]{Construction}
\newtheorem{conjecture}[theorem]{Conjecture}
\newtheorem{discussion}[theorem]{Discussion}
\newaliascnt{corollary}{theorem}
\newtheorem{corollary}[corollary]{Corollary}
\aliascntresetthe{corollary}
\newaliascnt{example}{theorem}
\newtheorem{example}[example]{Example}
\aliascntresetthe{example}
\newtheorem{conclusion}[theorem]{Conclusion}
\newtheorem{sketch}[theorem]{Sketch}
\newtheorem{triviality}[theorem]{Triviality}
\newtheorem{proto}[theorem]{Prototype Quasifibration}
\newtheorem{cauex}[theorem]{Cautionary Example}
\newaliascnt{hypo}{theorem}
\newtheorem{hypo}[hypo]{Hypothesis}
\aliascntresetthe{hypo}
\newtheorem{subth}{ }[theorem]
\newtheorem{case}{Case}[theorem]
\newtheorem{ssubth}{ }[subth]
\newtheorem{facts}[theorem]{Facts}
\newtheorem{history}[theorem]{Historical Survey}
\newtheorem{proofs}[theorem]{Discussion of the Proofs, Old and New}
\newtheorem{heuristic}[theorem]{Heuristic}

\newcommand{\st}{\mid} 
\nc\tri[1]{\begin{triviality}
\label{#1}}
\nc\fac[1]{\begin{facts}
\label{#1}
\begin{em}}
\nc\heu[1]{\begin{heuristic}
\label{#1}
\begin{em}}
\nc\app[1]{\begin{apl}
\label{#1}
\begin{em}}
\nc\skt[1]{\begin{sketch}
\label{#1}
\begin{em}}
\nc\hst[1]{\begin{history}
\label{#1}
\begin{em}}
\nc\pfs[1]{\begin{proofs}
\label{#1}
\begin{em}}
\nc\cas[1]{\begin{case}
\label{#1}
\begin{em}}
\nc\rfn[1]{\begin{refinement}
\label{#1}}
\nc\prt[1]{\begin{proto}
\label{#1}}
\nc\lem[1]{\begin{lemma}
\label{#1}}
\nc\pro[1]{\begin{prop}
\label{#1}}
\nc\thm[1]{\begin{theorem}
\label{#1}}
\nc\dis[1]{\begin{discussion}
\label{#1}
\begin{em}}
\nc\cor[1]{\begin{corollary}
\label{#1}}
\nc\dfn[1]{\begin{defin}
\label{#1}}
\nc\sthm[1]{\begin{subth}
\label{#1}}
\nc\exm[1]{\begin{example}
\label{#1}
\begin{em}}
\nc\obs[1]{\begin{observation}
\label{#1}
\begin{em}}
\nc\plm[1]{\begin{prblm}
\label{#1}
\begin{em}}
\nc\rmk[1]{\begin{remark}
\label{#1}
\begin{em}}
\nc\rmd[1]{\begin{reminder}
\label{#1}
\begin{em}}
\nc\ntn[1]{\begin{notation}
\label{#1}
\begin{em}}
\nc\exe[1]{\begin{exercise}
\label{#1}
\begin{em}}
\nc\xpl[1]{\begin{explanation}
\label{#1}
\begin{em}}
\nc\smr[1]{\begin{summary}
\label{#1}
\begin{em}}
\nc\cau[1]{\begin{caution}
\label{#1}
\begin{em}}
\nc\hyp[1]{\begin{hypo}
\label{#1}}
\nc\imn[1]{\begin{importnota}
\label{#1}
\begin{em}}
\nc\rdn[1]{\begin{redu}
\label{#1}
\begin{em}}
\nc\cax[1]{\begin{cauex}
\label{#1}
\begin{em}}
\nc\cmt[1]{\begin{com}
\label{#1}
\begin{em}}
\nc\con[1]{\begin{construction}
\label{#1}
\begin{em}}
\nc\cnj[1]{\begin{conjecture}
\label{#1}
\begin{em}}
\nc\ill[1]{\begin{illustration}
\label{#1}
\begin{em}}
\nc\ssthm[1]{\begin{ssubth}
\label{#1}
\begin{em}}
\nc\cnc[1]{\begin{conclusion}
\label{#1}
\begin{em}}

\nc\elem{\end{lemma}}
\nc\erdn{\end{em}\end{redu}}
\nc\erfn{\end{refinement}}
\nc\eprt{\end{proto}}
\nc\ethm{\end{theorem}}
\nc\ecor{\end{corollary}}
\nc\edfn{\end{defin}}
\nc\esthm{\end{subth}}
\nc\epro{\end{prop}}
\nc\etri{\end{triviality}}
\nc\eexm{\end{em}
\end{example}}
\nc\eobs{\end{em}
\end{observation}}
\nc\ecmt{\end{em}
\end{com}}
\nc\efac{\end{em}
\end{facts}}
\nc\eheu{\end{em}
\end{heuristic}}
\nc\eapp{\end{em}
\end{apl}}
\nc\ermk{\end{em}
\end{remark}}
\nc\ermd{\end{em}
\end{reminder}}
\nc\eill{\end{em}
\end{illustration}}
\nc\eplm{\end{em}
\end{prblm}}
\nc\ecas{\end{em}
\end{case}}
\nc\eskt{\end{em}
\end{sketch}}
\nc\ecau{\end{em}
\end{caution}}
\nc\ecax{\end{em}
\end{cauex}}
\nc\eimn{\end{em}
\end{importnota}}
\nc\entn{\end{em}
\end{notation}}
\nc\eexe{\end{em}
\end{exercise}}
\nc\expl{\end{em}
\end{explanation}}
\nc\edis{\end{em}
\end{discussion}}
\nc\econ{\end{em}
\end{construction}}
\nc\ecnj{\end{em}
\end{conjecture}}
\nc\esmr{\end{em}
\end{summary}}
\nc\ehst{\end{em}
\end{history}}
\nc\epfs{\end{em}
\end{proofs}}
\nc\ehyp{
\end{hypo}}
\nc\ecnc{\end{em}
\end{conclusion}}
\nc\essthm{\end{em}
\end{ssubth}}

\nc\sst{\scriptstyle}
\newcommand{\comment}[1]{}
\newcommand{\ri}{\longrightarrow}

\newcommand{\zz}{{\mathbb Z}}

\newcommand{\nn}{{\mathbb N}}
\newcommand{\K}{{\mathbf K}}
\newcommand{\D}{{\mathbf D}}
\newcommand{\qq}{{\mathbb Q}}

\nc\op{^{\hbox{\rm\tiny op}}}
\nc\mth{^{\hbox{\rm\tiny th}}}

\nc\script{\mathscr}
\nc\z{\zeta}
\nc\bc{\mathrm{bc}}
\nc\ct{{\script T}}
\nc\cf{{\script F}}
\nc\cg{{\script G}}
\nc\ch{{\script H}}
\nc\ck{{\script K}}
\nc\cl{{\script L}}
\nc\cv{{\script V}}
\nc\ce{{\script E}}
\nc\cs{{\script S}}
\nc\car{{\script R}}
\nc\cd{{\script D}}
\nc\cc{{\script C}}
\nc\ca{{\script A}}
\nc\ci{{\script I}}
\nc\cj{{\script J}}
\nc\co{{\script O}}
\nc\cu{{\script U}}
\nc\cx{{\script X}}
\nc\Cp{{\script P}}
\nc\cq{{\script Q}}
\nc\cy{{\script Y}}
\nc\cz{{\script Z}}
\nc\bd{\begin{description}}
\nc\ed{\end{description}}
\nc\ctob{{\script C}at\big(\ci^{op},\ca\big)}
\nc\clim{{\ds\mathop{\rm lim}_{\ds\longleftarrow}}\,}
\nc\climi{\clim_{\!i}\,}
\nc\climn{\clim^{\!n}\,}
\nc\colim{{\ds\mathop{\rm colim}_{\ds\la}}}
\nc\colimj{{\ds\mathop{\rm colim}_{\ds\la}}{}_{j\,}}
\nc\oa{\overline{\ca}}
\nc\s{\sigma}
\nc\ta{\tau}
\nc\os{\overline\sigma}
\nc\ot{\overline\tau}
\nc\T{\Sigma}
\nc\Tm{\Sigma^{-1}}
\nc\de[1]{{\mathop{\rm deg(#1)}}}
\nc\Ad[1]{\mathop{\rm Ad}(#1)}
\nc\ad[1]{\mathop{\rm ad}(#1)}
\nc\kth{{\it K}--theory}
\nc\loc[1]{{\text{\rm Loc(#1)}}}
\nc\coloc[1]{{\text{\rm Coloc}(#1)}}

\newcommand{\sh}[2][]{\Sigma^{#1}#2}

\def\der #1 {D\left(#1\right)}
\nc\prf{\begin{proof}}
\nc\eprf{\end{proof}}
\nc\ds{\displaystyle}
\nc\Tor{\text{\rm Tor}}

\nc\cb{{\script B}}
\nc\ab{{\script A}b}

\nc\be{\begin{roenumerate}}
\nc\ee{\end{roenumerate}}

\nc\cat[1]{{\script C}at\Big({\big\{#1\big\}}\op\,\,,\,\,\ab\Big)}
\nc\csab{{\script C}at\big(\cs^{op},\ab\big)}
\nc\ctab{{\script C}at\Big({\{\ct^\alpha\}}^{op},\ab\Big)}
\nc\csex{{\script E}x\big(\cs^{op},\ab\big)}
\nc\ctex{{\script E}x\Big({\{\ct^\alpha\}}^{op},\ab\Big)}
\nc\sub{\qquad\subseteq\qquad}
\nc\ctr[1]{{\left.\ct\left(-,#1\right)\right|}_{\cs}}
\nc\ctrf[2]{{\left.\ct\left(#1,#2\right)\right|}_{\cs}}
\nc\Ctr[1]{{\left.\ct\left(-,#1\right)\right|}_{\ct^\alpha}}
\nc\Ctrf[2]{{\left.\ct\left(#1,#2\right)\right|}_{\ct^\alpha}}

\nc\la{\longrightarrow}
\nc\nin{\noindent}
\nc\cad[1]{\text{card}(#1)}
\nc\eq{\quad=\quad}
\nc\BA{\begin{array}{c}}
\nc\EA{\end{array}}
\nc\barr{
\[
\begin{array}{cccccccccccccccc}
}
\nc\earr{
\end{array}
\]
}
\nc\as[1]{{\langle S\rangle}^{#1}}
\nc\shi{\text{\it shift}}

\nc\yy[1]{{\left.\ct\left(-,#1\right)\right|}_{\ct^c}}
\nc\vrep[2]{{\left.\ct\left(#1,#2\right)\right|}_{\ct^\alpha}}
\nc\da{\downarrow}
\nc\Hom{{\mathop{\rm Hom}}}
\nc\HHom{{\script H}{\mathop{\rm om}}}
\nc\End{{\mathop{\rm End}}}
\nc\Ext{{\mathop{\rm Ext}}}
\nc\mr\Modtc
\nc\PExt{{\mathop{\rm PExt}}}
\nc\stm{\text{\rm stmod}(kG)}
\nc\stM{\text{\rm StMod}(kG)}
\nc\e{\varepsilon}
\nc\p{\varphi}

\nc\rs{\s^{-1}A}
\nc\br{{\{\s^{-1}A\}}}
\nc\ra\ri
\nc\y[1]{\mathbf{y}#1}
\nc\x[1]{\mathbf{z}#1}
\nc\mmod[1]{#1\text{--\rm mod}}
\nc\Mod[1]{#1\text{--\rm Mod}}
\nc\MMod[1]{\text{\rm Mod--}#1}
\nc\Md {\ensuremath{\mathop{\textup{Mod}}}}
\rnc\mod[1]{\ensuremath{\mathop{#1\textup{--mod}}}\xspace}
\nc\Modtc{\Mod{\ct^c}}
\nc\pgldim[1]{\mathop{\rm pgldim}\,#1}
\nc\tf{{\rm [TR5]}}
\nc\tfs{{\rm [TR5$^*$]}}
\nc\Fun{\text{\rm Funct}(F\op,\ab)}
\nc\sym{\text{\rm Sym}}
\nc\sgn{\text{\rm sgn}}
\nc\Pro{\text{\rm Prod}^{}_\alpha(F\op,\ab)}
\nc\Yt[1]{{\left.\Hom_\ct^{}\left(-,#1\right)\right|}_F^{}}
\nc\dl{\delta}
\nc\Proj[1]{#1\text{--\rm Proj}}
\nc\proj[1]{#1\text{--\rm proj}}
\nc\Flat[1]{#1\text{--\rm Flat}}
\nc\Inj[1]{#1\text{--\rm Inj}}
\nc\Ima{\mathrm{Im}}
\nc\Ker{\mathrm{Ker}}
\nc\ov{\overline}
\nc\wt{\widetilde}
\nc\wh{\widehat}
\nc\ph{\varphi}
\nc\tstr{{\it t}--structure}
\nc\tstrs{{\it t}--structure }
\nc\spec[1]{{\text{\rm Spec}\left(#1\right)}}

\nc\EProd{\text{\rm EProd}}
\nc\ECoprod{\text{\rm ECoprod}}
\nc\Prod{\text{\rm Prod}}
\nc\ldimp{\text{\rm LDim}^{\prod}}
\nc\ldimc{\text{\rm LDim}^{\coprod}}
\nc\gen[2]{{\langle#1\rangle}^{}_{#2}}
\nc\genu[3]{{\langle#1\rangle}^{[#3]}_{#2}}
\nc\ogen[1]{\ov{\langle#1\rangle}}
\nc\ogenun[2]{\ov{\langle#1\rangle}_{#2}^{}}
\nc\ogenu[3]{\ov{\langle#1\rangle}^{[#3]}_{#2}}
\nc\ogenul[3]{\ov{\langle#1\rangle}^{(-\infty,#3]}_{#2}}
\nc\ogenuf[3]{\ov{\langle#1\rangle}^{[#3,\infty)}_{#2}}
\nc\genuf[3]{{\langle#1\rangle}^{[#3,\infty)}_{#2}}
\nc\genul[3]{{\langle#1\rangle}^{(-\infty,#3]}_{#2}}
\nc\dperf[1]{\D^{\mathrm{perf}}(#1)}
\nc\dcoh{\mathbf{D}^b_{\mathrm{coh}}}
\newcommand{\Dqc}{{\mathbf D_{\text{\bf qc}}}}

\newcommand{\Dqcmis}[1]{{\mathbf D_{\text{\bf qc},#1}^-}}
\newcommand{\Dqcpls}[1]{{\mathbf D_{\text{\bf qc},#1}^+}}
\newcommand{\Dqcbs}[1]{{\mathbf D_{\text{\bf qc},#1}^b}}
\newcommand{\Dqcps}[1]{{\mathbf D_{\text{\bf qc},#1}^p}}
\newcommand{\Dqcpbs}[1]{{\mathbf D_{\text{\bf qc},#1}^{p,b}}}
\nc\dmcoh{\mathbf{D}^-_{\mathrm{coh}}}
\nc\dscoh{\mathbf{D}^{}_{\mathrm{coh}}}
\nc\RHHom{{\script{RH}}{\mathrm{om}}}
\nc\Coprod{\mathrm{Coprod}}
\nc\COprod{\mathrm{coprod}}
\nc\add{\mathrm{add}}
\nc\Add{\mathrm{Add}}
\nc\Smr{\mathrm{smd}}
\nc\id{\mathrm{id}}
\nc\LL{\mathbf{L}}
\nc\R{\mathbf{R}}
\nc\wi{\wt{\text{\it\i}}}
\nc\exal{\ce\text{\it x}(\ct^\alpha,\ab)}
\nc\exalz{\ce\text{\it x}_{\aleph_0}^{}(\ct^\alpha,\ab)}
\nc\fc{\mathfrak{C}}
\nc\fl{\mathfrak{L}}
\nc\fs{\mathfrak{S}}
\nc\Prf{\text{\bf Perf}}
\nc\Enh{\mathbf{Enh}}
\nc\Tri{\mathbf{Tri}}
\nc\fgt{\mathsf{Fgt}}
\newcommand{\Dqcs}[1]{{\mathbf D_{\text{\bf qc},#1}}}
\nc\dcohs[1]{\mathbf{D}^b_{\mathrm{coh},#1}}
\nc\coh{\mathbf{Coh}}
\nc\qc{\mathbf{Qcoh}}
\nc\vect{\mathbf{Vect}}
\nc\dperfs[2]{\D_{#1}^{\mathrm{perf}}(#2)}

\newcommand{\fun}[1]{\mathsf{#1}}

\newcommand{\fF}{\fun{F}}

\newcommand{\fH}{\fun{H}}
\newcommand{\fI}{\fun{I}}
\newcommand{\fJ}{\fun{J}}

\nc\hoco{
\begin{picture}(40,10)
\put(20,0){\makebox(0,0)[b]{\text{\rm Hocolim}}}
\put(5,-2){\vector(1,0){30}}
\end{picture}\,\,}

\renewcommand{\leq}{\leqslant}
\renewcommand{\geq}{\geqslant}
\nc\tst[1]{\left({#1}^{\leq0},{#1}^{\geq1}\right)}
\nc\tstv[2]{\left({#1}_{#2}^{\leq0},{#1}_{#2}^{\geq1}\right)}
\nc\tsth[2]{{#1}_{#2}^{\heartsuit}}

\nc\holim{
\begin{picture}(40,10)
\put(20,0){\makebox(0,0)[b]{\text{\rm Holim}}}
\put(35,-2){\vector(-1,0){30}}
\end{picture}}

\makeatletter
\let\@wraptoccontribs\wraptoccontribs
\makeatother
\contrib[with an appendix by]{Christian Haesemeyer}

\begin{document}

\author{Alberto Canonaco, Amnon Neeman, and Paolo Stellari}

\address{A.C.: Dipartimento di Matematica ``F.\ Casorati''\\
        Universit{\`a} degli Studi di Pavia\\
        Via Ferrata 5\\
        27100 Pavia\\
        ITALY}
\email{alberto.canonaco@unipv.it}

\address{C.H.: School of Mathematics and Statistics\\The University of Melbourne\\Parkville VIC
3010\\ AUSTRALIA}
\email{christian.haesemeyer@unimelb.edu.au}

\address{A.N.: Dipartimento di Matematica ``F.\ Enriques''\\
        Universit{\`a} degli Studi di Milano\\
        Via Cesare Saldini 50\\
	20133 Milano\\
        ITALY}
\email{amnon.neeman@unimi.it}

\address{P.S.: Dipartimento di Matematica ``F.\ Enriques''\\
        Universit{\`a} degli Studi di Milano\\
        Via Cesare Saldini 50\\
	20133 Milano\\
        ITALY}
\email{paolo.stellari@unimi.it}
\urladdr{\url{https://sites.unimi.it/stellari}}

 \thanks{A.~C.~was partially supported by the research project PRIN 2022 ``Moduli spaces and special varieties''. A.~N.~was partly supported 
   by Australian Research Council Grants DP200102537 and DP210103397,
   and by ERC Advanced Grant 101095900-TriCatApp.
    	P.~S.~was partially supported by the ERC Consolidator Grant ERC-2017-CoG-771507-StabCondEn, by the research project PRIN 2022 ``Moduli spaces and special varieties'', and by the research project FARE 2018 HighCaSt (grant number R18YA3ESPJ)}

 \title[Subcategories of weakly approximable categories]
       {The passage among the subcategories of
   weakly approximable triangulated categories}

\begin{abstract}
In this article we prove that all the inclusions between the `classical' and naturally defined full triangulated subcategories of a weakly approximable triangulated category are intrinsic (in one case under a technical condition). This extends all the existing results about subcategories of weakly approximable triangulated categories. 

Together with a forthcoming paper about uniqueness of enhancements, our result allows us to generalize a celebrated theorem by Rickard which asserts that if $R$ and $S$ are left coherent rings, then a derived equivalence of $R$ and $S$ is
``independent of the decorations''.
That is, if $D^?(R\text{--}\square)$ and $D^?(S\text{--}\square)$
are equivalent as triangulated categories for some choice of decorations $?$ and $\square$, then they are equivalent
for every choice of decorations. But our theorem is much more general, and applies also to quasi-compact and quasi-separated schemes---even to the relative version, in which the derived categories consist of complexes with cohomology supported on a given closed subscheme with quasi-compact complement.
\end{abstract}



\maketitle

\setcounter{tocdepth}{1}
\tableofcontents

\setcounter{section}{0}

\section*{Introduction}
\label{S0}

Derived categories are by now old and well-established, but it helps
occasionally to remember the difficulties that the people who first
used them had to overcome. Let us therefore go back in time
to the 1960s and 70s, when the subject was still new.

The first printed introduction to them was in Hartshorne's 1966
exposition
of Grothendieck's duality theory~\cite{Hartshorne66}, but the fuller
and more thorough treatment in Verdier's PhD thesis circulated among
the experts at the time---it was finally published in 1996 (see~\cite{Verdier96}).
And the reader interested in the historical
perspective is encouraged to also look at Illusie's expos\'es
in SGA6~\cite{Illusie71A,Illusie71B,Illusie71C}, which
all address foundational questions about derived categories.
There was an enormous effort in these early days to try
to understand these objects, and figure out how best to
work with them.

Recall that, associated to a scheme $X$, it is
customary to attach three exact or abelian
categories: the category of quasi-coherent sheaves on
$X$, denoted here
$\qc(X)$, the subcategory of vector bundles on $X$,
denoted here $\vect(X)$,
and (if $X$ is noetherian) the category of coherent
sheaves on $X$, which
we will denote $\coh(X)$. Now when we pass from abelian
(or exact) categories to their derived categories,
the objects
of interest are cochain complexes of objects in the
abelian (or exact)
category, and the natural question becomes which
cochain complexes
should be permitted. And in the early days this was not
clear---a great deal of thought and effort went into
studying the advantages and disadvantages
of the various options.

After several decades, the consensus is as follows.
Let $X$ be a scheme, and assume $Z\subseteq X$ is
a Zariski-closed subset. For the derived categories
of (quasi-)coherent sheaves to work well, it is
best to impose finiteness hypotheses on the
schemes that arise in the constructions
we wish to make, and
the finiteness hypotheses that yield a
useful theory postulate that the
scheme $X$ is quasi-compact and quasi-separated\footnote{Recall that
a scheme $X$ is quasi-separated if
the intersection of any two quasi-compact open
subsets of $X$ is quasi-compact.}, and that the
open set $X\setminus Z$ is quasi-compact. And
the objects of the seven useful
(relative) derived categories, associated
to the pair $Z\subseteq X$, should all be
cochain complexes of sheaves of $\co_X^{}$--modules,
whose restriction to the open set $X\setminus Z$ is acyclic.
And then we may wish to impose
further restrictions, for
example on the allowed cohomology sheaves.
It so happens that the most
useful of the various derived categories one could
consider turn out to be
the following seven:
\be
\item[(a)]
  $\Dqcs Z(X)$. The only restriction is that all the
  cohomology sheaves
must be quasi-coherent.
\item[(b)]
$\Dqcmis Z(X)\subseteq\Dqcs Z(X)$. This is the
full subcategory of all
complexes $\cc$ with $\ch^n(\cc)=0$ for $n\gg0$.
\item[(c)]
$\Dqcpls Z(X)\subseteq\Dqcs Z(X)$. This is the
full subcategory of all
complexes $\cc$ with $\ch^n(\cc)=0$ for $n\ll0$.
\item[(d)]
$\Dqcbs Z(X)\subseteq\Dqcs Z(X)$. This is the
full subcategory of all
complexes $\cc$ with $\ch^n(\cc)=0$ for
all but finitely many $n\in\zz$.
\item[(e)]
$\dperfs ZX\subseteq\Dqcs Z(X)$. This is the
full subcategory of all
\emph{perfect} complexes supported on $Z$.
A complex is perfect if it is locally
isomorphic to a bounded complex of
finite-rank vector
bundles.
\item[(f)]
$\Dqcps Z(X)\subseteq\Dqcs Z(X)$. This is the
full subcategory of all
\emph{pseudocoherent} complexes supported on $Z$.
A complex is pseudocoherent if it is locally
isomorphic to a bounded-above complex of
finite-rank vector
bundles.
\item[(g)]
$\Dqcpbs Z(X)\subseteq\Dqcps Z(X)$. This is the
full subcategory of all
objects $\cc\in\Dqcps Z(X)$
such that $\ch^n(\cc)=0$ for all but finitely
many $n\in\zz$.  
\ee

\nin
Note that the finiteness conditions on the complex,
imposed in (e), (f) and (g), are all due to
Illusie~\cite{Illusie71C}. The definition of
(f) in Illusie is slightly more complicated than
what we presented, but the two definitions
are equivalent---this can easily be seen by using
\cite[Theorem~5.1]{Bokstedt-Neeman93}. In fact, any
of the three finiteness conditions (e), (f) and (g)
is local in the flat topology. And for an
affine scheme $X=\spec R$,
where by \cite[Theorem~5.1]{Bokstedt-Neeman93}
we know that the natural map $\D(\Mod R)\la\Dqc(X)$
is an equivalence, the complex
$\cc\in\Dqcs Z(X)\subseteq\Dqc(X)$ will belong to
the subcategories (e), (f) or (g) provided it is
isomorphic in $\D(\Mod R)$ to
a complex of finitely
generated, projective $R$--modules which is
(respectively) bounded, bounded-above, or bounded-above
with only finitely many nonvanishing cohomology groups.

And one issue that concerned the early workers in
the field was that the derived Hom involved injective
resolutions, while the derived tensor product involved
flat resolutions. And in the early days injective
resolutions were only known to exist in $\Dqcpls Z(X)$,
and flat resolutions only in $\Dqcmis Z(X)$. To put things
in historical perspective: it was not until Spaltenstein's
1988 article~\cite{Spaltenstein88} that anyone came up with
an adequate approach to forming derived tensor products and
derived Homs in $\Dqcs Z(X)$ (although in some sense
the homotopy theorists arrived at a different
satisfactory approach to a
parallel problem earlier, as
discussed in \cite{Bokstedt-Neeman93}).
Now especially in the case
of a subject like Grothendieck's duality theory, where
derived tensor products and
derived Homs both occur and intermingle,
this caused headaches.

To the early workers in the field it seemed crucial
to find the right derived category for the problem at
hand. Using derived categories was viewed as an art,
and a good artist displayed her competence by choosing
wisely the derived category to work with. In fact: to an extent
this attitude persists to the present day, in birational
geometry. The people using derived category techniques
have been known to argue over the relative merits
of the categories $\dperf X$ and $\dcoh(X)$.

Against this background
one can imagine how surprising was
the 
work of Rickard's, which appeared in 1989 and 1991
in the two articles~\cite{Rickard89b,Rickard91},
and proves:

\begin{thmInt}[Rickard]\label{T0.1}
Let $R$ and $S$ be two rings. Then,
in the standard notation for the various
derived categories associated to the two rings,
the following are
equivalent:
\be
\item
There exists a triangulated equivalence $\D(\Mod R)\cong\D(\Mod S)$.
\item
There exists a triangulated equivalence $\D^-(\Mod R)\cong\D^-(\Mod S)$.
\item
There exists a triangulated equivalence $\D^+(\Mod R)\cong\D^+(\Mod S)$.
\item
There exists a triangulated equivalence $\D^b(\Mod R)\cong\D^b(\Mod S)$.
\item
There exists a triangulated equivalence $\D^-(\proj R)\cong\D^-(\proj S)$.
\item
There exists a triangulated equivalence $\D^b(\proj R)\cong\D^b(\proj S)$.
\setcounter{enumiv}{\value{enumi}}
\ee
If we assume further that the rings
$R$ and $S$ are both left coherent, then the six
equivalent conditions above are also equivalent to:
\be
\setcounter{enumi}{\value{enumiv}}
\item
There exists a triangulated equivalence $\D^b(\mmod R)\cong\D^b(\mmod S)$.%
\footnote{Rickard's paper claims that it suffices for just \emph{one} of
the two rings $R,S$ to be left coherent---with
suitable adjustments made to the definition of $\D^b(\mmod R)$
for the non-left coherent ring. But the authors couldn't
follow the argument given in the old, published paper,
and neither could Rickard.
}
\ee
\end{thmInt}

\nin
If $R$ and $S$ are commutative this says that
a triangulated equivalence
$\D_{\Box,W}^?(X)\cong\D_{\Box,Z}^?(Y)$,
in the case
where $W=X=\spec R$ and $Z=Y=\spec S$ and
$\Box$ and $?$ are any of Illusie's decorations in the list (a)--(g) above, implies the equivalence of
all other pairs.

It should be explained that Rickard's theory works by
studying what an equivalence might look like. An
equivalence $\D^?(R\text{--}\Box)\cong\D^?(S\text{--}\Box)$
must take the object $R\in\D^?(R\text{--}\Box)$ to
some cochain complex, an object in $\D^?(S\text{--}\Box)$.
And the basic idea of the theory is that this complex
must be very special---it has come to be known as
a \emph{tilting complex.} And roughly the idea is that
the tilting complexes are the same,  independent
of the decoration $?$ and $\Box$. 

\subsection*{The result}

In this article and its sequel we give a vast
improvement and a vast
generalization of Rickard's remarkable
result. First of all the
improvement: we show that each of the seven derived
categories on Rickard's
list determines all the others. We mean
this in the
precise sense that, for each ordered pair $\ca,\cb$
of the derived categories, on Rickard's list
of seven, there is an
explicit recipe that takes
the triangulated category $\ca$ as
input and outputs the
triangulated category $\cb$. Thus the seven
derived categories
are interchangeable, each of them knows all
about the other six.

So much for the improvement, now the time has come
for the generalization. Our result is not only
about derived categories of rings, it is about
weakly approximable triangulated categories
with unique enhancements. And now it is time to
remind the reader of the terminology.

\smallskip

We begin by recalling some work by the second
author. In a series of recent articles he developed the
notion of ``approximable'' and ``weakly approximable''
triangulated categories---there will be a brief review
in \autoref{S1}. And for the purposes of this article, the
important features are that every
weakly approximable triangulated
category $\ct$ comes with intrinsically defined
\be
\item[(1)]
A preferred equivalence class of \tstr s, see
\autoref{R1.1} for more detail.
\item[(2)]
Thick subcategories $\ct^-$, $\ct^+$, $\ct^b$, $\ct^-_c$,
$\ct^c$, $\ct^b_c=\ct^-_c\cap\ct^b$, and
$\ct^{c,b}=\ct^c\cap\ct^b_c$. We will recall the definition
of each of these subcategories in
\autoref{the subcategories we care about}.
\setcounter{enumiv}{\value{enumi}}
\ee
And the relevance of all of this to Rickard's old theorem
comes from two facts.
\be
\setcounter{enumi}{\value{enumiv}}
\item[(3)]
If $R$ is any ring, then the category $\ct=\D(\Mod R)$
is (weakly) approximable. Also, if $X$ is a quasi-compact and quasi-separated scheme and $Z\subseteq X$
is a Zariski-closed subset with quasi-compact complement, then the
category $\ct=\Dqcs Z(X)$ is weakly approximable.
We will give references in
\autoref{ex:wa}.
\item[(4)]
The general intrinsic subcategories of (2),
of any weakly approximable
triangulated category $\ct$, 
can be computed when $\ct$ is one the two special cases
given in (3).
They turn out to be
\begin{center}
	\begin{tabular}{|c|c|c|}
		\hline
		& $Z\subseteq X$ as in (a)--(g)  & $R$ a ring\\
		\hline
		$\ct$ & $\Dqcs Z(X)$ & $\D(\Mod{R})$\\
		$\ct^-$ & $\Dqcmis Z(X)$ & $\D^-(\Mod{R})$\\
		$\ct^+$ & $\Dqcpls Z(X)$ & $\D^+(\Mod{R})$\\
		$\ct^b$ & $\Dqcbs Z(X)$ & $\D^b(\Mod{R})$\\
		$\ct^c$ & $\dperfs ZX$ & $\K^b(\proj R)$\\
		$\ct^-_c$ & $\Dqcps Z(X)$ & $\K^-(\proj{R})$\\
		$\ct^b_c$ & $\Dqcpbs Z(X)$ & $\K^{-,b}(\proj{R})$\\
  $\ct^{c,b}$ & $\dperfs ZX$ & $\dperf{R}$\\
		\hline
	\end{tabular}
\end{center}
In other words, in the case $\ct=\Dqcs Z(X)$ we
recover Illusie's old list (a)--(g), and in
the case $\ct=\D(\Mod R)$ we recover the list of
categories in Rickard's old \autoref{T0.1}.
\setcounter{enumiv}{\value{enumi}}
\ee
An immediate consequence, of the existence
of a recipe that produces these subcategories
out of $\ct$, is that, in \autoref{T0.1},
(i) implies all of (ii), (iii), (iv), (v), (vi) and (vii).

\smallskip

In summary, if $\ct$ is a weakly approximable triangulated category, then the subcategories listed in (2) sit in the following commutative diagram:
\begin{equation}\label{eq:incl}
\xymatrix@C+40pt@R-5pt{
&\ct&\\
\ct^-\ar@{^{(}->}[ur]&\ct^b\ar@{^{(}->}[u]\ar@{_{(}->}[l]\ar@{^{(}->}[r]&\ct^+\ar@{_{(}->}[ul]\\
\ct^-_c\ar@{^{(}->}[u]&\ct^b_c\ar@{^{(}->}[u]\ar@{_{(}->}[l]&\\
\ct^c\ar@{^{(}->}[u]&\ct^{c,b}.\ar@{^{(}.>}[u]\ar@{_{(}->}[l]
 \ar@/_2pc/@{^{(}->}[uu]&
}
\end{equation}
The main aim of this paper is to show that all the
solid inclusions $\ca\hookrightarrow\cb$ in \eqref{eq:incl} are intrinsic, in the sense that there is a recipe, depending only on the triangulated structure on $\cb$, that describes which objects in $\cb$ belong to the full subcategory $\ca$. Since we found the notion of ``recipe'' difficult to formulate precisely, we state our main theorem as follows.

\begin{thmInt}\label{thm:main1}
All the inclusions in diagram \eqref{eq:incl}, represented by the solid arrows $\ca\hookrightarrow\cb$, are \emph{invariant under triangulated equivalences}. By this we mean: given a pair of weakly approximable triangulated
categories $\ct,\ct'$, as well as matching inclusions
$\ca\hookrightarrow\cb\hookrightarrow\ct$ and
$\ca'\hookrightarrow\cb'\hookrightarrow\ct'$ from diagram \eqref{eq:incl},
then any triangulated equivalence $\cb\la\cb'$ must restrict to a triangulated equivalence $\ca\la\ca'$.

Furthermore, the same is true for the (unique) inclusion in diagram \eqref{eq:incl} represented by a dotted arrow, provided we further assume that one of the two conditions below holds.
\be
\item $\ct,\ct'$ are \emph{coherent,} as in \autoref{dfn:coherent}, or
\item $\ct^c\subseteq\ct^b_c$, ${\ct'}^c\subseteq{\ct'}^b_c$ and
$^\perp(\ct^b_c)\cap\ct^-_c=\{0\}={^\perp({\ct'}^b_c)}\cap{\ct'}^-_c$.
\ee
\end{thmInt}

\nin
Note that  \autoref{thm:main1} does not
mention enhancements,
this entire article is enhancement-free and
uses only triangulated category techniques.

The
uniqueness of enhancements enters in
the sequel, when we want to
go back. This means that,
when in the inclusion
$\ca\hookrightarrow\cb$ we want a recipe that
constructs $\cb$ out of $\ca$, then to
the best of the authors' knowledge such a recipe
needs enhancements. And when the enhancements are
unique, then so is the triangulated output. This will be discussed much more thoroughly and carefully in the sequel to the current article. Here we just observe that, for inclusions $\ca\hookrightarrow\cb$ and $\ca'\hookrightarrow\cb'$ as above, we can prove that the existence of a triangulated equivalence $\ca\to\ca'$ implies the existence of a triangulated equivalence $\cb\to\cb'$. But we do not know if every triangulated equivalence $\ca\to\ca'$ extends to a triangulated equivalence $\cb\to\cb'$, or if such extensions (when they exist) are unique. This problem is deeply related to the issue of comparing the autoequivalence groups of the various intrinsic subcategories which play a role in \autoref{thm:main1}. We will investigate this in future work.

\subsection*{Structure of the paper}

We should finish the introduction by giving an outline
of the structure of the article.

\autoref{S314} is a brief reminder of the notation
used to set up the theory of approximable triangulated
categories. \autoref{betterlabelthis} is a less-brief reminder
of compactly generated \tstr{s}. The reason that the
section on \tstr{s} is less terse is that the results
we refer to are scattered over an extensive
literature, whose focus is applications
irrelevant to us here.
Therefore the portion of this
vast literature that the reader needs,
for the current manuscript, is tiny.
\autoref{S69}
reproves  
Saor{\'{\i}}n and
{\v{S}}{\v{t}}ov{\'{\i}}{\v{c}}ek~\cite[Theorem~8.31]{Saorin-Stovicek20},
asserting that the heart of a compactly generated \tstrs is always
a locally finitely presented Grothendieck abelian category.
The reason we go to the trouble of
providing this new proof is that not only will
we be using the result, but the lemmas in our
proof will turn out to come up again in the course
of later arguments. And the preliminaries of
the paper end with \autoref{S1}, which outlines
the parts of the theory of weakly approximable
triangulated categories relevant to this manuscript.

And now we finally begin giving the recipes promised
in \autoref{thm:main1}. We even go overboard: given
any category $\cb$ in the hierarchy of diagram in
\eqref{eq:incl}, we occasionally give explicit
recipes for
one of the subcategories not immediately
below it in the hierarchy.
This will happen for example in the case of $\ct^-$.
We have inclusions $\ct^c\subseteq\ct^-_c\subseteq\ct^-$,
but it so happens that we will first give a recipe
for computing $\ct^c\subseteq\ct^-$, and then use it
to concoct a recipe for $\ct^-_c\subseteq\ct^-$.

We note also that, when there are multiple paths
in the hierarchy connecting a pair $\ca\subseteq\cb$,
the recipes they yield are not necessarily equally
complicated. We have inclusions
\[
\xymatrix{
\ct^-                     &       \ct^b\ar@{_{(}->}[l]\\
\ct^-_c\ar@{^{(}->}[u] & \ct^b_c\ar@{_{(}->}[l]\ar@{^{(}->}[u]
}
\]
and they combine to give two recipes for
$\ct^b_c\subseteq\ct^-$. As it happens the
path $\ct^b_c\subseteq\ct^-_c\subseteq\ct^-$
is much less involved than the path
$\ct^b_c\subseteq\ct^b\subseteq\ct^-$.

Back to the structure of the paper, \autoref{sect:ct-c} studies
the category $\ct^-_c$ and its subcategories.
The recipe for obtaining the
subcategory $\ct^c\subseteq\ct^-_c$ is by a trick
contained in \autoref{L2.1}, and the recipe for
$\ct^b_c\subseteq\ct^-_c$ is not difficult once
we know $\ct^c$, see the
proof of \autoref{P2.5}(ii). The same formula
also works to give a recipe for 
$\ct^{c,b}$ as a subcategory of $\ct^c$,
see \autoref{P2.5}(iii). This provides all the solid
arrows in the bottom-left square of the
hierarchy in the diagram \eqref{eq:incl}.

Now for the category $\ct^-$. In \autoref{P3.5} we
give a recipe for the preferred equivalence
class of \tstr{s} on $\ct^-$, and \autoref{L4.1} gives
a recipe for $\ct^c\subseteq\ct^-$. \autoref{cor:iota694}
combines this information to give recipes for
$\ct^b\subseteq\ct^-$ and $\ct^-_c\subseteq\ct^-$,
completing our task with $\ct^-$.

In the case of $\ct^+$, the recipe for the preferred
equivalence class of \tstr{s} can be found in
\autoref{P3.11}, and the recipe for $\ct^b\subseteq\ct^+$
that follows from it is in \autoref{cor:iota694}.

And now we come to the hardest part of the article,
providing recipes for the subcategories of $\ct^b$.
In order to study this we introduce the notion of
strongly pseudocompact objects in a triangulated category,
the reader can see the subject developed in
\autoref{S5}. And \autoref{C5.905} tells us that,
in the categories $\ct^b$ and $\ct^+$, the
strongly pseudocompact objects are precisely
$\ct^b_c$. This delivers the promised recipe
for $\ct^b_c$ as a subcategory of $\ct^b$ (and as
it happens also as a subcategory of $\ct^+$).
And the recipe giving $\ct^{c,b}$ as a subcategory
of $\ct^b$ is to be found in \autoref{L107.1}. Once
again, this recipe also works to give $\ct^{c,b}$ as
a subcategory of $\ct^+$.

And now for the dotted arrow in the hierarchy
of the diagram \eqref{eq:incl}. We have already said that we
do not have a general recipe, we only know recipes
that work in special cases. In this article
we study the special case in which $\ct^c$ is assumed
to be contained in $\ct^b$, and we find a recipe
that works as long as $^\perp(\ct^b_c)\cap\ct^-_c=\{0\}$.
See \autoref{P112.791} for the precise recipe for $\ct^c$.

In \autoref{sec:examples} we treat examples to which
\autoref{P112.791} applies. Any coherent, weakly
approximable triangulated category satisfies
$^\perp(\ct^b_c)\cap\ct^-_c=\{0\}$, in particular if $R$ is a
left coherent ring
then $\D(\Mod R)$ satisfies the hypothesis.
Applying  \autoref{P112.791}
to the situation of Rickard's \autoref{T0.1}, this
provides
a new proof that, as long as the rings $R,S$ are
left coherent, then a triangulated equivalence
as in \autoref{T0.1}(vii) implies a
triangulated equivalence as
in \autoref{T0.1}(vi). And this is one of the rare
cases where going in the reverse direction does not
require enhancements, see
\cite[Example~4.2 and Proposition~4.8]{Neeman18A}.

But perhaps more remarkable is that
$^\perp(\ct^b_c)\cap\ct^-_c=\{0\}$
holds for $\Dqcs Z(X)$ unconditionally, for
any quasi-compact, quasi-separated $X$ and any
closed subset $Z\subseteq X$ with quasi-compact
complement. Therefore in the category $\ct=\Dqcs Z(X)$
all the arrows of the diagram \eqref{eq:incl}
are solid. For the proof see \autoref{P1000.1},
which depends crucially on Haesemeyer's
\autoref{thm:appendix-main}.

\bigskip

{\small\noindent{\bf Acknowledgements.} It is our great pleasure to thank Lidia Angeleri H\"ugel and Rosanna Laking for inspiring conversations concerning the contents of \autoref{S69}, and to Jan {\v{S}}{\v{t}}ov{\'{\i}}{\v{c}}ek for
  his careful reading and many improvements on an earlier version of the section. Thanks also to Sergio Pavon for
pointing us to the beautiful recent article \cite{Hrbek-Pavon23},
which treats related topics from a completely
different perspective. We are also grateful to Jeremy Rickard for patiently answering our questions and to Evgeny Shinder for pointing out possible future developments about the relations between the autoequivalence groups of the intrinsic subcategories. We also thank Rudradip Biswas and Kabeer Menali  for suggestions and corrections and improvements to earlier versions.  

The second author did much of
his work during the trimester program Spectral Methods in
Algebra, Geometry, and Topology at the Hausdorff Institute in Bonn. It is a pleasure to thank the Hausdorff Institute for its hospitality, and acknowledge the funding by the Deutsche Forschungsgemeinschaft
(DFG) under Excellence Strategy EXC-2047/1-390685813.
A later part of the work was carried out while the second author was Mercator Fellow at the University of Bielefeld and the first and third author were visiting the same institution. We are pleased to thank the institution for its warm hospitality and financial support,
and also the Deutsche Forschungsgemeinschaft grant (SFB-TRR 358/1
2023 - 491392403) that helped underwrite the activity.}

\section{Reminder of some basic notation}
\label{S314}

In this short section we briefly summarize some basic constructions involving various notions of generation.

\subsection{Notation and constructions}\label{subsec:constr}

We start by recalling that, if $\ct$ is a triangulated category, then $\sh\colon\ct\la\ct$ denotes the shift functor.

\ntn{N314.1}
Let $\ct$ be a triangulated category, and let $\ca,\cb$ be full subcategories of $\ct$. The following conventions will be used throughout:
\be
\item
$\ca*\cb\subseteq\ct$ is the full subcategory of all objects $x\in\ct$ for which there exists a distinguished triangle $a\la x\la b$ with $a\in\ca$ and $b\in\cb$.
\item
$\add(\ca)\subseteq\ct$ is the full subcategory whose objects are all finite direct sums of objects of $\ca$.
\item
Assume $\ct$ has coproducts. Then $\Add(\ca)\subseteq\ct$ is the full subcategory whose objects are all the small coproducts of objects of $\ca$. 
\item
$\Smr(\ca)\subseteq\ct$ is the full
subcategory with objects all direct summands of objects of $\ca$.
\setcounter{enumiv}{\value{enumi}}
\ee
With this notation, we will furthermore adopt the conventions:
\be
\setcounter{enumi}{\value{enumiv}}
\item
$\COprod(\ca)$ is the
smallest full subcategory $\cs\subseteq\ct$ satisfying
\[
\ca\subseteq\cs,\qquad\cs*\cs\subseteq\cs,\qquad\add(\cs)\subseteq\cs.
\]
\item
$\gen\ca{}$ is given by the formula  
\[
\gen\ca{}:=\Smr\left(\COprod(\ca)\right).
\]
\item
Assume $\ct$ has coproducts. Then $\Coprod(\ca)$ is the
smallest full subcategory $\cs\subseteq\ct$ satisfying
\[
\ca\subseteq\cs,\qquad\cs*\cs\subseteq\cs,\qquad\Add(\cs)\subseteq\cs.
\]
\item
Still assuming $\ct$ has coproducts, then $\ogen\ca{}$ is given by
the formula  
\[
\ogen\ca{}:=\Smr\left(\Coprod(\ca)\right).
\]
\ee
\entn

The following special case will interest us later in the article.

\ntn{N314.3}
Let $\ct$ be a triangulated category, and let $G\in\ct$ be an object.
\be
\item
If $m\leq n$ are integers, possibly infinite, then $G[m,n]$ is defined
to be the full subcategory of $\ct$ whose objects are given by
\[
G[m,n]:=\{\sh[i]{G}\st i\in\zz\text{ and } m\leq -i\leq n\}.
\]
In the rest of the paper when $m=-\infty$ (resp.\ $n=+\infty$) we use the notation $G(-\infty,n]$ (resp.\ $G[m,+\infty)$) instead of the above one.
\item
The subcategory $\genu G{}{m,n}$ is defined by the formula
\[
\genu G{}{m,n}:=
\Smr\Big(\COprod\big(G[m,n]\big)\Big)\ .
\]
\item
Assuming $\ct$ has coproducts, the subcategory $\ogenu G{}{m,n}$ is defined by the formula
\[
\ogenu G{}{m,n}:=
\Smr\Big(\Coprod\big(G[m,n]\big)\Big)\ .
\]
\ee
\entn

The elementary properties of these constructions are discussed in
\cite[Section~1]{Neeman17}. We will freely use the results 
proved there when needed---the reader may wish to have a
quick look before proceeding with the rest of the
current article.

\subsection{Basic properties}\label{subsec:properties}

We end the section with reminders of facts well-known to the experts.
Since all of these facts are easy enough, we include the (short) proofs
for the reader's convenience.

Recall that a subcategory $\cc$ of
a category $\cd$ is \emph{strictly full} if it is full and,
given an object
$x\in\cc$ and an isomorphism $x\la y$ in $\cd$, then $y$ is in $\cc$.

As a matter of notation, if $\ct$ is a triangulated category and $\ca$ is a full subcategory of $\ct$, we will denote by $\ca^\perp$ (resp.\ ${}^\perp\ca$) the full subcategory of $\ct$ consisting of all objects $x$ such that $\Hom(a,x)=0$ (resp.\ $\Hom(x,a)=0$) for every $a\in\ca$.

\lem{R701.1}
If $\ct$ is
a triangulated category with coproducts, and if $\ca\subseteq\ct$ is any
full subcategory, then
\be
\item
The subcategory $\Coprod(\ca)\subseteq\ct$ is strictly full.
\item
The subcategory $\Coprod(\ca)^\perp\subseteq\ct$ is also
strictly full, and moreover $\ca^\perp=\Coprod(\ca)^\perp$.
\ee
\elem

\prf
We start with (i). By \autoref{N314.1}(vii) the full subcategory
$\cs=\Coprod(\ca)$ is closed in $\ct$ under all small coproducts, and coproducts are only defined up to isomorphism.
Thus any isomorph $B$
of an object $A\in\cs$ is a coproduct of objects
of $\cs$, in this case $B$ is the
coproduct of all objects in the (singleton) set $\{A\}$.

Now for (ii), the fact that $\cb^\perp$ is strictly full, for any full subcategory $\cb$ of $\ct$, is obvious by definition. Furthermore, \autoref{N314.1}(vii)
guarantees that $\ca\subseteq\Coprod(\ca)$, and hence $\Coprod(\ca)^\perp\subseteq\ca^\perp$. We need to prove the reverse inclusion.

Clearly $\Hom(\ca,\ca^\perp)=0$, giving the inclusion
$\ca\subseteq{^\perp\big(\ca^\perp\big)}$. But for any full subcategory
$\cb$ of $\ct$, the subcategory $\cs={^\perp\cb}$ is closed
under coproducts and satisfies $\cs*\cs\subseteq\cs$. Applying
this to $\cb=\ca^\perp$ we obtain, with
$\cs={^\perp\cb}={^\perp\big(\ca^\perp\big)}$,
\[
\ca\subseteq\cs,\qquad\cs*\cs\subseteq\cs,\qquad\Coprod(\cs)\subseteq\cs.
\]
Hence ${^\perp\big(\ca^\perp\big)}$ contains the minimal subcategory
satisfying the conditions, that is $\Coprod(\ca)\subset{^\perp\big(\ca^\perp\big)}$.
And this inclusion rewrites as $\Hom\big(\Coprod(\ca),\ca^\perp\big)=0$, that is $\ca^\perp\subseteq\Coprod(\ca)^\perp$.
\eprf

\lem{R701.3}
Let $\ct$ be a triangulated category with coproducts, let $\ca\subseteq\ct$
be a full subcategory, and assume $\sh\ca\subseteq\ca$. Then
\be
\item
$\sh\Coprod(\ca)\subseteq\Coprod(\ca)$.
\item
$\sh[-1]\Coprod(\ca)^\perp\subseteq\Coprod(\ca)^\perp$.
\item
The subcategory $\Coprod(\ca)$ contains the homotopy colimit of any  
countable sequence of its morphisms.
\ee
\elem

\prf
To see (i), observe that $\sh\Coprod(\ca)=\Coprod(\sh\ca)\subseteq\Coprod(\ca)$,
where the equality is because $\sh$ is an autoequivalence of the category
$\ct$ preserving extensions and (of course) coproducts, and the
containment is by applying $\Coprod(-)$ to the inclusion $\sh\ca\subseteq\ca$.

We deduce (ii) from (i) formally, in the following easy way. By combining (i) with the fact that $\Hom\big(\Coprod(\ca),\Coprod(\ca)^\perp\big)=0$, we get 
\[
\Hom\big(\sh\Coprod(\ca),\Coprod(\ca)^\perp\big)=\Hom\big(\Coprod(\ca),\sh[-1]\Coprod(\ca)^\perp\big)=0.
\]
This proves the inclusion in (ii).

For (iii) recall that the homotopy colimit of a sequence of composable morphisms
in $\ct$
\[\xymatrix{
X_1\ar[r]^-{f_1} & 
X_2\ar[r]^-{f_2} & 
X_3\ar[r]^-{f_3} & 
X_4\ar[r]^-{f_4} & 
\cdots
}\]
is defined (up to noncanonical isomorphism) by the distinguished triangle
\[\xymatrix@C+30pt{
\ds\coprod_{n=1}^\infty X_n \ar[r]^-{1-\shi} &
\ds\coprod_{n=1}^\infty X_n \ar[r]&
\hoco X_n\ar[r] &
\ds\sh\left(\coprod_{n=1}^\infty X_i\right)
}\]
where $1\colon\coprod_{n=1}^\infty X_n\la\coprod_{n=1}^\infty X_n$ is the identity
map and $\shi\colon\coprod_{n=1}^\infty X_n\la\coprod_{n=1}^\infty X_n$ is the
map whose only nonzero components are $f_n\colon X_n\la X_{n+1}$.

Since $X_n$ is assumed to belong to $\Coprod(\ca)$ for all $n>0$,
we deduce that $\coprod_{n=1}^\infty X_n$ belongs to
$\Coprod(\ca)$. By (i) we have that $\sh\Coprod(\ca)\subseteq\Coprod(\ca)$,
and therefore $\sh\big(\coprod_{n=1}^\infty X_n\big)$ belongs to
$\Coprod(\ca)$. And now the fact that $\hoco X_n$ belongs to
$\Coprod(\ca)$ is from its triangle of definition
and because $\Coprod(\ca)*\Coprod(\ca)\subseteq\Coprod(\ca)$.
\eprf

\section{The basics of \tstr s}
\label{betterlabelthis}

This section introduces the notion of \tstr\ and some of its basic properties. The presentation in \autoref{S699} is somewhat non-traditional and we devote \autoref{S700} to a comparison with the `standard' presentation of the subject going back to Be{\u\i}linson, Bernstein
and Deligne \cite{BeiBerDel82}. Finally, in \autoref{S701} we introduce the main example of interest in this paper: compactly generated \tstr s.

\subsection{Definitions and basic properties}\label{S699}

We begin with the following.

\dfn{D699.1}
Let $\ct$ be a triangulated category. A strictly full subcategory
$\cs\subseteq\ct$ is called a \emph{pre-aisle} if
\[
\sh\cs\subseteq\cs\qquad\text{ and }\qquad\cs*\cs\subseteq\cs\ .
\]
The pre-aisle $\cs\subseteq\ct$ is called an \emph{aisle} if the
inclusion $\fI\colon\cs\la\ct$ has a right adjoint $\fI_\rho\colon\ct\la\cs$.
\edfn

\rmk{R699.2}
Let $\ct$ be a triangulated category. A strictly full subcategory
$\cs\subseteq\ct$ is a \emph{co-pre-aisle} if $\cs\op\subseteq\ct\op$ is a
pre-aisle.
And $\cs\subseteq\ct$ is a \emph{co-aisle} if $\cs\op\subseteq\ct\op$ is an aisle.
\ermk

\lem{L699.3}
Let $\ct$ be a triangulated category, let $\cs\subseteq\ct$ be an aisle,
let $\fI\colon\cs\la\ct$ be the inclusion, let $\fI_\rho\colon\ct\la\cs$ be
its right adjoint and let $\e\colon\fI\circ\fI_\rho\la\id$ be the counit of the adjunction.

Let $B\in\ct$ be any object and complete the morphism
$\e_B\colon\fI\circ\fI_\rho(B)\la B$ to a distinguished triangle
$\fI\circ\fI_\rho(B)\stackrel{\e_B}{\la} B\stackrel g\la C\stackrel h\la \sh\circ\fI\circ\fI_\rho(B)$. Then $C$ must belong to the subcategory $\cs^\perp\subseteq\ct$.
\elem

\prf
Consider the morphism $\e_C\colon\fI\circ\fI_\rho(C)\la C$. It allows us first to build up a diagram
\[\xymatrix@C+30pt{
    &    &    \fI\circ\fI_\rho(C)\ar[d]_{\e_C} \ar[r]^-{h\circ\e_C} & \sh\circ\fI\circ\fI_\rho(B)\ar@{=}[d] \\
\fI\circ\fI_\rho(B)\ar[r]^-{\e_B} & B\ar[r]^-g & C\ar[r]^-h &  \sh\circ\fI\circ\fI_\rho(B)
}\]
where the square is commutative, and then
extend to a morphism of distinguished triangles
\begin{equation}\label{eq:commdist}
\xymatrix@C+30pt{
\fI\circ\fI_\rho(B)\ar[r]^-\alpha\ar@{=}[d] & E\ar[d]_\beta \ar[r]^-\gamma
  &   \fI\circ\fI_\rho(C)\ar[d]_-{\e_C} \ar[r]^-{h\circ\e_C} & \sh\circ\fI\circ\fI_\rho(B)\ar@{=}[d] \\
\fI\circ\fI_\rho(B)\ar[r]^-{\e_B} & B\ar[r]^-g & C\ar[r]^-h &  \sh\circ\fI\circ\fI_\rho(B).
}
\end{equation}
From the triangle on the top row, coupled with the facts
that both $\fI\circ\fI_\rho(B)$ and $\fI\circ\fI_\rho(C)$ belong to
$\cs$ and that $\cs*\cs\subseteq\cs$, we
learn that $E$ must belong to $\cs$.
Therefore the morphism $\beta\colon E\la B$
must factor (uniquely) as $E\stackrel\s\la \fI\circ\fI_\rho(B)\stackrel{\e_B}{\la} B$.
But then the equality $\beta\circ\alpha=\e_B$ becomes $\e_B\circ\s\circ\alpha=\e_B$, that
is the two composites
\[\xymatrix@C+30pt{
\fI\circ\fI_\rho(B)\ar@<0.5ex>[r]^-{\s\circ\alpha}
\ar@<-0.5ex>[r]_-{\id} & 
\fI\circ\fI_\rho(B)\ar[r]^-{\e_B} & B
}\]
must be equal. But any morphism $S\la B$, where $S\in\cs$, has a \emph{unique}
factorization as $S\la \fI\circ\fI_\rho(B)\stackrel{\e_B}{\la}B$, and applying this uniqueness to $S=\fI\circ\fI_\rho(B)$, we deduce that $\s\circ\alpha=\id$. That is the map $\alpha\colon\fI\circ\fI_\rho(B)\la E$ must be
a split monomorphism, and then from the fact that the top row in \eqref{eq:commdist} is a distinguished triangle we deduce $h\circ\e_C=0$. Since also the bottom row in \eqref{eq:commdist} is a distinguished triangle, $\e_C$ must factor (in some way) as $\fI\circ\fI_\rho(C)\stackrel f\la B\stackrel g\la C$. But as $\fI\circ\fI_\rho(C)$ belongs to
$\cs$ the map $f\colon\fI\circ\fI_\rho(C)\la B$ must factor (uniquely) as
$\fI\circ\fI_\rho(C)\stackrel\tau\la \fI\circ\fI_\rho(B)\stackrel{\e_B}\la B$. Combining
these factorizations we have that $\e_C\colon\fI\circ\fI_\rho(C)\la C$ factors as
\[\xymatrix@C+30pt{
\fI\circ\fI_\rho(C)\ar[r]^\tau & 
\fI\circ\fI_\rho(B)\ar[r]^-{\e_B} & B\ar[r]^-g & C\ ,
}\]
and, as $g\circ\e_B=0$, this composite must vanish.

Finally, every morphism $S\la C$, where $S\in\cs$ is any object, factors
(uniquely) as $S\la \fI\circ\fI_\rho(C)\stackrel{\e_C}\la C$. And now that we know the
vanishing of $\e_C\colon\fI\circ\fI_\rho(C)\la C$, it follows that any map $S\la C$ must
vanish. That is $C$ belongs to the subcategory $\cs^\perp\subseteq\ct$.
\eprf

And now we give the older, more symmetric version of the same.

\dfn{D699.5}
Let $\ct$ be a triangulated category. A \emph{\tstr}\ on $\ct$ is a pair
of strictly full subcategories $(\cs,\cs')$ satisfying:
\be
\item
$\sh\cs\subseteq\cs$ and $\sh[-1]\cs'\subseteq\cs'$.
\item
$\Hom(\cs,\cs')=0$.
\item
For any object $B\in\ct$ there exists a distinguished triangle $A\la B\la C\la\sh{A}$, 
with $A\in\cs$ and $C\in\cs'$.  
\ee
\edfn

\rmk{R699.6}
We note that the definition is self-dual: the pair $(\cs,\cs')$ is
a \tstr\ on $\ct$ if and only if the pair $\big((\cs')\op,\cs\op\big)$
is a \tstr\ on $\ct\op$.
\ermk

\exm{E699.7}
If $\ct$ is a triangulated category and $\cs\subseteq\ct$ is an aisle, then
the pair $(\cs,\cs^\perp)$ is a \tstr\ on $\ct$.

Indeed, because $\cs$ is an aisle, we are given that $\cs\subseteq\ct$ is a strictly full subcategory. The fact that
$\cs^\perp\subseteq\ct$ is strictly full is by its definition.
Because $\cs\subseteq\ct$ is an aisle, we are also given the inclusion
$\sh\cs\subseteq\cs$. The other inclusion in \autoref{D699.5}(i) follows from the same formal argument as in the proof of \autoref{R701.3}(ii).
The fact that $\Hom(\cs,\cs^\perp)=0$ is obvious, giving that the pair
$(\cs,\cs^\perp)$ satisfies \autoref{D699.5}(ii).
The fact that the pair
$(\cs,\cs^\perp)$ satisfies \autoref{D699.5}(iii) was proved
in \autoref{L699.3}.
\eexm

The next proposition shows, among
other things, that all \tstr{s} are as in
\autoref{E699.7}. This justifies the assertion we made, just
before \autoref{D699.5}, that \tstr{s} are just a more
symmetric presentation of the information contained in giving
an aisle.

\pro{P699.9}
Let $\ct$ be a triangulated category, and let $(\cs,\cs')$ be a \tstr\
on $\ct$. Then the following assertions are all true:
\be
\item
$\cs\subseteq\ct$ is an aisle and $\cs'\subseteq\ct$ is a co-aisle.
\item
$\cs'=\cs^\perp$ and $\cs={^\perp\cs'}$.
\setcounter{enumiv}{\value{enumi}}
\ee

Furthermore, let $\fI\colon\cs\la\ct$ and $\fJ\colon\cs'\la\ct$ be the inclusion functors, let
$\fI_\rho\colon\ct\la\cs$ be the right adjoint of $\fI$, let $\fJ_\lambda\colon\ct\la\cs'$
be the left adjoint of $\fJ$, let $\e\colon\fI\circ\fI_\rho\la\id$ be the counit of
the adjunction $\fI\dashv\fI_\rho$, and let $\eta\colon\id\la\fJ\circ\fJ_\lambda$ be the
unit of the adjunction $\fJ_\lambda\dashv\fJ$. Then we have:
\be
\setcounter{enumi}{\value{enumiv}}
\item
For each object $B\in\ct$, there exists a unique morphism
$\ph_B\colon\fJ\circ\fJ_\lambda(B)\la\sh\circ\fI\circ\fI_\rho(B)$ such that
\[\xymatrix@C+30pt{
\fI\circ\fI_\rho(B)\ar[r]^-{\e_B}
& B
\ar[r]^-{\eta_B} & \fJ\circ\fJ_\lambda(B) \ar[r]^-{\ph_B} &
\sh\circ\fI\circ\fI_\rho(B)
}\]
is a distinguished triangle. 
\item
The assignment taking an object $B\in\ct$, to the unique
morphism $\ph_B$ satisfying the condition in \emph{(iii)}, delivers 
a natural transformation $\ph\colon\fJ\circ\fJ_\lambda\la\sh\circ\fI\circ\fI_\rho$.
\ee
\epro

\prf
Let us begin by proving (ii), and note that it suffices to prove one of the
equalities in (ii) as the other is its dual. Hence, we are reduced
to proving
the equality $^\perp\cs'=\cs$.

In \autoref{D699.5}(ii) we are given that $\Hom(\cs,\cs')=0$,
that is the inclusion $\cs\subseteq{^\perp\cs'}$ is given. What needs proof is
the reverse inclusion---we need to show that $^\perp\cs'\subseteq\cs$.

Choose therefore any object $B\in{^\perp\cs'}$. \autoref{D699.5}(iii) gives the existence of a distinguished triangle
$A\overset{\alpha}{\la} B\la C\la\sh A$ with $A\in\cs$ and $C\in\cs'$. Note that
\[
\sh A\in\sh\cs\subseteq\cs\subseteq{^\perp\cs'},
\]
where the first inclusion is by \autoref{D699.5}(i) and
the second by \autoref{D699.5}(ii). The part $B\la C\la\sh A$
of the triangle, coupled with the fact that we now know that both
$B$ and $\sh A$ belong to $^\perp\cs'$, tells us that $C$ must belong
to $({^\perp\cs'})*({^\perp\cs'})\subset{^\perp\cs'}$. Thus $C$ must
belong to ${^\perp\cs'}\cap\cs'$, and hence must be isomorphic to zero.
Thus the map $\alpha$ must be an isomorphism. Since $A\in\cs$
and the subcategory $\cs\subseteq\ct$ is strictly full, we have that
$B\in\cs$. This completes the proof of (ii).

Next we prove (i), and note that by duality it suffices to prove that
$\cs'\subseteq\ct$ is a co-aisle.

The inclusion $\sh[-1]\cs'\subseteq\cs'$ is part of \autoref{D699.5}(i).
By part (ii) of the current proposition we know that $\cs'=\cs^\perp$, and
the inclusion $\cs'*\cs'\subseteq\cs'$ follows immediately. Thus the
information we have so far easily shows that $\cs'\subseteq\ct$ is
a co-pre-aisle. What needs proof is that the inclusion $\fJ\colon\cs'\la\ct$
has a left adjoint. We have to show that, for every object $B\in\ct$,
the functor $\Hom_\ct^{}\big(B,\fJ(-)\big)$ is a representable functor
$\cs'\la\ab$. Choose therefore an object $B\in\ct$.

\autoref{D699.5}(iii) allows us to find some distinguished triangle
$A\la B\stackrel\eta\la\fJ(C)\la\sh A$ with $A\in\cs$ and $C\in\cs'$.
Now choose any object $X\in\cs$, and apply $\Hom_\ct^{}\big(-,\fJ(X)\big)$ to this triangle. This gives us a long exact sequence, but as
$\Hom\big(-,\fJ(X)\big)$ kills $\cs\subseteq\ct$ and $A$ and $\sh A$ both belong
to $\cs$, this long exact sequence simplifies to the
second isomorphism in
\begin{equation}\label{eq:isos1}
\xymatrix@C+30pt{
\Hom_{\cs'}^{}(C,X)
\ar[r]^-{\fJ} &
\Hom_\ct^{}\big(\fJ(C),\fJ(X)\big)
\ar[r]^-{\Hom\big(\eta,\fJ(X)\big)} &
\Hom_\ct^{}\big(B,\fJ(X)\big)
}
\end{equation}
The first isomorphism is by the full faithfulness of $\fJ$.

So far we have produced a natural isomorphism
$\Hom_{\cs'}^{}(C,-)\cong\Hom_\ct^{}\big(B,\fJ(-)\big)$. The functor
$\Hom_\ct^{}\big(B,\fJ(-)\big)$ is therefore representable for every
$B\in\ct$, allowing us to form the left adjoint $\fJ_\lambda\colon\ct\la\cs'$
to the natural inclusion. In fact, for $B\in\ct$ we can choose
$\fJ_\lambda(B)\in\cs'$ to be $\fJ_\lambda(B)=C$, and the extension of this to
a functor is unique. This completes the proof of (i).

But while we are at it let us also compute the counit of adjunction.
By definition, given any object $B\in\ct$, the counit of adjunction (evaluated
at the object $B$) is the image of
$1\in\Hom_{\cs'}^{}\big(\fJ_\lambda(B),\fJ_\lambda(B)\big)$
under the natural isomorphism
\[
\xymatrix@C+30pt{
\Hom_{\cs'}^{}\big(\fJ_\lambda(B),\fJ_\lambda(B)\big)
\ar[r]^-{\cong} &
\Hom_\ct^{}\big(B,\fJ\circ\fJ_\lambda(B)\big)
}
\]
By our choice we have that $\fJ_\lambda(B)=C$, and the isomorphism above
comes down to the composite in \eqref{eq:isos1} with $X=C$.
To put it concretely, given an object $B\in\ct$, we choose an
object $C\in\cs'$ and a distinguished triangle
$A\la B\stackrel\eta\la\fJ(C)\la\sh A$ with $A\in\cs$. The existence
of such a triangle comes from
\autoref{D699.5}(iii). And, when we choose $\fJ_\lambda(B)$ to be
$\fJ_\lambda(B)=C$ as above, the image of $1\in\Hom_{\cs'}^{}(C,C)$
under the natural composite computes to be $\eta\colon B\la\fJ(C)$.

Now we proceed with the proof of (iii) and (iv). We begin by reminding the
reader that the right adjoint $\fI_\rho$ of $\fI$ is only defined up to canonical
isomorphism, as is the left adjoint $\fJ_\lambda$ of $\fJ$. That is, given
an object $B\in\ct$, the objects $\fI_\rho(B)$ and $\fJ_\lambda(B)$ are defined
up to canonical isomorphism. If $\wt\fI_\rho(B)$ and $\wt\fJ_\lambda(B)$ are a
different pair of choices, then there are canonical isomorphisms
$\alpha\colon\fI_\rho(B)\la\wt\fI_\rho(B)$ and $\beta\colon\fJ_\lambda(B)\la\wt\fJ_\lambda(B)$,
rendering commutative the diagram
\[\xymatrix@C+30pt{
\fI\circ\fI_\rho(B)\ar[r]^-{\e_B}\ar[d]_{\fI(\alpha)}
& 
B\ar[r]^-{\eta_C}\ar@{=}[d] & \fJ\circ\fJ_\lambda (B) \ar[d]^-{\fJ(\beta)}\\
\fI\circ\wt\fI_\rho (B)\ar[r]^-{\wt\e_B}
& 
B\ar[r]^-{\wt\eta_C} &\fJ\circ\wt\fJ_\lambda (B)
}\]
where $\e,\wt\e$ are the counits and $\eta,\wt\eta$ the units of the
respective adjunctions. But the existence assertion in (iii)
is clearly independent of the choices.

Let $B$ be any object of $\ct$;
\autoref{D699.5}(iii) allows us to choose objects $A\in\cs$ and
$C\in\cs'$ and in $\ct$ a distinguished triangle
$\fI(A)\stackrel\e\la B\stackrel\eta\la\fJ(C)\stackrel\ph\la\sh\circ\fI(A)$.
We saw above that the left adjoint $\fJ_\lambda\colon\ct\la\cs'$ of
the inclusion $\fJ\colon\cs'\la\ct$ can be chosen so that $\fJ_\lambda(B)=C$,
and with this choice the unit of adjunction $\eta_B\colon B\la\fJ\circ\fJ_\lambda(B)$
becomes the map $\eta\colon B\la\fJ(C)$. Dually, the right adjoint
$\fI_\rho\colon\ct\la\cs$ to the inclusion $\fI\colon\cs\la\ct$ can be chosen so
that $\fI_\rho(B)=A$, and with that choice the counit of adjunction
$\e_B\colon\fI\circ\fI_\rho(B)\la B$ becomes the map $\e\colon\fI(A)\la B$. Thus the triangle
$\fI(A)\stackrel\e\la B\stackrel\eta\la\fJ(C)\stackrel\ph\la\sh\circ\fI(A)$
rewrites as
\[\xymatrix@C+30pt{
\fI\circ\fI_\rho(B)\ar[r]^-{\e_B}
& B
\ar[r]^-{\eta_B} & \fJ\circ\fJ_\lambda(B) \ar[r]^-{\ph_B} &
\sh\circ\fI\circ\fI_\rho(B)\ ,
}\]
proving the existence assertion of (iii).

It remains to prove the uniqueness in (iii), and to prove (iv).
Choose therefore any morphism $f\colon B\la\wt B$ in
$\ct$. By the existence part of (iii)  we can
choose morphisms
$\ph_B\colon\fJ\circ\fJ_\lambda (B) \la
\sh\circ\fI\circ\fI_\rho (B)$
and
$\ph_{\wt B}\colon\fJ\circ\fJ_\lambda (\wt B) \la
\sh\circ\fI\circ\fI_\rho (\wt B)$
so that in the diagram
\[\xymatrix@C+30pt{
\fI\circ\fI_\rho (B)\ar[r]^-{\e_B}\ar[d]_{\fI\circ\fI_\rho(f)}
& 
B\ar[r]^-{\eta_B}\ar[d]^f & \fJ\circ\fJ_\lambda (B) \ar[r]^-{\ph_B} &
\sh\circ\fI\circ\fI_\rho (B)\\
\fI\circ\fI_\rho (B)\ar[r]^-{\e_{\wt B}}
& 
\wt B\ar[r]^-{\eta_{\wt B}} & \fJ\circ\fJ_\lambda (\wt B) \ar[r]^-{\ph_{\wt B}} &
\sh\circ\fI\circ\fI_\rho (B)
}\]
the rows are distinguished triangles. The square commutes by the naturality of
$\e\colon\fI\circ\fI_\rho\la\id$. The axioms of triangulated categories allow
us to extend this to some morphism of triangles
\[\xymatrix@C+30pt{
\fI\circ\fI_\rho (B)\ar[r]^-{\e_B}\ar[d]_{\fI\circ\fI_\rho(f)}
& 
B\ar[r]^-{\eta_B}\ar[d]^f & \fJ\circ\fJ_\lambda (B) \ar[r]^-{\ph_B}\ar[d]^-{\pi} &
\sh\circ\fI\circ\fI_\rho (B)\ar[d]^-{\sh\circ\fI\circ\fI_\rho(f)}\\
\fI\circ\fI_\rho (B)\ar[r]^-{\e_{\wt B}}
& 
\wt B\ar[r]^-{\eta_{\wt B}} & \fJ\circ\fJ_\lambda (\wt B) \ar[r]^-{\ph_{\wt B}} &
\sh\circ\fI\circ\fI_\rho (B)
}\]
Next observe that both of the squares in the diagram below
\[\xymatrix@C+60pt{ 
B\ar[r]^-{\eta_B}\ar[d]^f & \fJ\circ\fJ_\lambda (B) \ar@<0.5em>[d]^{\fJ\circ\fJ_\lambda(f)}  \ar@<-0.5em>[d]_\pi
  \\
\wt B\ar[r]^-{\eta_{\wt B}} & \fJ\circ\fJ_\lambda (\wt B) 
}\]
commute; the one involving $\pi$ by the above, and the one
involving $\fJ\circ\fJ_\lambda(f)$ by the naturality of $\eta$. This gives two
factorizations of the composite
$B\stackrel f\la\wt B\stackrel{\eta_{\wt B}}\la\fJ\circ\fJ_\lambda(\wt B)$
through $\eta_B\colon B\la\fJ\circ\fJ_\lambda(B)$, which must agree. Thus $\pi=\fJ\circ\fJ_\lambda(f)$,
and we deduce the commutativity of the square
\[\xymatrix@C+30pt{
\fJ\circ\fJ_\lambda (B) \ar[r]^-{\ph_B}\ar[d]_{\fJ\circ\fJ_\lambda(f)} &
\sh\circ\fI\circ\fI_\rho (B)\ar[d]^{\sh\circ\fI\circ\fI_\rho(f)}\\
\fJ\circ\fJ_\lambda (\wt B) \ar[r]^-{\ph_{\wt B}} &
\sh\circ\fI\circ\fI_\rho (\wt B).
}\]

If we consider the special case where $f\colon B\la\wt B$ is the identity map
$\id\colon B\la B$, the commutativity
of the square above
tells us that any two choices $\ph_B,\ph_{\wt B}$ for the
map $\fJ\circ\fJ_\lambda(B)\la\sh\circ\fI\circ\fI_\rho(B)$ must agree---that is we have proved the uniqueness part of
(iii). And now that we know both the existence and
the uniqueness assertions in (iii),
the commutativity of the square
above, for any morphism $f\colon B\la\wt B$, completes the proof of (iv).
\eprf

\cor{C699.11}
Let $\ct$ be a triangulated category, and let $\cs\subseteq\ct$ be an aisle. Then  $\cs$ is closed in $\ct$ under direct summands and coproducts. Dually, any co-aisle $\cs'\subseteq\ct$ is closed in $\ct$ under direct summands and products.
\ecor

Just to clarify the statement above, let us stress that it says that any direct summand in $\ct$ of an object
of $\cs$ must belong to $\cs$, and given any collection of objects
in $\cs$ whose
coproduct exists in $\ct$, that coproduct must belong to $\cs$.

\prf
It clearly suffices to prove the assertion about aisles. Let $\cs\subseteq\ct$
be an aisle; by \autoref{E699.7} there exists a subcategory $\cs'\subseteq\ct$
such that the pair $(\cs,\cs')$ is a \tstr\ on $\ct$; more explicitly
\autoref{E699.7} tells us that $\cs'=\cs^\perp$ works. But what is important
for us here is \autoref{P699.9}(ii); it gives us the equality
$\cs={^\perp\cs'}$. The corollary now follows because $^\perp\cs'$ is
always closed in $\ct$ under direct summands and coproducts. 
\eprf

\subsection{The more traditional notation for \tstr s}
\label{S700}

The notion of \tstrs originated with
Be{\u\i}linson, Bernstein
and Deligne~\cite[Section~1.3]{BeiBerDel82}. The
original motivation came
from applications irrelevant to this manuscript,
and hence in~\cite[Section~1.3]{BeiBerDel82}
the subject is
pursued from a perspective different from our
\autoref{S699}---our
treatment is much
closer in spirit to Keller
and Vossieck~\cite{Keller-Vossieck88}. The very least we
owe the reader is a glossary, explaining how to pass
back and forth between our notation and the traditional
one.

\ntn{N700.1}
Let $\ct$ be a triangulated category, and let $(\cs,\cs')$
be a \tstr\ on $\ct$ as in \autoref{D699.5}.
One sets
\be
\item
$\ct^{\leq n}:=\sh[-n]\cs$.
\item
$\ct^{\geq n}:=\sh[-n+1]\cs'$.
\item
The triangle
$\fI\circ\fI_\rho(B)\stackrel{\e_B}\la B\stackrel{\eta_B}\la\fJ\circ\fJ_\lambda(B)\stackrel{\ph_B}\la\sh\circ\fI\circ\fI_\rho(B)$, of
\autoref{P699.9}, is traditionally written
\[\xymatrix{
B^{\leq0}\ar[r]^-{\e} & B\ar[r]^-{\eta} &B^{\geq1}\ar[r]^-{\ph} &
\sh B^{\leq0}\ .
}\]
In other words, the functor which we have written as
$\fI\circ\fI_\rho\colon\ct\la\ct$ is normally written $(-)^{\leq0}$,
and the functor which we have written as
$\fJ\circ\fJ_\lambda\colon\ct\la\ct$ is normally written $(-)^{\geq1}$.
\item
More generally, for any integer $n\in\zz$ we define the
functors $(-)^{\leq n}$ and $(-)^{\geq n}$ by the formulas
\[
(-)^{\leq n}:=\sh[-n]\circ\fI\circ\fI_\rho\circ\sh[n],\qquad
(-)^{\geq n}:=\sh[-n+1]\circ\fJ\circ\fJ_\lambda\circ\sh[n-1].
\]
\item
As a pure matter of notation, in \autoref{D699.5}
a \tstr\ was defined to be a pair of subcategories
$(\cs,\cs')$ satisfying some conditions. In the
comparison with the standard notation, we let
$\cs=\ct^{\leq0}$ and $\cs'=\ct^{\geq1}$. In reading
the literature the reader should just note that
it is traditional to give a \tstrs as a pair of
subcategories $(\ct^{\leq0},\ct^{\geq0})$. To be consistent with \autoref{D699.5}, in the rest of the paper, a \tstrs will be denoted $\tst\ct$. 
\ee
\entn

As the reader can easily check, with the definitions as
above, we have, for any integer $n\in\zz$ and any object
$B\in\ct$, a distinguished triangle
\[\xymatrix{
B^{\leq n}\ar[r]^-{\e_n} & B\ar[r]^-{\eta_n} &B^{\geq n+1}\ar[r]^-{\ph_n} &
\sh B^{\leq n}\ ,
}\]
with $B^{\leq n}\in\ct^{\leq n}\subseteq\ct$ and with
$B^{\geq n+1}\in\ct^{\geq n+1}\subseteq\ct$, and this triangle
is functorial in $B$. This is just the appropriate
translation of the triangle of (iii).

Next, we recall the following definition.
The reader can find a more complete treatment in
Be{\u\i}linson, Bernstein
and Deligne~\cite[Section~1.2]{BeiBerDel82}.

\dfn{D700.2}
Let $\ct$ be a triangulated category. A strictly full subcategory $\ca\subseteq\ct$ is called an \emph{admissible abelian subcategory of $\ct$} if
the following conditions hold:
\be
\item
$\Hom(\ca,\sh[n]\ca)=0$ for all $n<0$.
\item
The category $\ca$ is abelian, and the short exact sequences 
in $\ca$ are precisely the pairs of morphisms
$A\stackrel f\la B\stackrel g\la C$ for which there exists
a morphism $h\colon C\la\sh A$ rendering
$A\stackrel f\la B\stackrel g\la C\stackrel h\la\sh A$ a
distinguished triangle in $\ct$.
\ee
\edfn

Next we want to refer the reader to
known facts we will use, without providing
our own proofs. In the case of
these particular known facts, the results and proofs
are available in a short section
of a single manuscript, and we deemed that giving
a reference was the way to go.

\rmd{R700.3}
Let $\ct$ be a triangulated category, and let
$\tst\ct$ be a \tstr\ on $\ct$, where the
  notation is as in \autoref{N700.1}(v). Then the following is true:
\be
\item
Let $m\leq n$ be integers. Then the functors 
$(-)^{\leq n}$ and $(-)^{\geq m}$ commute with each other
up to natural isomorphism. Hence the composite functor
\[
\Big((-)^{\leq n}\Big)^{\geq m}\cong 
\Big((-)^{\geq m}\Big)^{\leq n}
\]
can be viewed as a functor 
$\tau^{[m,n]}\colon\ct\la\ct^{\geq m}\cap\ct^{\leq n}$.
\item
The category $\ct^\heartsuit:=\ct^{\geq0}\cap\ct^{\leq0}$
is an admissible abelian subcategory of $\ct$ called the \emph{heart} of
the \tstr\ $\tst\ct$, and the functor $\tau^{[0,0]}\colon\ct\la\ct^\heartsuit$
is homological. It will be often denoted $\ch^0\colon\ct\la\tsth\ct{}$. As usual, we write $\ch^n=\ch^{0}\circ\sh[n]$, for every integer $n$.
\ee
We refer to \cite[Proposition~1.3.5]{BeiBerDel82} for (i) and to \cite[Th\'eor\`eme~1.3.6]{BeiBerDel82} for (ii).
\ermd

\subsection{Compactly generated \tstr s}
\label{S701}

As in \cite{Rouquier08}, we introduce the following:

\dfn{D701.5}
Let $\ct$ be a triangulated category. An
object $c\in\ct$ is declared to be 
\emph{compact} if $\Hom(c,-)$ respects those
coproducts that exist in $\ct$.
\edfn

When $\ct$ has coproducts this goes back to
\cite[Definition~0.1]{Neeman92A}, and there are
powerful techniques allowing us to work out what the
compact objects are. But we will care about the more
general situation.

\rmk{R701.7}
Let $\ct$ be a triangulated category. Then the full
subcateory $\ct^c\subseteq\ct$ is defined to have for
objects all the compacts in $\ct$.
Clearly $\ct^c$ is strictly full, and one easily checks that
\be
\item
$\ct^c*\ct^c\subseteq\ct^c$.
\item
$\sh\ct^c=\ct^c$.
\item
$\Smr(\ct^c)=\ct^c$
\setcounter{enumiv}{\value{enumi}}  
\ee
Or, to summarize the above properties more densely,
$\ct^c$ is a thick, triangulated subcategory of $\ct$.
But we also remind the reader of
\be
\setcounter{enumi}{\value{enumiv}}  
\item
For any compact object $c\in\ct$ and any
sequence of objects and composable
morphisms
in $\ct$
\[\xymatrix{
X_1\ar[r]^-{f_1} & 
X_2\ar[r]^-{f_2} & 
X_3\ar[r]^-{f_3} & 
X_4\ar[r]^-{f_4} & 
\cdots
}\]
the natural map
\[\xymatrix{
\colim\,\Hom(c,X_n)\ar[r] &\Hom(c,\hoco X_n)
}\]
is an isomorphism.
\ee

As we already said: (i), (ii) and (iii) are easy and left
to the reader. To prove (iv) apply the homological
functor $\Hom(c,-)$ to the triangle defining
$\hoco X_n$ (see the proof of \autoref{R701.3}(iii) for
a reminder of this triangle), and then use the fact that the compactness
of $c$ allows us to compute what $\Hom(c,-)$ does
to the coproduct terms in this triangle.
\ermk

The next result is the combination of
\cite[Theorem~A.1 and Proposition~A.2]{Alonso-Jeremias-Souto03}.

\thm{T701.9}
Let $\ct$ be a triangulated category with coproducts,
let $\ca\subseteq\ct^c$ be an essentially small full subcategory,
and assume $\sh\ca\subseteq\ca$. Then the pair
$\tau_\ca:=\big(\Coprod(\ca),\Coprod(\ca)^\perp\big)$ is
a \tstr\ on $\ct$, and the subcategory
$\Coprod(\ca)^\perp$ is closed in $\ct$ under coproducts.
\ethm

\prf
The fact that $\Coprod(\ca)$ and $\Coprod(\ca)^\perp$
are strictly full was proved in \autoref{R701.1},
the fact that $\sh\Coprod(\ca)\subseteq\Coprod(\ca)$
and $\sh[-1]\Coprod(\ca)^\perp\subseteq\Coprod(\ca)^\perp$
was proved in \autoref{R701.3}, and the
fact that $\Hom\big(\Coprod(\ca),\Coprod(\ca)^\perp\big)=0$
is obvious.

To show that $\tau_\ca$ is a \tstr\ on $\ct$, 
it remains to check \autoref{D699.5}(iii):
for every $X\in\ct$ we must produce a distinguished triangle
$A\la X\la C\la\sh A$, with $A\in\Coprod(\ca)$ and
$C\in\Coprod(\ca)^\perp$. We propose to do it as follows:
starting with $X\in\ct$ we will produce a
sequence
\[\xymatrix{
X_1\ar[r]^-{f_1} & 
X_2\ar[r]^-{f_2} & 
X_3\ar[r]^-{f_3} & 
X_4\ar[r]^-{f_4} & 
\cdots
}\]
of objects and morphisms in $\Coprod(\ca)$, all
mapping to $X$, and we will then let $A\la X$ be any
compatible choice of the
map $\hoco X_n\la X$.

Let us start. Because the subcategory $\ca$ is essentially
small, there is (up to isomorphism) only a set $\Lambda$ of
morphisms $f_\lambda\colon A_\lambda\la X$ with $A_\lambda\in\ca$.
We begin with
\be
\item
Let $X_1=\coprod_{\lambda\in\Lambda}A_\lambda$, and
let the morphism $h_1\colon X_1\la X$ be the obvious map.
We observe that, for every object $A\in\ca$, the functor
$\Hom(A,-)$ takes the map $h_1$ to an epimorphism.
\setcounter{enumiv}{\value{enumi}}  
\ee
Now suppose we are given a morphism $h_n\colon X_n\la X$, with
$X_n\in\Coprod(\ca)$. Complete this to a distinguished triangle
$Y_n\stackrel{\alpha_n}\la X_n\stackrel{h_n}\la X$, and let $\Lambda_n$ be
a set containing isomorphism classes of all
possible maps $k\colon B\la Y_n$ with
$B\in\sh[-1]\ca$. Next we define
\be
\setcounter{enumi}{\value{enumiv}}  
\item
The object $Z_n\in\sh[-1]\Coprod(\ca)$ is given
by the formula $Z_n:=\coprod_{\lambda\in\Lambda_n}B_\lambda$,
and the map $\beta_n\colon Z_n\la Y_n$ is the obvious. By construction
$\Hom(B,-)$ takes the map $\beta_n\colon Z_n\la Y_n$ to an epimorphism
whenever $B\in\sh[-1]\ca$.
\setcounter{enumiv}{\value{enumi}}  
\ee
Now we have composable morphisms
$Z_n\stackrel{\beta_n}\la Y_n\stackrel{\alpha_n}\la X_n$; then we can complete
the morphism $\gamma_n=\alpha_n\circ\beta_n\colon Z_n\la X_n$ to a distinguished triangle
$Z_n\stackrel{\gamma_n}\la X_n\stackrel{f_n}\la X_{n+1}\la \sh Z_n$.
As both $X_n$ and $\sh Z_n$ lie in
$\Coprod(\ca)$, the triangle gives that
$X_{n+1}\in\Coprod(\ca)$. And now the
composite
\[\xymatrix@C+20pt{
Z_n\ar[r]_{\beta_n}\ar@/^2pc/[rr]^{\gamma_n} & Y_n\ar[r]_{\alpha_n}& X_n\ar[r]_{h_n}& X
}\]
vanishes because $h_n\circ\alpha_n=0$, allowing us to factor the map
$h_n\colon X_n\la X$ as
$X_n\stackrel{f_n}\la X_{n+1}\stackrel{h_{n+1}}\la X$.
And the notable features of the construction are
\be
\setcounter{enumi}{\value{enumiv}}  
\item
We have factored the map $h_n\colon X_n\la X$ as
$X_n\stackrel{f_n}\la X_{n+1}\stackrel{h_{n+1}}\la X$,
with $X_{n+1}\in\Coprod(\ca)$, and in such a way
that, for every $B\in\sh[-1]\ca$, the kernel
of $\Hom(B,h_n)\colon\Hom(B,X_n)\la\Hom(B,X)$ is
annihilated by the morphism
$\Hom(B,f_n)\colon\Hom(B,X_n)\la\Hom(B,X_{n+1})$.
\setcounter{enumiv}{\value{enumi}}  
\ee
To put it more succinctly: for every
$B\in\sh[-1]\ca$, the sequence
\[\xymatrix{
\Hom(B,X_1)\ar[r]^-{f_1} & 
\Hom(B,X_2)\ar[r]^-{f_2} & 
\Hom(B,X_3)\ar[r]^-{f_3} & 
\Hom(B,X_4)\ar[r]^-{f_4} & 
\cdots
}\]
is Ind-isomorphic to the sequence of monomorphisms
\[\xymatrix{
I_1\ar[r]^-{\wt f_1} & 
I_2\ar[r]^-{\wt f_2} & 
I_3\ar[r]^-{\wt f_3} & 
I_4\ar[r]^-{\wt f_4} & 
\cdots
}\]
where $I_n$ is the image of the map
$\Hom(B,h_n)\colon\Hom(B,X_n)\la\Hom(B,X)$. If $B$ belongs
to $\ca\subseteq\sh[-1]\ca$, this is actually a sequence of
isomorphisms; after all by (i) the map $I_1\la\Hom(B,X)$ is
epi, and hence so are all the maps $I_n\la\Hom(B,X)$.
But being both epi and mono they must all be isomorphisms.
Now let $A:=\hoco X_n$. By
\autoref{R701.3}(iii) the object $A$ belongs
to $\Coprod(\ca)$ and, with any choice of
the morphism $h\colon A\la X$ compatible with all the $h_n$,
the computation above (taking into account also \autoref{R701.7}(iv)) tells us that $\Hom(B,h)$ is a
monomorphism for all $B\in\sh[-1]\ca$ and an isomorphism
if $B\in\ca$. Now complete to a distinguished triangle
$A\stackrel h\la X\la C\la\sh A$. Take any object
$G\in\ca$ and apply the homological
functor $\Hom(G,-)$ to this triangle. We deduce
an exact sequence
\[\xymatrix{
\Hom(G,A)\ar[r]^-{\alpha} & 
\Hom(G,X)\ar[r] & 
\Hom(G,C)\ar[r] & 
\Hom(G,\sh A)\ar[r]^-{\beta} & \Hom(G,\sh X) 
}\]
where, by the above, $\alpha$ is an isomorphism
and $\beta$ is a monomorphism. Hence
$\Hom(G,C)=0$, and as this is true for all
$G\in\ca$ we deduce that $C\in\ca^\perp$, which is
equal to $\Coprod(\ca)^\perp$ 
by \autoref{R701.1}(ii). This completes
the proof that $\tau_\ca$ is a \tstr\ on $\ct$.

And now the equality $\Coprod(\ca)^\perp=\ca^\perp$
of \autoref{R701.1}(ii), coupled with the fact
that $\ca\subseteq\ct^c$, immediately
tells us that $\Coprod(\ca)^\perp$
is closed in $\ct$ under coproducts.
\eprf

\rmk{R701.11}
As we have already said, \autoref{T701.9} is
the union of 
\cite[Theorem~A.1 and Proposition~A.2]{Alonso-Jeremias-Souto03}.
The reason we have gone to the trouble of giving a
self-contained, complete
proof is the following.

In the proof of 
\cite[Theorem~A.1]{Alonso-Jeremias-Souto03},
the triangle $A\la X\la C\la\sh A$ with $A\in\Coprod(\ca)$
and $C\in\Coprod(\ca)^\perp$ is produced
as a homotopy colimit of distinguished triangles
$X_n\stackrel{h_n}\la X\la C_n\la\sh X_n$, where as it
happens the maps $h_n\colon X_n\la X$ are the same as
those of our proof of \autoref{T701.9}.
But without enhancements it is not known that this
homotopy limit can always be assumed to be a distinguished triangle,
and hence the argument of
\cite[Theorem~A.1]{Alonso-Jeremias-Souto03}
has a gap.

Our argument fixes the problem by letting $h\colon A\la X$
be any choice of the homotopy colimit of the maps
$h_n\colon X_n\la X$, computing what the functor $\Hom(B,-)$
does to $h$ for any $B\in\sh[-1]\ca$, and then,
after completing $h\colon A\la X$ to a distinguished triangle
$A\la X\la C\la\sh A$, using the information to
show that $C\in\Coprod(\ca)^\perp$.

Of course, in the presence of enhancements
the proof of 
\cite[Theorem~A.1]{Alonso-Jeremias-Souto03}
works just fine.
But in this article the results are all
enhancement-free, and hence we deemed it worthwhile
going to a little trouble to give a
proof of \autoref{T701.9} that is
evidently enhancement-free.

While still on the subject of
enhancement-free proofs of
\autoref{T701.9}: an alternative approach goes as
follows. Still starting with
the above sequence of triangles
$X_n\stackrel{h_n}\la X\stackrel{k_n}\la C_n\la\sh X_n$,
we can focus instead on the map $k_n$. Let $C:=\hoco C_n$, and as
in the proof of 
\cite[Theorem~A.1]{Alonso-Jeremias-Souto03}
it is easy to see that $C\in\Coprod(\ca)^\perp$.
The challenge becomes to show that the homotopy
colimit map $k\colon X\la C$ can be so chosen that, when
we complete to a distinguished triangle $A\la X\stackrel k\la C\la\sh A$,
the object $A$ belongs to $\Coprod(\ca)$. This can be done, but we have not included it here. One
complete, enhancement-free proof of 
\autoref{T701.9} was deemed enough.
\ermk

\rmk{R701.13}
We should also mention that \autoref{T701.9}
is only a baby case of a whole genre of results.
There are several results
in the literature giving other $\ca$, with
$\sh\ca\subseteq\ca$,
for which it is known that
$\big(\Coprod(\ca),\Coprod(\ca)^\perp\big)$ is a \tstr\ on
$\ct$. The basic idea is that, by the results of
\autoref{S699}, this is equivalent to
showing that the inclusion $\Coprod(\ca)\la\ct$ has
a right adjoint, and representability theorems
can be used to produce the right adjoint.
\ermk

The \tstr s constructed in \autoref{T701.9} deserve a special name.

\dfn{R701.15}
Let $\ct$ be a triangulated category with coproducts.
The \tstr\ of the form
$\tau_\ca:=\big(\Coprod(\ca),\Coprod(\ca)^\perp\big)$, where
$\ca\subseteq\ct^c$ is an essentially small full subcategory
satisfying $\sh\ca\subseteq \ca$, will be called the
\emph{\tstr\ on $\ct$ compactly generated by $\ca$}.
\edfn

\rmk{rmk:propgen}
Following \autoref{N700.1}(iv), we will denote by $(-)_\ca^{\leq n}$ the functors associated to such a \tstr . Note that, since $\Coprod(\ca)^\perp$ is closed under coproducts by \autoref{T701.9}, such a functor commutes with coproducts. Furthermore we denote by $\tsth\ct \ca$ its heart and by $\ch_\ca^0$ the corresponding homological functor, in accordance with \autoref{R700.3}(ii).

Furthermore, one can consider the following subcategories:
\[
\ct^{\leq0}_{\ca,c}:=\ct^c\cap\ct^{\leq0}_\ca\qquad \tsth\ct{\ca,c}:=\ch_\ca^{0}\big(\ct^{\leq0}_{\ca,c}\big)\subseteq\tsth\ct \ca.
\]
Then \cite[Proposition~1.9(ii)]{Neeman17} gives us the final
equality in the string
\[
\ct^{\leq0}_{\ca,c}=\ct^c\cap\ct^{\leq0}_\ca=\ct^c\cap\Coprod(\ca)=
\Smr\big(\COprod(\ca)\big).
\]
This clearly implies that $\ct^{\leq0}_{\ca,c}$ is essentially
small, hence the same is true for $\tsth\ct{\ca,c}=\ch_\ca^{0}\big(\ct^{\leq0}_{\ca,c}\big)$.
\ermk

As the reader will see, the existence and basic
properties of these \tstr s will play an important
role in the rest of the article.

\comment{
\exm{ex:compactgentstr}
Let $X$ be a quasi-compact and quasi-separated scheme and
let $Z\subseteq X$ be a closed subscheme such that
$X\setminus Z$ is quasi-compact. Let $\ct=\Dqcs{Z}(X)$ be as
in (a), and let $\ch^n$ be the functor
taking an object $C\in\Dqcs Z(X)$ to its $n\mth$ cohomology
sheaf. It is known that the following formula gives a 
\tstrs on $\ct$:
\[
\begin{array}{ccccc}
\ct^{\leq0}&=&\Dqcs{Z}^{\leq 0}(X)&=&\{C\in\Dqcs Z(X)\mid
\ch^n(C)=0\text{ for all }n>0\}\ ,\\
\ct^{\geq0}&=&\Dqcs{Z}^{\geq 0}(X)&=&\{C\in\Dqcs Z(X)\mid
\ch^n(C)=0\text{ for all }n<0\}\ .
\end{array}
\]
This is usually referred to as the \emph{standard}
\tstrs on $\ct$.

Now let 
$dperfs{Z}{X}\subset\ct$ be the
full subcategory of perfect
complexes, as in (e).
By \cite[Theorem 6.8]{Rouquier08} and \cite[Theorem 3.2(i)]{Neeman22A}, the category $dperfs{Z}{X}$ coincides with
$\ct^c$, the subcategory of compact objects. And
then \cite[Theorem 6.8]{Rouquier08} and \cite[Theorem 3.2(ii)]{Neeman22A} go on to tell us that there exists
an object $G\in\dperfs{Z}{X}$ which generates $\Dqcs{Z}(X)$;
with the $\cb:=G(-\infty,\infty)$ as in \autoref{N314.3}(i), this means
$\ct=\Coprod(\cb)$.

Still with the conventions as in \autoref{N314.3}(i),
let $\ca:=G(-\infty,0]$. Then, by \cite[Theorem 3.2(iii)]{Neeman22A}, the two \tstrs $\tau_\ca$ and $\tau_{\mathrm{stan}}$ are in the same equivalence class. Namely, there is a positive integer $N$ such that the inclusions
\[
(\Dqcs{Z}(X))_\ca^{\leq -N}\subseteq \Dqcs{Z}^{\leq 0}(X)\subseteq (\Dqcs{Z}(X))_\ca^{\leq N}
\]
hold true.
\eexm
}

\section{The heart of a compactly generated \tstr}
\label{S69}

This section is devoted to proving important properties of hearts of \tstr{s}  generated by sets of compact objects as in \autoref{R701.15}.
The main theorem of this section can be found in 
Saor{\'{\i}}n and
{\v{S}}{\v{t}}ov{\'{\i}}{\v{c}}ek~\cite[Theorem~8.31]{Saorin-Stovicek20}.
The proof there is different---we present here a self-contained argument. The main reason is that the technical lemmas in our proof will be needed later in the paper.

\thm{thm:maintstr}
Let $\ct$ be a triangulated category with coproducts,
and let $\ca\subseteq\ct^c$ be an essentially small full subcategory satisfying $\sh{\ca}\subseteq \ca$.
Then the heart $\tsth\ct A$ of the \tstr\ $\tau_\ca=\big(\Coprod(\ca),\Coprod(\ca)^\perp\big)$ is a locally
finitely presented Grothendieck abelian category.
Moreover: $\tsth\ct{A,c}\subset\tsth\ct A$
is precisely the full subcategory of
finitely presented objects in $\tsth\ct A$.
\ethm

Recall that an abelian category $\cb$ is a \emph{Grothendieck category} if it satisfies axiom [AB5] (meaning that $\cb$ has small coproducts and filtered colimits of exact sequences are exact in $\cb$), and
it has a set of generators.  An object $b\in\cb$ is \emph{finitely presented} if the functor $\Hom_\cb(b,-)$ commutes with filtered colimits, and $\cb$ is \emph{locally finitely presented} if the full subcategory of finitely presented objects in $\ca$ is essentially small and all objects of $\cb$ are filtered colimits of finitely presented objects.

\subsection{Some technical results}\label{subsec:techn}

In this section, the standing assumption is that $\ct$ is a triangulated category with coproducts and $\ca\subseteq\ct^c$ is an essentially small full subcategory satisfying $\sh{\ca}\subseteq\ca$.

The following little lemma may also be found in
Saor{\'{\i}}n and
{\v{S}}{\v{t}}ov{\'{\i}}{\v{c}}ek~\cite[part 2b of Proposition~8.27]{Saorin-Stovicek20}. The interested reader is encouraged to compare the
approaches---the short summary is that in \cite{Saorin-Stovicek20}
the ideas are explored much more fully, including a converse to
the lemma below.

\lem{L69.5}
If $f\colon c\la X$ is a morphism in $\ct$ with $c\in\ct^c$ and
$X\in\ct^{\leq0}_\ca$, then there is a factorization of $f$ as $c\la t\la X$
with $t\in\ct^{\leq0}_{\ca,c}$.
\elem

\prf
By construction, $f$ is a morphism
from $c\in\ct^c$ to $X\in\Coprod(\ca)$. By \cite[Lemma~1.8(ii)]{Neeman17}
this map must factor through an object $t\in\COprod(\ca)\subseteq\ct^c\cap\ct^{\leq0}_\ca=\ct_{\ca,c}^{\leq0}$.
\eprf

\lem{L69.7}
Let $f\colon c\la X^{\ge0}$ be a morphism in $\ct$, 
where $c\in\ct^c$ and the truncation $X^{\ge0}$ of $X\in\ct$ is with respect to the \tstr\
$\tstv\ct \ca$. Then there are  morphisms $\ph\colon b\la c$ and
$g\colon b\la X$, with $b\in\ct^c$, such that
\be
\item
The map $\ph^{\geq0}\colon b^{\geq0}\la c^{\geq0}$ is an isomorphism;
\item
The triangle below commutes
\[\xymatrix@C+40pt@R-20pt{
  & c^{\geq0}\ar[dd]^{f^{\geq0}}\\
b^{\geq0} \ar[ru]^-{\ph^{\geq0}}\ar[rd]_{g^{\geq0}} &  \\
 & X^{\geq0} .
}\]
\ee
\elem

\prf
Consider the distinguished triangle
$X^{\leq-1}\la X\la X^{\geq0}\stackrel\delta\la\sh{X^{\leq-1}}$.
The composite $c\stackrel f\la X^{\geq0} \stackrel\delta\la\sh{X^{\leq-1}}$
is a morphism from $c\in\ct^c$ to the object
$\sh{X^{\leq-1}}\in\ct^{\leq-2}_\ca$, and \autoref{L69.5} permits us
to factor it as $c\la d\la \sh{X^{\leq-1}}$ with $d\in\sh[2]{\ct^{\leq0}_{\ca,c}}$.
And now form the distinguished triangle $\sh[-1]{d}\la b\stackrel\ph\la c\la d$.
Since $c$ and $d$ are both compact so is $b$, and
as $d$ and $\sh[-1]{d}$ both belong to $\ct^{\leq-1}_\ca$ the functor
$(-)^{\geq0}$ must take $\ph\colon b\la c$ to an isomorphism.

Next consider the commutative diagram where the rows are distinguished triangles
\[\xymatrix@C+20pt{
b\ar[r]^\ph & c\ar[r]\ar[d]_{f} & d\ar[d]\ar[r] &\sh{b}\\
X\ar[r] &   X^{\geq0}\ar[r]^{\delta} & \sh{X^{\leq-1}}\ar[r] &\sh{X}.
}\]
We may complete it to a morphism of triangles, and in particular
this allows us to produce a commutative square
\[\xymatrix@C+20pt{
b\ar[r]^\ph\ar[d]_g & c\ar[d]_{f}  \\
X\ar[r] &   X^{\geq0}
}\]
We have now constructed the object $b\in\ct^c$ and the morphisms
$\ph\colon b\la c$ and $g\colon b\la X$. The last paragraph showed that $\ph^{\geq0}$
is an isomorphism, proving (i). And (ii) follows from applying
the functor $(-)^{\geq0}$ to the commutative square above.
\eprf

\lem{L69.9}
In the abelian category $\tsth\ct \ca$, the full subcategory
$\tsth\ct{\ca,c}$ satisfies the following
\be
\item
$\tsth\ct{\ca,c}$ is closed in $\tsth\ct \ca$ under extensions.
\item
$\tsth\ct{\ca,c}$ is closed in $\tsth\ct \ca$ under cokernels.
\item
For any object $X\in\tsth\ct{\ca,c}$, the functors $\Hom_{\tsth\ct{\ca}}(X,-)$
and $\Ext_{\tsth\ct{\ca}}^1(X,-)$ commute with coproducts in $\tsth\ct \ca$.  
\setcounter{enumiv}{\value{enumi}}
\ee
\elem

\prf
We begin with (i): let $X,Z\in\tsth\ct{\ca,c}$
and suppose we are given in $\tsth\ct \ca$ a short exact sequence
$0\la X\la Y\la Z\la 0$. This means that in
$\ct$ there must be a morphism $Z\la\sh{X}$ such that
$X\la Y\la Z\la\sh{X}$ is a distinguished triangle (see \autoref{R700.3}).

Now because $X$ and $Z$ belong to $\tsth\ct{\ca,c}$, there must exist
objects $a,c$ in $\ct^{\leq0}_{\ca,c}$ with $X=a^{\geq0}$ and $Z=c^{\geq0}$.
In particular the \tstr\ truncation provides us with a morphism
$c\la c^{\geq0}=Z$, and composing with
the morphism $Z\la \sh{X}$ produces a map
$c\la c^{\geq0}\la \sh{a^{\geq0}}=(\sh{a})^{\geq-1}$
which we will denote $f\colon c\la (\sh{a})^{\geq-1}$. Our construction
so far has produced a commutative diagram
\[
\xymatrix@C+20pt{
c\ar[r]\ar[dr]_f & c^{\geq-1}\ar[r]\ar[d]^{f^{\geq-1}} & c^{\geq0}\ar@{=}[r]\ar[d] & Z\ar[d]\\
   & (\sh{a})^{\geq-1}\ar@{=}[r] & (\sh{a})^{\geq-1}\ar@{=}[r] &\sh{X}.
}
\]
To the morphism $f$ we now apply \autoref{L69.7} and conclude that there  exists
a compact object $\wt c\in\ct^c$ and morphisms
$\ph\colon \wt c\la c$ and $g\colon \wt c\la\sh{a}$, satisfying
\be
\setcounter{enumi}{\value{enumiv}}
\item
The map $\ph^{\geq-1}\colon \wt c^{\geq-1}\la c^{\geq-1}$ is an isomorphism,
and hence $\wt c$ must belong to $\ct^{\leq0}_{\ca,c}$ and
$\ch_\ca^{0}(\ph)\colon \ch_\ca^{0}(\wt c)\la\ch_\ca^{0}(c)$ is an isomorphism.
\item
The triangle below commutes
\[
\xymatrix@C+40pt@R-20pt{
  &  c^{\geq-1}\ar[dd]^{f^{\geq-1}}\\
\wt c^{\geq-1} \ar[ru]^-{\ph^{\geq-1}}\ar[rd]_{g^{\geq-1}} &  \\
 & (\sh{a})^{\geq-1} 
}
\]
\setcounter{enumiv}{\value{enumi}}
\ee
In view of this we deduce the following commutative diagram 
\[\xymatrix@C+20pt{
\wt c\ar[r]\ar[d]_g & \wt c^{\geq-1}\ar[r]\ar[d]^{g^{\geq-1}} & c^{\geq0}\ar@{=}[r]\ar[d] & Z\ar[d]\\
\sh{a}\ar[r]   & (\sh{a})^{\geq-1}\ar@{=}[r] & (\sh{a})^{\geq-1}\ar@{=}[r] &\sh{X}.
}\]
Now complete the commutative square
\[\xymatrix@C+20pt{
\wt c\ar[r]^g \ar[d] & 
\sh{a}\ar[d] \\
Z\ar[r] & \sh{X}
}\]
to a morphism of distinguished triangles
\[\xymatrix@C+20pt{
a\ar[r]\ar[d] & b\ar[r]\ar[d] &\wt c\ar[r]^g \ar[d] & 
\sh{a}\ar[d] \\
X\ar[r] & Y\ar[r] &Z\ar[r] & \sh{X}
}\]
As $a$ and $\wt c$ both belong
to $\ct^{\leq0}_{\ca,c}$ so does $b$. And, by construction, the morphism
$g^{\geq-1}\colon\wt c^{\geq-1}\la (\sh{a})^{\geq-1}$ factors
as $\wt c^{\geq-1}\la \wt c^{\geq0}\la(\sh{a})^{\geq-1}$,
forcing the composite $\wt c^{\leq-1}\la\wt c \la (\sh{a})^{\geq-1}$ to
vanish. Hence the morphism $\ch^{-1}_\ca(g)\colon\ch_\ca^{-1}(\wt c)\la\ch_\ca^{-1}(\sh{a})=\ch_\ca^0(a)$ vanishes.
Thus
the functor $\ch^0_\ca$ takes the morphism
of triangles
to a commutative diagram in $\tsth\ct \ca$ with exact rows
\[\xymatrix@C+20pt{
0\ar[r] &
\ch^0_\ca(a)\ar[r]\ar[d]^\wr & \ch^0_\ca(b)\ar[r]\ar[d] &\ch^0_\ca(\wt c)\ar[r] \ar[d]^\wr & 
0 \\
0\ar[r] & X\ar[r] & Y\ar[r] &Z\ar[r] & 0
}\]
and the fact that the outside maps are isomorphisms forces
the map $\ch^0_\ca(b)\la Y$ to be an isomorphism. Hence
$Y$ is isomorphic to $\ch^0_\ca(b)\in\tsth\ct{\ca,c}$,
proving (i).

Now for (ii). Take a morphism $X\la Y$ in $\tsth\ct{\ca,c}$. There
exist objects $a,b\in\ct^{\leq0}_{\ca,c}$ with $X=a^{\geq0}$ and $Y=b^{\geq0}$.
The map $a\la a^{\geq0}=X$ can be composed with the morphism
$X\la Y$ to produce a map $a\la a^{\geq0}\la b^{\geq0}$ which
we denote $f\colon a\la b^{\geq0}$. And now
\autoref{L69.7} allows us to find an object $\wt a\in\ct^c$
and morphisms $\ph\colon\wt a\la a$ and $g\colon\wt a\la b$ such that
\be
\setcounter{enumi}{\value{enumiv}}
\item
The map $\ph^{\geq0}\colon\wt a^{\geq0}\la a^{\geq0}$ is an isomorphism,
and hence $\wt a$ must belong to $\ct^{\leq0}_{\ca,c}$ and
$\ch_\ca^{0}(\ph)\colon\ch_\ca^{0}(\wt a)\la\ch_\ca^{0}(a)$ is an isomorphism.
\item
The triangle below commutes
\[\xymatrix@C+40pt@R-20pt{
  &  a^{\geq0}\ar[dd]^{f^{\geq0}}\\
\wt a^{\geq0} \ar[ru]^-{\ph^{\geq0}}\ar[rd]_{g^{\geq0}} &  \\
 & b^{\geq0} 
}\]
\setcounter{enumiv}{\value{enumi}}
\ee
Completing $g\colon\wt a\la b$ to a distinguished triangle $\wt a\la b\la c\la\sh{\wt a}$
produces for us an object $c$ which is compact, and as $b$ and $\sh{\wt a}$
both lie in $\ct^{\leq0}_\ca$ we have that so does $c$. Thus
$c\in\ct^{\leq0}_{\ca,c}$, and in cohomology we have an exact sequence
$\ch_\ca^0(\wt a)\la\ch_\ca^0(b)\la\ch_\ca^0(c)\la\ch_\ca^1(\wt a)=0$. This identifies
with $X\la Y\la\ch_\ca^0(c)\la 0$, proving that the cokernel
in $\tsth\ct \ca$ of the map $X\la Y$ is isomorphic to
$\ch_\ca^0(c)\in\tsth\ct{\ca,c}$, establishing (ii).

It remains to prove (iii). Let $c$ be an object in $\ct^{\leq0}_{\ca,c}$
and let $\{X_\lambda,\,\lambda\in\Lambda\}$ be any set
of objects in $\ct^{\geq-1}$. Suppose we are given a map
\[
\xymatrix@C+30pt@R-20pt{
c^{\geq0}\ar[r]^-f &\ds\coprod_{\lambda\in\Lambda}X_\lambda\ .
}
\]
The composite
\[
\xymatrix@C+30pt{
c\ar[r] &
c^{\geq0}\ar[r]^-f &\ds\coprod_{\lambda\in\Lambda}X_\lambda
}
\]
is a morphism from the compact object $c$ to a coproduct, and
must factor through a finite subcoproduct. Therefore there is
a subset $\Lambda'\subseteq\Lambda$, with $\Lambda\setminus\Lambda'$ finite and such that the longer composite
\[
\xymatrix@C+30pt{
c\ar[r] &
c^{\geq0}\ar[r]^-f &\ds\coprod_{\lambda\in\Lambda}X_\lambda
\ar[r]^-\pi &\ds\coprod_{\lambda\in\Lambda'}X_\lambda
}
\]
vanishes, where $\pi$ is the projection to the
direct summand. But the vanishing of
$c\la c^{\geq0}\stackrel{\pi\circ f}\la\coprod_{\lambda\in\Lambda'}X_\lambda$
coupled with the distinguished triangle $c\la c^{\geq0}\la\sh{c^{\leq-1}}$
allows us to produce a commutative
square
\[
\xymatrix@C+40pt{
c^{\geq0}\ar[r]^-f\ar[d] &\ds\coprod_{\lambda\in\Lambda}X_\lambda
\ar[d]^-\pi\\
\sh {c^{\leq-1}}\ar[r]^-\ph &\ds\coprod_{\lambda\in\Lambda'}X_\lambda.
}
\]
Now because $\ct^{\geq-1}_\ca$ is closed under coproducts (see \autoref{T701.9}), the object
$\coprod_{\lambda\in\Lambda'}X_\lambda$ belongs to $\ct^{\geq-1}_\ca$, while
$\sh{c^{\leq-1}}$ clearly belongs to $\ct^{\leq-2}_\ca$. Hence
the map $\ph$ must vanish, and thus so does the composite
from top left to bottom right. This forces the map $f$ to factor
through $\coprod_{\lambda\in\Lambda\setminus\Lambda'}X_\lambda$.

The assertions in (iii) are the special cases where either all
the $X_\lambda$ belong to $\tsth\ct \ca\subseteq\ct_\ca^{\geq-1}$ or they
all belong to $\sh{\tsth\ct \ca}\subseteq\ct_\ca^{\geq-1}$.
\eprf

\cor{C69.10}
The subcategory $\tsth\ct{\ca,c}$ is closed in $\ct$ under
direct summands.
\ecor

\prf
Let $X$ be a direct summand of an object $Y\in\tsth\ct{\ca,c}$.
Then there must exist an idempotent $e\colon Y\la Y$ and an exact
sequence in the abelian category $\tsth\ct \ca$
of the form $Y\stackrel e\la Y\la X\la 0$.
But then \autoref{L69.9}(ii) tells us that $X$ must belong to
$\tsth\ct{\ca,c}$.
\eprf

\lem{L69.11}
Suppose we are given in the abelian category
$\tsth\ct \ca$ a diagram
\[
\xymatrix@C+20pt{
 & & c\ar[d] & \\
X\ar[r] & Y\ar[r] & Z\ar[r] & 0
}
\]
in which the row is exact and $c$ belongs to the subcategory
$\tsth\ct{\ca,c}\subseteq\tsth\ct \ca$. Then we may extend it
in $\tsth\ct \ca$ to a commutative diagram with exact
rows
\[
\xymatrix@C+20pt{
a\ar[r]\ar[d] &b\ar[r]\ar[d] & c\ar[d]\ar[r] & 0\\
X\ar[r] & Y\ar[r] & Z\ar[r] & 0
}
\]
with $a,b\in\tsth\ct{\ca,c}$.
\elem

\prf
Complete the morphism $X\la Y$ to a distinguished triangle $X\la Y\la M\la\sh{X}$ in $\ct$.
Then $M$ is an object of $\ct^{\leq0}_\ca$ and $\ch_\ca^0(M)=Z$. The vertical
morphism in the given diagram therefore produces a map
$c\la Z= M^{\geq0}$. Since $c$ is assumed to belong to
$\tsth\ct{\ca,c}$ there must exist an object $z\in\ct^{\leq0}_{\ca,c}$ with
$c=\ch_\ca^0(z)=z^{\geq0}$. Thus we can form the composite
$z\la z^{\geq0}=c\la Z=M^{\geq0}$. Call this map $f\colon z\la M^{\geq0}$, and
note that $f^{\geq0}$ identifies with the given map $c\la Z$. \autoref{L69.7} now applies: there exist morphisms $\ph\colon\wt z\la z$ and $g\colon\wt z\la M$ with $\wt z\in\ct^c$ such that
\be
\item
The map $\ph^{\geq0}\colon\wt z^{\geq0}\la z^{\geq0}$ is an isomorphism.
This implies that $\wt z\in\ct^{\leq0}_{\ca,c}$ and that
$\wt z^{\geq0}\cong z^{\geq 0}\cong c$.
\item
The triangle below commutes
\[\xymatrix@C+40pt@R-20pt{
  & z^{\geq0}\ar[dd]^{f^{\geq0}}\\
\wt z^{\geq0} \ar[ru]^-{\ph^{\geq0}}\ar[rd]_{g^{\geq0}} &  \\
 & M^{\geq0} 
}\]
\ee

Next apply \autoref{L69.5} to the composite $\wt z\la M\la\sh{X}$: since
$\wt z\in\ct^c$ and $\sh{X}\in\ct^{\leq-1}_\ca$, there is a commutative square
\[\xymatrix@C+20pt{
\wt z\ar[r]\ar[d] & \sh{x}\ar[d] \\
M\ar[r] &\sh{X}
}\]
with $\sh{x}$ a compact object in $\ct^{\leq-1}_\ca$.
Now we complete the commutative square to a morphism of distinguished triangles
\[
\xymatrix@C+20pt{
x\ar[r]\ar[d] & y\ar[r]\ar[d] &\wt z\ar[r]\ar[d] & \sh{x}\ar[d] \\
X\ar[r] & Y\ar[r] & M\ar[r] &\sh{X}.
}
\]
Since $x$ and $\wt z$ both belong to $\ct^{\leq0}_{\ca,c}$ so does $y$. Therefore applying the functor $\ch_\ca^0$ we deduce a commutative diagram in
$\tsth\ct \ca$ with exact rows
\[\xymatrix@C+20pt{
\ch_\ca^0(x)\ar[r]\ar[d] & \ch_\ca^0(y)\ar[r]\ar[d] & c\ar[r]\ar[d] & 0 \\
X\ar[r] & Y\ar[r] & Z\ar[r] &0
}\]
with $\ch_\ca^0(x)$ and $\ch_\ca^0(y)$ both in $\tsth\ct{\ca,c}$.
\eprf

\subsection{The proof of \autoref{thm:maintstr}}\label{subsect:proofthmtstr}

We are now ready to combine the technical results of the previous section and prove our first main theorem.  We keep the assumptions in \autoref{thm:maintstr}, that is $\ct$ is a triangulated category with coproducts and $\ca\subseteq\ct^c$ is an essentially small full subcategory satisfying $\sh{\ca}\subseteq \ca$.

We first prove the existence of generators in the abelian category  $\tsth\ct \ca$.

\lem{L69.13}
The subcategory $\tsth\ct{\ca,c}$ is essentially small and it generates the abelian
category $\tsth\ct \ca$.
\elem

\prf
The fact that $\tsth\ct{\ca,c}$ is essentially small
was observed \autoref{rmk:propgen}. The
next step in the proof is to note
\be
\item
Every nonzero object $X$ in the abelian category $\tsth\ct \ca$ admits
a nonzero map $x\la X$ with $x\in \ca$.
\setcounter{enumiv}{\value{enumi}}
\ee
Assume $X\in \ca^\perp$. Now use the equalities
\[
\ca^\perp=\Coprod(\ca)^\perp=\ct_\ca^{\geq1}\ ,
\]
where the first equality is by \autoref{R701.1}(ii) and
the second is by combining
\autoref{R701.15} with \autoref{N700.1}(ii).
We deduce that $X\in\ct_\ca^{\geq1}$, but as $X$ is assumed to also belong to
$\tsth\ct \ca\subseteq\ct^{\leq0}_\ca$ this forces $X=0$, completing
the proof of (i).

Our next easy observation is
\be
\setcounter{enumi}{\value{enumiv}}
\item
Every nonzero object $X$ in the abelian category $\tsth\ct \ca$ admits
a nonzero map $a\la X$ with $a\in\tsth\ct{\ca,c}$.
\setcounter{enumiv}{\value{enumi}}
\ee
Take a nonzero object $X\in\tsth\ct \ca$. By (i) we have a nonzero
morphism $x\la X$ with
$x\in \ca\subseteq\ct^c\cap\ct^{\leq0}_\ca=\ct^{\leq0}_{\ca,c}$.
This morphism admits a factorization as $x\la x^{\geq0}\la X$,
and as $\ch_\ca^0(x)=x^{\geq0}$ belongs to $\tsth\ct{\ca,c}$ we have
produced our nonzero map $x^{\geq0}\la X$, proving (ii).

And now for the proof of the Lemma. Let $Y$ be an arbitrary object
of the category $\tsth\ct \ca$. Because the subcategory
$\tsth\ct{\ca,c}\subseteq\tsth\ct \ca$
is essentially small, there is only a set $\Lambda$
of maps $y\la Y$ with
$y\in\tsth\ct{\ca,c}$. We need to prove that the morphism
\[\xymatrix@C+40pt{
X:=\ds\coprod_{\lambda\in\Lambda}y_\lambda^{}\ar[r]^-\ph &
Y
}\]
is an epimorphism. In the abelian category $\tsth\ct \ca$ we can complete the morphism
$\ph\colon X\la Y$ to an exact sequence
$X\stackrel\ph\la Y\la Z\la 0$. We will assume $Z\neq0$ and prove a
contradiction. By (ii) there must exist a nonzero morphism
$c\la Z$ with $c\in\tsth\ct{\ca,c}$, and we have a diagram
\[\xymatrix@C+20pt{
 & & c\ar[d] & \\
 X\ar[r] & Y\ar[r] & Z\ar[r] & 0
}\]
to which we can apply \autoref{L69.11}.
It produces for us a commutative diagram with exact rows
\[\xymatrix@C+20pt{
 a\ar[r]\ar[d] &b\ar[r]\ar[d] & c\ar[d]\ar[r] & 0\\
 X\ar[r] & Y\ar[r] & Z\ar[r] & 0,
}\]
with $a,b\in\tsth\ct{\ca,c}$. Now the morphism $b\la Y$ is a morphism
from $b\in\tsth\ct{\ca,c}$ to $Y$, and by construction it must
factor through $X$. Therefore the composite $b\la Y\la Z$ must
vanish. Hence so does the equal composite $b\la c\la Z$, but as
$b\la c$ is an epimorphism we deduce that $c\la Z$ must vanish.
This is a contradication.
\eprf

The next lemma is known. From Positselski and {\v{S}}{\v{t}}ov{\'{\i}}{\v{c}}ek~\cite[Corollary~9.6 and Remark~9.7]{Positselski-Stovicek21} the reader will learn that, if $\ca$ is an abelian category with
coproducts, if $G\subseteq\ca$ is a set of generators, and if for all objects $g\in G$ the functor $\Hom(g,-)$
respects coproducts, then $\ca$ satisfies [AB5].

In \autoref{L69.15} below we include the short proof---we present it
in the special case where $\ca=\tsth\ct \ca$, but the argument is easy to modify to cover the
more general statement. The proof we give is included both for the reader's convenience,
and because the technique will be useful. We will use it again in the
proof of \autoref{L69.917}, and it will play a key role in \autoref{S67}, culminating in the proof of \autoref{P67.13}.

\lem{L69.15}
The abelian category $\tsth\ct \ca$ is a Grothendieck abelian category.
\elem

\prf
\autoref{L69.13} tells us that the abelian category $\tsth\ct \ca$
has a set of generators $\tsth\ct{\ca,c}$. Axiom [AB4] is easy: the categories
$\ct^{\leq0}_\ca$ and $\ct^{\geq0}_\ca$ are both closed in $\ct$ under
coproducts, and hence so is their intersection $\tsth\ct \ca$. And
as coproducts of distinguished triangles are distinguished triangles, if we have a set of
distinguished triangles $X_\lambda\la Y_\lambda\la Z_\lambda\la\sh{X_\lambda}$ with
$X_\lambda ,Y_\lambda, Z_\lambda\in\tsth\ct \ca$, then the coproduct
is a triangle with the same property. That is the coproduct
in $\tsth\ct \ca$ of short exact sequences
$0\la X_\lambda\la Y_\lambda\la Z_\lambda\la0$ is a short exact
sequence in $\tsth\ct \ca$.

Now let us prove [AB5]. Let $I$ be a small filtered category
and let $\fF\colon I\la\tsth\ct \ca$ be a functor; we need to prove that
left derived functors $\colimj\fF$ vanish when $j>0$. Because
the category $\tsth\ct \ca$ satisfies [AB4] these left derived
functors are computed as the cohomology of the complex
obtained from realizing the nerve of the category $I$; there
is a standard cochain complex $\cf$ of the form
\[\xymatrix@C+20pt{
 \cdots \ar[r] & \cf^{-3}\ar[r] &
 \cf^{-2}\ar[r] &\cf^{-1}\ar[r] &\cf^{0}\ar[r] &0
}\]
whose $(-j)\mth$ cohomology is $\colimj \fF$, and where $\cf^{-n}$ is the
coproduct over sequences of composable morphisms in $I$
\[\xymatrix@C+0pt{
i_0^{}\ar[r] & i_1^{}\ar[r] &\cdots\ar[r]&i_{n-1}^{}\ar[r] &i_n^{}
}\]
of $\fF(i_0^{})$. And the differentials of this complex are standard.

Now pick any object $a\in\tsth\ct{\ca,c}$, and apply the functor
$\Hom(a,-)$ to $\cf$. 
\autoref{L69.9}(iii) tells us that
$\Hom(a,-)$ respects coproducts, and therefore it takes
the cochain complex $\cf$ to the standard cochain complex computing
$\colimj \Hom\big(a,\fF(-)\big)$, for the functor
$\Hom\big(a,\fF(-)\big)\colon I\la\ab$. As the category $\ab$ satisfies
[AB5], this filtered colimit has vanishing $\colimj$ for $j>0$.
Therefore the only nonvanishing cohomology of the complex
$\Hom(a,\cf)$ is in degree 0. 

\autoref{L69.13} tells us that the subcategory $\tsth\ct{\ca,c}$
generates $\tsth\ct \ca$. Hence any nonzero cohomology of the
complex $\cf$ would be detected by $\Hom(a,-)$ for some
$a\in\tsth\ct{\ca,c}$. Thus it certainly follows that the cohomology
of $\cf$ is concentrated in degree 0, and the category $\tsth\ct \ca$
satisfies [AB5].
\eprf

The two lemmas above clearly prove the first part of the statement of \autoref{thm:maintstr}: $\tsth\ct \ca$ is a Grothendieck abelian category with a set of generators given by the objects of $\tsth\ct{\ca,c}$ .

\lem{L69.917}
Every object in $\tsth\ct{\ca,c}$ is finitely presented in the abelian category $\tsth\ct \ca$
\elem

One comment before proving the lemma above is in order here. Let $I$ be a filtered category, let $F\colon I\la\tsth\ct \ca$ be a functor, and let $a\in\tsth\ct{\ca,c}$ be an object.
To establish \autoref{L69.917} we need to
prove two things.
\be
\item
Every morphism $a\la\colim\,F$ factors through $F(i)$ for some
$i\in I$.
\item
If $a\la F(i)$ is a morphism such that the composite
$a\la F(i)\la\colim\,F$ vanishes, then there exists in $I$ a  
morphism $i\la j$ with $a\la F(i)\la F(j)$ 
already vanishing.
\ee
Using an argument similar to that in
Positselski and
{\v{S}}{\v{t}}ov{\'{\i}}{\v{c}}ek~\cite[proof of Lemma~9.2(ii)]{Positselski-Stovicek22}
it is possible to unify the two parts. We leave this to the
interested reader, observing that in the proof given below the
two halves are similar.

\prf
Let $I$ be a small filtered category and let $F\colon I\la\tsth\ct \ca$ be
a functor. With the complex $\cf$ as in the proof of
\autoref{L69.15}, meaning it is the standard complex
\[\xymatrix@C+20pt{
 \cdots \ar[r] & \cf^{-3}\ar[r] &
 \cf^{-2}\ar[r] &\cf^{-1}\ar[r] &\cf^{0}\ar[r] &0
}\]
computing $\colimj F$, the fact that the category $\tsth\ct \ca$
is Grothendieck gives the exactness in the category
$\tsth\ct \ca$ of the sequence
\[\xymatrix@C+20pt{
\cf^{-2}\ar[r] &
\cf^{-1}\ar[r] &\cf^{0}\ar[r] &
\colim\, F\ar[r] &0.
}\]
Now let $c$ be an object of $\tsth\ct{\ca,c}$ and suppose
we are given a morphism $c\la \colim\, F$.
Applying \autoref{L69.11} to the diagram
\[\xymatrix@C+20pt{
 & & c\ar[d] & \\
\cf^{-1}\ar[r] &\cf^{0}\ar[r] &
\colim\, F\ar[r] &0
}\]
permits us to construct in $\tsth\ct \ca$ a commutative diagram
with exact rows
\[\xymatrix@C+20pt{
a\ar[r]\ar[d] &b\ar[r]\ar[d] & c\ar[d]\ar[r] &0 \\
\cf^{-1}\ar[r] &\cf^{0}\ar[r] &
\colim\, F\ar[r] &0
}\]
with $a,b\in\tsth\ct{\ca,c}$. The maps $a\la\cf^{-1}$ and $b\la\cf^0$ are
morphisms from
$a,b\in\tsth\ct{\ca,c}$ to coproducts of objects in $\tsth\ct \ca$,
and \autoref{L69.9}(iii) tells us that 
both maps 
must factor through finite subcoproducts. There is a finite
subcategory $I'\subseteq I$ such that the maps from
$a,b$ only see coproducts over $I'$, and as $I$ is filtered we may assume
that the category $I'$ has a terminal object $t$. Letting $G$ be the
composite functor $I'\hookrightarrow I\stackrel F\la\tsth\ct \ca$, we deduce a
factorization
\[\xymatrix@C+20pt{
a\ar[r]\ar[d] &b\ar[r]\ar[d] & c\ar[d]\ar[r] &0 \\
\cg^{-1}\ar[r]\ar[d] &\cg^{0}\ar[r]\ar[d] &
\colim\, G\ar[r]\ar[d] &0\\
\cf^{-1}\ar[r] &\cf^{0}\ar[r] &
\colim\, F\ar[r] &0
}\]
where $\cg$ is the standard complex computing $\colimj G$. But
$\colim\, G=G(t)=F(t)$, and we have factored $c\la\colim\, F$ as
$c\la F(t)\la\colim\, F$.

Next suppose we are given an object $c\in\tsth\ct{\ca,c}$,
an object $s\in I$, and a morphism $c\la F(s)$ such that
the composite
$c\la F(s)\la\colim\, F$ vanishes. As $F(s)$ is a direct
summand of $\cf^0$, the vanishing of the
composite $c\la\cf^0\la\colim\, F$ allows us to factor
the map $c\la\cf^0$ as $c\la\ck\la\cf^0$, where $\ck=\ker(\cf^0\la\colim\,F)$.
Combining this with the exact sequence
\[
\xymatrix@C+20pt{
\cf^{-2}\ar[r] &
\cf^{-1}\ar[r] &\cf^{0}\ar[r] &
\colim\, F\ar[r] &0
}
\]
produces for us a diagram
\[
\xymatrix@C+20pt{
 & & c\ar[d] & \\
\cf^{-2}\ar[r] &
\cf^{-1}\ar[r] &\ck\ar[r] &0
}
\]
to which we can apply \autoref{L69.11}. Thus we construct a commutative
diagram with exact rows
\[
\xymatrix@C+20pt{
a\ar[r]\ar[d] &b\ar[r]\ar[d] & c\ar[d]\ar[r] &0 \\
\cf^{-2}\ar[r] &
\cf^{-1}\ar[r] &\ck\ar[r] &0
}
\]
with $a,b\in\tsth\ct{\ca,c}$. The maps $a\la\cf^{-2}$ and $b\la\cf^{-1}$ are
morphisms from
$a,b\in\tsth\ct{\ca,c}$ to coproducts of objects in $\tsth\ct \ca$,
and \autoref{L69.9}(iii) tells us that 
both maps 
must factor through finite subcoproducts. There is a finite
subcategory $I'\subseteq I$ such that the maps from
$a,b$ only see coproducts over $I'$, and as $I$ is filtered we may assume
that the category $I'$ has a terminal object $t$ and that $s\in I'$.
If $G$ is the composite functor $I'\hookrightarrow I\stackrel F\la\tsth\ct \ca$,
then we obtain a factorization through
\[
\xymatrix@C+20pt{
a\ar[r]\ar[d] &b\ar[r]\ar[d] & c\ar[d]\ar[r] &0 \\
\cg^{-2}\ar[r]\ar[d] &
\cg^{-1}\ar[r]\ar[d] &\ck'\ar[r]\ar[d] &0\\
\cf^{-2}\ar[r] &
\cf^{-1}\ar[r] &\ck\ar[r] &0
}
\]
where
\[\xymatrix@C+20pt{
\cg^{-2}\ar[r] &
\cg^{-1}\ar[r] &\cg^{0}\ar[r] &
\colim\, G\ar[r] &0
}\]
is the standard cochain complex computing $\colimj$ for the functor $G$,
and the sequence $0\la\ck'\la\cg^0\la\colim\, G\la 0$ is exact.
But $\colim\,G=G(t)=F(t)$, and we deduce a commutative diagram
with exact rows
\[\xymatrix@C+20pt{
a\ar[r]\ar[d] &b\ar[r]\ar[d] & c\ar[d]\ar[r] &0\ar[r]\ar[d] & 0 \\
\cg^{-2}\ar[r]\ar[d] &
\cg^{-1}\ar[r]\ar[d] &\cg^0\ar[r]\ar[d] & G(t)\ar[r]\ar[d] &0\\
\cf^{-2}\ar[r] &
\cf^{-1}\ar[r] &\cf^0\ar[r] &\colim\, F\ar[r] &0.
}\]
Thus we have found a morphism $s\la t$ in $I$ such that the composite
$c\la F(s)\la F(t)$ vanishes.

This completes the proof that every object in
$\tsth\ct{\ca,c}$ is finitely
presented in $\tsth\ct \ca$.
\eprf

Next, we establish the following result.

\lem{L69.999}
Every object $Z\in\tsth\ct \ca$ is the filtered colimit
of the objects in $\tsth\ct{\ca,c}$ mapping to it.
\elem

\prf
Let $I$ be the category whose objects are morphisms $f\colon c\la Z$
with $c\in\tsth\ct{\ca,c}$, and whose morphisms are commutative
triangles
\[\xymatrix@C+40pt@R-20pt{
c\ar[rd]^f\ar[dd]_\ph    &   \\ 
&  Z    \\
c'.\ar[ur]_g &
}\]
Using \autoref{L69.9} it is easy to show that the category
$I$ is filtered. Let $F\colon I\la\tsth\ct{\ca,c}$ be the functor
taking the object $f\colon c\la Z$ in $I$ to the object
$c\in\tsth\ct{\ca,c}$, put $Y=\colim\,F$ and let $\rho\colon Y\la Z$
be the obvious map. \autoref{L69.13} guarantees that $\rho$
is an epimorphism. We assert that $\rho$ is an isomorphism.

The proof is by contradiction. Let $X$ be the kernel of
$\rho$ and assume $X$ nonzero. \autoref{L69.13}
establishes that there must exist a nonzero map $a\la X$
with $a\in\tsth\ct{\ca,c}$. Since $X\la Y$ is a monomorphism
the composite $a\la X\la Y=\colim\,F$ is nonzero, but
\autoref{L69.917} tells us that $\Hom(a,-)$ commutes with
filtered colimits. Therefore the map $a\la \colim\,F$ must
factor as $a\la F(i)\la\colim\,F$ for some $i\in I$.
If $i\in I$ is the object $f\colon b\la Z$, then we deduce
in $\tsth\ct{\ca,c}$ a morphism
$a\la b$, and the composite $a\la b\la Z$ must vanish since
it is equal to the vanishing composite $a\la X\la Y\la Z$.
If we complete $a\la b$ to the exact
sequence $a\la b\stackrel\ph\la c\la 0$,
then \autoref{L69.9}(ii) tells us
that $c\in\tsth\ct{\ca,c}$ and by the above the 
map $f\colon b\la Z$ must factor through $\ph\colon b\la c$.
We deduce in $\tsth\ct \ca$
a commutative triangle
\[\xymatrix@C+40pt@R-20pt{
b\ar[rd]^f\ar[dd]_\ph    &   \\ 
&  Z    \\
c\ar[ur]_g &
}\]
This can be viewed as a morphism in $I$ from $\{f\colon b\la Z\}$
to $\{g\colon c\la Z\}$, whose image under the functor $F$
is $\ph\colon b\la c$. Thus we have produced in $I$ a morphism
$i\la j$ such that the composite $a\la F(i)\la F(j)$ vanishes,
proving the vanishing of the composite $a\la X\la Y$.
This completes our contradiction.
\eprf

The proof of \autoref{thm:maintstr}]
is complete if we show that the objects of $\tsth\ct{\ca,c}$ are the only
finitely presented
objects.

Assume therefore
that $P$ is a finitely presented object in $\tsth\ct \ca$.
By \autoref{L69.999} we can find a filtered
category $I$, a functor $F\colon I\la\tsth\ct{\ca,c}$, and an
isomorphism $P\cong\colim\,F$.

But $P$ is finitely presented, and hence the identity map
$P\stackrel\id\la P=\colim\,F$ must factor through some
$F(i)\in\tsth\ct{\ca,c}$. Thus $P$ is a direct summand of an object
$F(i)\in\tsth\ct{\ca,c}$, and \autoref{C69.10} establishes that
$P\in\tsth\ct{\ca,c}$.

\section{Weakly approximable triangulated categories and their subcategories}
\label{S1}

In this section we introduce weakly approximable triangulated categories and some relevant full triangulated subcategories. This is carried out in \autoref{subsec:wa}. In the rest of the paper, we will need a general discussion about products and coproducts in weakly approximable triangulated categories. This is the content of \autoref{S22}.

\subsection{Definitions, examples and the relevant subcategories}\label{subsec:wa}

We recall the main object of study in this paper,
first introduced in \cite[Definition~0.21]{Neeman17A}.

\dfn{dfn:wa}
A triangulated category with coproducts $\ct$ is
\emph{weakly approximable} if there exist an integer $A>0$,
a compact generator $G\in\ct^c$ and a \tstrs
$\tst\ct$ on $\ct$, such that
\be
\item
$G\in\ct^{\leq A}$ and $\Hom(G,\ct^{\leq-A})=0$.
\item
For any object $F\in\ct^{\leq0}$ there exists in $\ct$ a distinguished triangle
$E\la F\la D\la\sh{E}$ with $E\in\ogenu G{}{-A,A}$ and with
$D\in\ct^{\leq-1}$, where $\ogenu G{}{-A,A}$ is as in \autoref{N314.3} (iii).
\setcounter{enumiv}{\value{enumi}}
\ee
\edfn

Axiom (i) tells us that the compact generator $G$ must be bounded above (i.e. $G\in\ct^-$) for the given \tstrs, while axiom (ii) provides the reason for the name `weakly approximable triangulated category': given a bounded above object in $\ct$, its top nontrivial cohomology can be approximated by objects generated by $G$ and its shifts in a prescribed interval.

\rmk{rmk:inveq}
Let $\fF\colon\ct\la\ct'$ be an exact equivalence such that $\ct$ is weakly approximable with integer $A$, generator $G$ and \tstrs $\tst\ct$. Then it is very easy to see from the definition that $\ct'$ is weakly approximable with integer $A':=A$, generator $G':=\fF(G)$ and \tstrs $\tst{(\ct')}$ with $(\ct')^{\leq 0}:=\fF(\ct^{\leq 0})$ and $(\ct')^{\geq 1}:=\fF(\ct^{\geq 1})$.
\ermk

\exm{ex:wa}
The reader should keep in mind the following examples of weakly approximable triangulated categories. In particular, those in (i) and (ii) will be discussed in \autoref{sec:examples}.
\be
\item Let $R$ be a ring. Then it is not difficult to prove that $\D(\Mod{R})$ is weakly approximable (see \cite[Example 3.1]{Neeman17A} for a reference).
\item It is a deeper result that, if $X$ is a quasi-compact and quasi-separated scheme and $Z\subseteq X$ is a closed subset such that $X\setminus Z$ is quasi-compact, then $\Dqcs{Z}(X)$ is weakly approximable. This can be found in \cite[Theorem 3.2(iv)]{Neeman22A}.
\item If $\ct$ is the homotopy category of spectra, then \cite[Example~3.2]{Neeman17A} shows that $\ct$ is weakly approximable. Actually more is known: this category, as well as any of the ones in (i), is not only weakly approximable but \emph{approximable.} We will not deal with the stronger version of approximability in this paper.
\ee
\eexm

\rmk{R1.1}
A weakly approximable triangulated category $\ct$ comes with two \tstr s. The first $\tst\ct$ is the
one forming part of \autoref{dfn:wa}. The second
is obtained from the given compact generator $G$ as follows:
let $\ca=G(-\infty,0]$
as in \autoref{N314.3}(i), that is take the minimal
$\ca$ with $G\in\ca$ and such that $\T\ca\subset\ca$.
And then form the \tstrs $\tstv\ct G$ with
$\ct^{\leq0}_G=\Coprod(\ca)$, using \autoref{T701.9}.
And \cite[Proposition~2.4]{Neeman17A} shows that the given \tstrs $\tst\ct$ is in the same equivalence class as $\tstv\ct G$,
meaning that there exists a positive integer $N>0$ such that $\ct^{\leq-N}\subseteq\ct^{\leq0}_G\subseteq\ct^{\leq N}$.

Note: in \cite[Remark~0.15]{Neeman17A} it is shown
that, if $G$ and $H$ are
both compact generators for a triangulated category $\ct$ with coproducts, then the \tstr s $\tstv\ct G$ and
$\tstv\ct H$ are equivalent. 
Thus for any
triangulated category $\ct$, with a
single compact generator,
we obtain a \emph{preferred equivalence class of \tstr s}---it
is the equivalence class that contains $\tstv\ct G$ for every compact generator $G$. See \cite[Definition~0.14]{Neeman17A}.
And in the case of a weakly approximable triangulated category $\ct$, the \tstr\ given in \autoref{dfn:wa}, for which
\autoref{dfn:wa}~(i) and (ii) both hold, must belong to the
preferred equivalence class.
\ermk

\rmk{the subcategories we care about}
Let $\ct$ be a weakly approximable triangulated category,
and let $\tst\ct$ be  a \tstrs in the preferred equivalence class. Then we can define the following full subcategories of $\ct$, which are clearly independent of the choice of
\tstrs in the preferred equivalence
class.
\begin{itemize}
\item \emph{Bounded above objects}: $\ct^-:=\bigcup_{m=1}^\infty\ct^{\leq m}$;
\item \emph{Bounded below objects}: $\ct^+:=\bigcup_{m=1}^\infty\ct^{\geq-m}$;
\item \emph{Bounded objects}: $\ct^b:=\ct^-\cap\ct^+$,
\item \emph{Compact objects}: $\ct^c$, as defined in
  \autoref{D701.5} and \autoref{R701.7};
\item \emph{Pseudo-compact objects}: $\ct^-_c$, where an object $F\in\ct$ belongs to $\ct^-_c$ if, for any integer $m>0$, there exists in $\ct$ a distinguished triangle $E\la F\la D$ with $E\in\ct^c$ and with $D\in\ct^{\leq-m}$. Thus
\[
\ct^-_c=\bigcap_{m=1}^\infty\big(\ct^c*\ct^{\leq-m}\big);
\]
\item \emph{Bounded pseudo-compact objects}: $\ct^b_c:=\ct^-_c\cap\ct^b$.
\item \emph{Bounded compact objects}: $\ct^{c,b}:=\ct^b_c\cap\ct^c$.
\end{itemize}
The very definitions, listed above, tell us how to obtain
out of $\ct$ all of the subcatergories on the list.
In other words: the part of \autoref{thm:main1} that
says there are recipes for passing from $\ct$ to
its various subcategories is by the very nature of their
definition, as above.
\ermk

The following easy result clarifies the relation between $\ct^-$, $\ct^c$ and $\ct^-_c$.

\lem{lem:compabove}
Let $\ct$ be a weakly approximable triangulated category. Then $\ct^c\subseteq\ct^-_c\subseteq\ct^-$.
\elem

\prf
By definition of $\ct^-_c$ it is clearly enough to show $\ct^c\subseteq\ct^-$. But this inclusion holds because $\ct^c$ is classically generated by a compact generator $G$ of $\ct$ and $G\in\ct_G^{\leq 0}\subseteq\ct_G^-=\ct^-$ (where the last equality holds by \autoref{R1.1}).
\eprf

In general, it is not always true that compact objects are bounded below as well, as illustrated by the example below. On the other hand, it is clear that $\ct^{c,b}=\ct^c\subseteq\ct^b_c$ whenever $\ct^c\subseteq\ct^+$.

\exm{R1.972}
If $\ct$ is either $\D(\Mod{R})$ as in \autoref{ex:wa}(i) or $\Dqcs Z(X)$ as in \autoref{ex:wa}(ii), we have $\ct^c\subseteq\ct^b_c$. Indeed, $\D(\Mod{R})^c=\dperf{R}$ and $\Dqcs Z(X)^c=\dperfs ZX$, and in both cases perfect complexes are (cohomologically) bounded.

But the inclusion $\ct^c\subseteq\ct^b_c$ does not hold when $\ct$ is the homotopy category of spectra as in \autoref{ex:wa}(iii). In this case it is known that $\ct^c$ is the homotopy category of finite spectra, and it is not difficult to compute the subcategory $\ct^b_c$. It turns out that in this example we have $\ct^b_c\neq\{0\}\ne\ct^c$ but $\ct^{c,b}=\{0\}$. There are also many examples of algebraic triangulated categories $\ct$ such that $\ct^{c,b}$ is a proper full subcategory of both $\ct^b_c$ and $\ct^c$.
\eexm

\lem{lem:inveq}
If $\ct$ is a weakly approximable triangulated category, then all the subcategories listed in \autoref{the subcategories we care about} are thick.
\elem

\prf
By \cite[Observation~0.12]{Neeman17A}, the subcategories $\ct^-$, $\ct^+$ and $\ct^b$ are thick. It is straightforward to prove that $\ct^c$ is thick using the definition. By \cite[Proposition~2.10]{Neeman17A}, $\ct^-_c$ is thick. It then follows that $\ct^b_c$ and $\ct^{c,b}$ are thick as well.
\eprf

The fact that each of these subcategories is thick will play a key role in passing between them.

Finally now, once all the terms have been explained, the
reader might wish to glance back at the table in the Introduction, which gives what the intrisic subcategories
of \autoref{the subcategories we care about} turn out to
be when $\ct=\D(\Mod R)$ and when $\ct=\Dqcs Z(X)$.
 
\subsection{Products and coproducts in weakly approximable triangulated categories}
\label{S22}

This section collects some nice but technical properties of products and coproducts in weakly approximable triangulated categories which will be used in the rest of the paper. This is the first instance where axiom (ii) in \autoref{dfn:wa} plays a crucial role.

\lem{L22.2}
Let $\ct$ be a compactly generated triangulated category, and let
$\tst\ct$ be any \tstr.
If $\{t_\lambda\st\lambda\in\Lambda\}$ is any set of objects in
$\ct^{\geq0}$,
then the product $\prod_{\lambda\in\Lambda}t_\lambda$ belongs to
$\ct^{\geq0}$.

Now assume that $\ct$ is a weakly approximable triangulated category,
and that the \tstr\ $\tst\ct$ belongs to the preferred
equivalence class. There exists an integer
$B>0$ such that,
if $\{t_\lambda\st\lambda\in\Lambda\}$ is any set of objects in
$\ct^{\leq0}$,
then the product $\prod_{\lambda\in\Lambda}t_\lambda^{}$ belongs to
$\ct^{\leq B}$.
\elem

\prf
First observe that products exists in $\ct$ (see \cite[Remark 5.1.2]{Krause10}). Now, for any \tstr\
$\tst\ct$ we have that $\ct^{\geq0}$ is equal to
$\left(\ct^{\leq-1}\right)^\perp$, and is therefore closed under whatever
products exist in $\ct$ by \autoref{C699.11}. Thus the first part of the statement follows.

Assume that the category $\ct$ is weakly approximable, and that the
\tstr\ $\tst\ct$ belongs to the preferred equivalence class. Let $G$ be
a compact generator for $\ct$. By \cite[Proposition~2.6]{Neeman17A}
we may choose an integer $A>0$ such that $\Hom(\sh[-A]{G},\ct^{\leq0})=0$,
and therefore $\Hom(\sh[-\ell]{G},-)$ is zero on $\ct^{\leq0}$, for every $\ell\geq A$. In other words,
\[
\ct^{\leq0}\subseteq G[A,\infty)^\perp
\]
where $G[A,\infty)$ is as in \autoref{N314.3}.
On the other hand, \cite[Lemma~3.9(iv)]{Burke-Neeman-Pauwels18}, which uses (ii) of \autoref{dfn:wa}, allows us to choose an
integer $B>A$ such that
\[
G[A,\infty)^\perp\subseteq
\ct^{\leq B}.
\]
Given any collection 
of objects
$\{t_\lambda\st\lambda\in\Lambda\}$, all belonging to
$\ct^{\leq0}$, they belong to
the bigger $G[A,\infty)^\perp$,
which is closed under products. Hence
$\prod_{\lambda\in\Lambda}t_\lambda^{}$ belongs
to $G[A,\infty)^\perp$,
and therefore to the bigger $\ct^{\leq B}$.
\eprf

\lem{L22.3}
Let $\ct$ be a weakly approximable triangulated category, and let
$\tst\ct$ be a \tstr\ on $\ct$ in the preferred equivalence class.
Suppose $\{t_m\st 1\leq m<\infty\}$ is a countable sequence of
objects in $\ct$ such that, for any
integer $n>0$, one of the two conditions below holds:
\be
\item
all but finitely many $t_m$ lie in $\ct^{\geq n}$;
\item
all but finitely many $t_m$ lie in $\ct^{\leq-n}$.
\ee
Then the natural map
\[\xymatrix@C+30pt{
	\ds\coprod_{m=1}^\infty t_m\ar[r] &
	\ds\prod_{m=1}^\infty t_m
}\]
is an isomorphism.
\elem

\prf
Let $G\in\ct$ be a compact generator, and let $\ell\in\zz$ be any integer.
As $\sh[\ell]{G}$ is a compact generator for $\ct$ and $\tst\ct$ is a
\tstrs in the preferred equivalence class, \cite[Proposition~2.6]{Neeman17A}
allows us to find an integer $A>0$ such that $\sh[\ell]{G}\in\ct^{\leq A}$ and
$\Hom(\sh[\ell]{G},-)$ is zero on
$\ct^{\leq -A}$. Combining these, we have that $\Hom(\sh[\ell]{G},-)$ is trivial
both on $\ct^{\leq -A-1}$ and on $\ct^{\geq A+1}$. And because either (i) or (ii) holds, this means that $\Hom(\sh[\ell]{G},-)$ is zero on all but finitely many
of the $t_m$. So we may choose an integer $N>0$ such that $\Hom(\sh[\ell]{G},t_m)=0$ for all $m>N$. But then $\Hom(\sh[\ell]{G},-)$ is trivial both on $\coprod_{m=N+1}^\infty t_m$ and on $\prod_{m=N+1}^\infty t_m$. And since the map
\[\xymatrix@C+30pt{
	\ds\coprod_{m=1}^N t_m\ar[r] &
	\ds\prod_{m=1}^N t_m
}\]
is an isomorphism, we discover that $\Hom(\sh[\ell]{G},-)$ must take the map
\[\xymatrix@C+30pt{
	\ds\coprod_{m=1}^\infty t_m\ar[r]^-\ph &
	\ds\prod_{m=1}^\infty t_m
}\]
to an isomorphism. Since $G$ is a compact generator and the above is true for
every $\ell$, the map $\ph$ must be an isomorphism.
\eprf

\lem{L22.5}
Let $\ct$ be a weakly approximable triangulated category,
and let $\tst\ct$ be a \tstrs in the preferred equivalence class.
Suppose $\{t_m\st 1\leq m<\infty\}$ is a countable sequence of
objects in $\ct^-_c$.
Assume further that, for any
integer $n>0$, all but finitely many $t_m$ lie in $\ct^{\leq-n}$.
Then the isomorphic objects
\[
Q=\coprod_{m=1}^\infty t_m\cong
\prod_{m=1}^\infty t_m
\]
lie in $\ct^-_c$.
\elem

\prf
The fact that the natural map in an isomorphism, from the coproduct
to the product, is by \autoref{L22.3}. What needs proof
is that one of the objects belongs to $\ct^-_c$.

Pick any integer
$n>0$, and then choose an integer $N>0$ such that, for all
$m\geq N$, the object $t_m$ belongs to $\ct^{\leq-n}$.
The finite coproduct
$\ov Q_n=\coprod_{m=1}^{N-1} t_m$ belongs
to $\ct^-_c$, and hence there exists a distinguished triangle
$\ov E_n\la \ov Q_n\la\ov D_n$ with $\ov E_n\in\ct^c$ and with
$\ov D_n\in\ct^{\leq-n}$.

Since $\ct^{\leq-n}$ contains $t_m$ for all $m\geq N$ and is closed under
coproducts, we have that $Q_n=\coprod_{m=N}^\infty t_m$ must belong to
$\ct^{\leq-n}$
But now form the direct sum of the two distinguished triangles
\[
\xymatrix@R-20pt{
	\ov E_n\ar[r] & \ov Q_n\ar[r] &\ov D_n\ar[r] & \sh{\ov E_n} \\
	0 \ar[r] & Q_n\ar[r] & Q_n\ar[r] & 0  
}
\]
to obtain the distinguished triangle
$\ov E_n\la\ov Q_n\oplus Q_n\la \ov D_n\oplus Q_n\la\sh{\ov E_n}$. This does
the trick; we have
that $\ov E_n\in\ct^c$ and $\ov D_n\oplus Q\in\ct^{\leq-n}$, hence $Q=\ov Q_n\oplus Q_n\in\ct^c*\ct^{\leq-n}$. And since $n>0$ is arbitrary
this proves that $Q\in\ct^-_c$.
\eprf

\lem{L27.9}
Let $\ct$ be a weakly approximable triangulated category. Then the inclusion
functors
\[
\ct^-\hookrightarrow\ct,\qquad\ct^b\hookrightarrow\ct,\qquad\ct^+\hookrightarrow\ct
\]
respect both products and coproducts. More precisely,
if $\{t_\lambda\st\lambda\in\Lambda\}$ is a set of objects
in one of the three subcategories $\ct^b$, $\ct^-$ or $\ct^+$,
and if either the coproduct or the product exists in
the subcategory, then both the coproduct and the product
exist in the subcategory and agree (respectively) with
the coproduct and product in $\ct$.
\elem

\prf
We will treat the inclusion $\ct^b\hookrightarrow\ct$, leaving to the reader the other two completely similar cases. Recall that if $\ct$ is a weakly approximable triangulated category, then, in
the preferred equivalence class of \tstr s, we may choose
a \tstrs $\tst\ct$ such that $\ct^{\geq0}$ is closed
under coproducts (by \autoref{T701.9} and \autoref{R1.1}). Furthermore, any \tstrs $\tst\ct$, with $\ct^{\geq0}$ closed under coproducts, is such that the truncation functors $(-)^{\leq0}$ and $(-)^{\geq0}$ both respect coproducts. Thus, let us choose and fix a
\tstr\ $\tst\ct$ in the preferred equivalence class, such that
$\ct^{\geq0}$ is closed under coproducts.

Suppose $\{t_\lambda\st\lambda\in\Lambda\}$ is a set of objects
in $\ct^b$, and suppose $P\in\ct^b$ is either the coproduct or the
product of the $t_\lambda$ in the category $\ct^b$. Then there
must exist an integer $n>0$ with $P\in\ct^{\leq n}\cap\ct^{\geq-n}$.
But every $t_\lambda^{}$ is a direct summand of $P$, and hence
every $t_\lambda^{}$ must belong to $\ct^{\leq n}\cap\ct^{\geq-n}$.

Because both $\ct^{\leq n}$ and $\ct^{\geq-n}$ are closed under
coproducts in $\ct$,
the coproduct $t_1^{}=\coprod_{\lambda\in\Lambda}t_\lambda$ in $\ct$
must belong to $\ct^{\leq n}\cap\ct^{\geq-n}\subseteq\ct^b$. Now
$\ct^{\geq-n}=\big(\ct^{<-n}\big)^\perp$ is closed under
products in $\ct$, while \autoref{L22.2}
produced an integer $B>0$ such
that the product in $\ct$ of objects in $\ct^{\leq n}$ must lie in
$\ct^{\leq n+B}$. This gives us that the product
$t_2^{}=\prod_{\lambda\in\Lambda}t_\lambda$ in $\ct$
must belong to $\ct^{\leq n+B}\cap\ct^{\geq-n}\subseteq\ct^b$.
\eprf

\section{The subcategories of $\ct^-_c$ and of $\ct^c$}
\label{sect:ct-c}

In this section we give the
recipes promised by the solid arrows
of \autoref{thm:maintstr},
in the case of
the subcategories of $\ct^-_c$ and of
$\ct^c$. We give these recipes in \autoref{S2}. As
a prelude we
include the short \autoref{S79}, which clarifies the
properties of $\ct^-_c$ as a set of generators for $\ct$.

\subsection{Preliminaries to understanding the subcategory $\ct^-_c\subseteq\ct$}
\label{S79}

The results of this section are not essential to the article, but
do help form a better understanding of the subcategory $\ct^-_c\subseteq\ct$.

We begin with a technical observation.

\lem{L79.3}
Let $\ct$ be a triangulated category with countable coproducts,
and let $\tst\ct$ be a \tstr\ on $\ct$.
Suppose $F\in\ct$
is an object, and assume we are given a sequence
$E_1\la E_2\la E_3\la\cdots$ of objects of $\ct$, with compatible maps
$E_m\la F$, and
such that
\be
\item
Any choice of a factorization through the homotopy colimit gives
an isomorphism
\[
F\cong\hoco E_m\ .
\]
\item
If we complete $E_m\la F$ to a distinguished triangle $E_m\la F\la D_m$,
then $D_m\in\ct^{\leq-m}$.
\ee
Then, in the distinguished triangle
\[\xymatrix@C+20pt{
	\ds\coprod_{m=1}^\infty E_m \ar[r]^-{\id-\shi} &
	\ds\coprod_{m=1}^\infty E_m \ar[r] &
	F\ar[r]^-\delta &
	\ds\coprod_{m=1}^\infty \sh{E_m}
}\]
defining the homotopy colimit, we have that the morphism $\delta$
must factor as
\[\xymatrix@C+20pt{
	F\ar[r] &\ds\coprod_{m=1}^\infty \sh{(E_m^{\leq-m+1})}\ar[r] &
	\ds\coprod_{m=1}^\infty \sh{E_m}.
}\]
\elem

\prf
Consider the commutative diagram
\[\xymatrix@C+10pt{
	F\ar[r]^-\delta &
	\ds\coprod_{m=1}^\infty \sh{E_m}
	\ar[r]^-{\id-\shi}\ar[d] &
	\ds\coprod_{m=1}^\infty \sh{E_m}\ar[d]&   \\
	& \ds\coprod_{m=1}^\infty \sh{(E_m^{\geq-m+2})}\ar[r]^-\alpha &
	\ds\coprod_{m=1}^\infty \sh{(E_m^{\geq-m+3})}\ar[r]^-\beta &
	\ds\coprod_{m=2}^\infty \sh{(E_m^{\geq-m+3})} ,
}\]
where the horizontal maps $\alpha$ and $\beta$ are the obvious ones. This
means that $\beta$ is the projection to the direct summand, where we
discard the term corresponding to $m=1$. And
$\alpha$ is the map whose components $E_m^{\geq-m+2}\la E_m^{\geq-m+3}$ and
$E_m^{\geq-m+2}\la E_{m+1}^{\geq-m+2}$ are the \tstr\ truncations of the
obvious maps. As $(\id-\shi)\circ\delta=0$ (being the composite
of two morphisms in a distinguished triangle), the composite from top left to bottom right in the diagram must vanish. 

Now observe that the map $\beta\circ\alpha$ is an isomorphism. Indeed, we first observe that our
hypotheses (i) and (ii) in the statement guarantee that the map
$E_m^{\geq-m+2}\la E_{m+1}^{\geq-m+2}$ must be an isomorphism, as both map
isomorphically to $F^{\geq-m+2}$. Call the coproduct over $m$ of
these isomorphisms $\ph$. Then
the map $\beta\circ\alpha$ can be written as $\ph\circ(\id-\tau)$, where $\tau$ is some
map we do not care about beyond the fact the infinite sum
$(\id+\tau+\tau^2+\cdots)$ exists, as the sum is finite on each summand
of the first coproduct. Therefore $\beta\circ\alpha=\ph\circ(\id-\tau)$ is an isomorphism as we claimed.

It then follows that the composite 
\[\xymatrix@C+20pt{
	F\ar[r]^-\delta &
	\ds\coprod_{m=1}^\infty \sh{E_m}
	\ar[r] &
	\ds\coprod_{m=1}^\infty \sh{(E_m^{\geq-m+2})}
}\]
must vanish. The conclusion of the lemma is then immediate.
\eprf

\lem{L79.5}
Let $\ct$ be a triangulated category with coproducts,
and let $\tst\ct$ be a \tstr\ on $\ct$.
Suppose $F\in\ct$
is an object, and assume we are given a sequence
$E_1\la E_2\la E_3\la\cdots$ of objects of $\ct$
as in \autoref{L79.3}.
Assume further that all the objects $E_m$ are compact in $\ct$ and that we are given in $\ct$ a set of objects
$\{X_\lambda\st\lambda\in\Lambda\}$ and a morphism
\[\xymatrix@C+20pt{
	F\ar[r]^-f &
	\ds\coprod_{\lambda\in\Lambda} X_\lambda
}\]
such that, for each integer $m$, the composite
$E_m\la F\stackrel f\la\coprod_{\lambda\in\Lambda} X_\lambda$ vanishes.

Then there is a countable subset $\Lambda'\subseteq\Lambda$ and
a factorization of $f$ as
\[\xymatrix@C+20pt{
	F\ar[r] &\ds\coprod_{\lambda\in\Lambda'}X_\lambda^{\leq-n_\lambda^{}}\ar[r] &
	\ds\coprod_{\lambda\in\Lambda'}X_\lambda\ar@{^{(}->}[r] &
	\ds\coprod_{\lambda\in\Lambda}X_\lambda
}\]
where the inclusion is of the countable subcoproduct, the maps induced
by the truncation are the obvious, and for any integer $m>0$ there exist only
finitely many $\lambda\in\Lambda'$ with $n_\lambda^{}<m$.
\elem

\prf
The vanishing of each composite $E_m\la F\stackrel f\la\coprod_{\lambda\in\Lambda} X_\lambda$ yields a vanishing composite
\[\xymatrix@C+20pt{
	\ds\coprod_{m=1}^\infty E_m \ar[r] &
	F\ar[r]^-f  &
	\ds\coprod_{\lambda\in\Lambda}X_\lambda
}\]
which, in the notation of \autoref{L79.3}, means that the morphism
$f$ must factor as
\[\xymatrix@C+20pt{
	F\ar[r]^-\delta  &
	\ds\coprod_{m=1}^\infty \sh{E_m} \ar[r] &
	\ds\coprod_{\lambda\in\Lambda}X_\lambda.
}\]
Since each $E_m$ is compact, for each $m$ the map from $E_m$ to the
coproduct must factor through a finite subcoproduct. Taking the
union over $m$ gives a countable subset $\Lambda'\subseteq\Lambda$.

And the remainder of the current Lemma comes from combining the above
with the conclusion of \autoref{L79.3}, which allows us to factor
the above further as
\[\xymatrix@C+20pt{
	F\ar[r]  &
	\ds\coprod_{m=1}^\infty \sh{(E_m^{\leq-m+1})} \ar[r] &
	\ds\coprod_{m=1}^\infty \sh{E_m} \ar[r] &
	\ds\coprod_{\lambda\in\Lambda}X_\lambda.
}\]
This concludes the proof.
\eprf

If we assume more about the \tstr\ on $\ct$, and about its
relation with compact objects, then the factorization obtained
in \autoref{L79.5} can be made without restrictions on the map $f$.

\lem{L79.7}
Let $F$ and the sequence $E_1\la E_2\la E_3\la\cdots$ mapping to it be
as in \autoref{L79.5}. Assume further that
\be
\item
The \tstrs $\tst\ct$ is such that $\ct^{\geq0}$ is closed under coproducts.
\item
For every compact object $C\in\ct$ there exists an integer $m>0$ with
$\Hom(C,\ct^{\leq-m})=0$.
\item
Given a countable collection of objects $\{Y_i\in\ct\st 1\leq i<\infty\}$,
such that for every integer $n>0$ we have that all but finitely many of  
the $Y_i$ belong to $\ct^{\leq-n}$, the natural map
\[\xymatrix@C+20pt{
	\ds\coprod_{i=1}^\infty Y_i \ar[r] & \ds\prod_{i=1}^\infty Y_i
}\]
is an isomorphism.
\ee

Suppose next that we are given in $\ct$ a set of objects
$\{X_\lambda\st\lambda\in\Lambda\}$ and a morphism
\[\xymatrix@C+20pt{
	F\ar[r]^-f &
	\ds\coprod_{\lambda\in\Lambda} X_\lambda.
}\]
Then there is a countable subset $\Lambda'\subseteq\Lambda$ and
a factorization of $f$ as
\[\xymatrix@C+20pt{
	F\ar[r] &\ds\coprod_{\lambda\in\Lambda'}X_\lambda^{\leq-n_\lambda^{}}\ar[r] &
	\ds\coprod_{\lambda\in\Lambda'}X_\lambda\ar@{^{(}->}[r] &
	\ds\coprod_{\lambda\in\Lambda}X_\lambda,
}\]
where the inclusion is of the countable subcoproduct, the maps induced
by the truncation are the obvious, and for any integer $m>0$ there exist only
finitely many $\lambda\in\Lambda'$ with $n_\lambda^{}<m$.
\elem

\prf
For every integer $m>0$ we have a distinguished triangle $E_m\la F\la D_m$ with
$E_m\in\ct^c$ and with $D_m\in\ct^{\leq-m}$. Now consider the
composite
\[\xymatrix@C+20pt{
	E_m\ar[r] & F\ar[r]^-f &
	\ds\coprod_{\lambda\in\Lambda} X_\lambda\ar[r] &
	\ds\coprod_{\lambda\in\Lambda} X_\lambda^{\geq-m+1}.
}\]
As $E_m$ is compact this composite factors through a finite subcoproduct,
which we denote $\Lambda_m\subseteq\Lambda$. Without loss of generality we may assume
$\Lambda_m\subseteq\Lambda_{m+1}$.
It follows that composite
\[\xymatrix@C+20pt{
	E_m\ar[r] & F\ar[r]^-f &
	\ds\coprod_{\lambda\in\Lambda} X_\lambda\ar[r] &
	\ds\coprod_{\lambda\in\Lambda} X_\lambda^{\geq-m+1} \ar[r] &
	\ds\coprod_{\lambda\in\Lambda\setminus\Lambda_m} X_\lambda^{\geq-m+1}
}\]
vanishes. But then the map from $F$ to the term on the right must
factor as
\[\xymatrix@C+20pt{
	F\ar[r] & D_m\ar[r] &
	\ds\coprod_{\lambda\in\Lambda\setminus\Lambda_m} X_\lambda^{\geq-m+1}\ .
}\]
By hypothesis (i) in the statement, the category $\ct^{\geq-m+1}$ is closed under coproducts,
hence the coproduct on the right is an object in $\ct^{\geq-m+1}$. But $D_m$ belongs
to $\ct^{\leq-m}$, and the map
$D_m\la
\coprod_{\lambda\in\Lambda\setminus\Lambda_m} X_\lambda^{\geq-m+1}$
must vanish. Hence, for each $\lambda\in\Lambda\setminus\Lambda_m$,
the map $F\la X_\lambda$ can be factored as $F\la X_\lambda^{\leq-m}\la X_\lambda$.

Now for each of the finitely many elements
$\lambda\in\Lambda_{m+1}\setminus\Lambda_m$
choose a factorization of the map $F\la X_\lambda$ as
$F\stackrel{\ph_\lambda^{}}\la X_\lambda^{\leq-m}\la X_\lambda$.
We declare that $\Lambda_0=\emptyset$, and when
$\lambda$ is in the finite set $\Lambda_1=\Lambda_1\setminus\Lambda_0$ we set
$\ph_\lambda^{}\colon F\la X_\lambda$ to be
untruncated. The morphisms $\ph_\lambda^{}$, where
$\lambda$ is in the union $\cup_{m=1}^\infty(\Lambda_m\setminus\Lambda_{m-1})$, combine to a
single morphism
to the product of these countably many trucations. And by hypothesis (iii) in the statement this product agrees with the coproduct. We have produced a countable subset
$\Lambda'\subseteq\Lambda$ and a composite
\[\xymatrix@C+20pt{
	F\ar[r] &\ds\coprod_{\lambda\in\Lambda'}X_\lambda^{\leq-n_\lambda^{}}\ar[r] &
	\ds\coprod_{\lambda\in\Lambda'}X_\lambda\ar@{^{(}->}[r] &
	\ds\coprod_{\lambda\in\Lambda}X_\lambda.
}\]
Call this composite $g$. The construction tells us that $(f-g)^{\geq-n+1}=0$ for
every $n>0$. Hence $f-g$ must factor through an object in $\ct^{\leq-n}$ for every
integer $n$.

Now each $E_m$ is compact, and hypothesis (ii) above says that
there exists an integer $n$ such that
$\Hom(E_m,\ct^{\leq-n})=0$. Hence the composite
$E_m\la F\stackrel{(f-g)}\la\coprod_{\lambda\in\Lambda}X_\lambda$ vanishes. This puts us
in the situation of \autoref{L79.5}, and the map $f-g$ must have
a factorization of the required form.
\eprf

We are about to specialize to the case of weakly approximable triangulated
categories. And for this, the next little lemma will help.

\lem{L79.8}
Let $\ct$ be a weakly approximable triangulated category, and let
$\tst\ct$ be a compactly generated 
\tstr\ in the preferred equivalence class. Then any morphism $F\la X$, with $F\in\ct^-_c$ and $X\in\ct^{\leq0}$, factors
as $F\la F'\la X$ with $F'\in\ct^-_c\cap\ct^{\leq0}$.
\elem

\prf
Because $F$ belongs to $\ct^-_c$, there exists a distinguished triangle $E\la F\la D$ with
$E\in\ct^c$ and $D\in\ct^{\leq0}$. And since $F\in\ct^-_c$ and
$E\in\ct^c\subseteq\ct^-_c$, the triangle tells us that
$D\in\ct^-_c\cap\ct^{\leq0}$.

Now form the composite $E\la F\la X$. It is a morphism from the compact
object $E$ to $X\in\ct^{\leq0}$, for the compactly generated \tstr\
$\tst\ct$. By \autoref{L69.5} we may factor the map $E\la X$
as $E\la E'\la X$ with $E'\in\ct^c\cap\ct^{\leq0}$. This gives
us a commutative square
\[\xymatrix{
	E\ar[r]\ar[d] & F\ar[d] \\
	E'\ar[r] & X.
}\]
Now form the homotopy pushout of the maps $E'\longleftarrow E\la F$. By \cite[Lemma 1.4.4]{Neeman99}
we obtain a morphism of distinguished triangles where the square on the left
is homotopy cartesian
\[\xymatrix{
	E\ar[r]\ar[d] & F\ar[d]\ar[r] & D\ar[r]\ar@{=}[d] & \sh{E}\ar[d] \\
	E'\ar[r] & F'\ar[r] & D\ar[r] & \sh{E'}.
}\]
The fact that the bottom row is a distinguished triangle, with
$E'\in\ct^c\cap\ct^{\leq0}\subseteq\ct^-_c\cap\ct^{\leq0}$ and with
$D\in\ct^-_c\cap\ct^{\leq0}$, tells us that $F'\in\ct^-_c\cap\ct^{\leq0}$.
And since the square on the left is a homotopy
pushout, we may fill in the dotted arrow
in the diagram
\[\xymatrix{
	E\ar[r]\ar[d] & F\ar[d]\ar@/^1pc/[rdd] \\
	E'\ar[r]\ar@/_1pc/[rrd] & F'\ar@{.>}[rd] \\
	& & X
}\]
to give the required factorization
$F\la F'\la X$.
\eprf

Following \cite[Definitions~3.3.1 and 4.1.1]{Neeman99}, recall that if $\ct$ is a triangulated category with small coproducts, a set $\cs$ of objects in $\ct$ is \emph{an $\aleph_1$--perfect class of $\aleph_1$--small objects} if $0\in\cs$ and, given an object $F\in\cs$, an arbitrary set of objects
$\{X_\lambda\st\lambda\in\Lambda\}$ in $\ct$, and any morphism
\begin{equation}\label{eq:coprod}
\xymatrix@C+20pt{
	F\ar[r]^-f &
	\ds\coprod_{\lambda\in\Lambda} X_\lambda\,
}
\end{equation}
then there exists a factorization
\begin{equation}\label{eq:ses}
\xymatrix@C+20pt{
	F\ar[r] &\ds\coprod_{\lambda\in\Lambda'}F_\lambda\ar[r] &
	\ds\coprod_{\lambda\in\Lambda'}X_\lambda\ar@{^{(}->}[r] &
	\ds\coprod_{\lambda\in\Lambda}X_\lambda
}
\end{equation}
with $\Lambda'\subseteq\Lambda$ a countable set, where $F_\lambda$ are all
objects of $\cs$.

We are now ready to prove the main result of this section which clarifies the role of $\ct_c^-$ as a set of generators for $\ct$.

\thm{T79.9}
Let $\ct$ be a weakly approximable triangulated category, and let $\tst\ct$ be a
\tstr\ in the preferred equivalence class. Then the full subcategory $\ct^-_c$
is an $\aleph_1$--perfect class of $\aleph_1$--small objects. Moreover, for every morphism as in \eqref{eq:coprod}, a factorization as in \eqref{eq:ses} can be chosen such that, for every integer $n>0$, all but
finitely many of the $F_\lambda$ lie in $\ct^{\leq-n}$.
\ethm

\prf
Since $\ct$ is weakly approximable, given $F\in\ct^-_c$, the hypotheses of \autoref{L79.7} are all fulfilled. Indeed, the existence of a sequence $E_1\la E_2\la E_3\la\cdots$ as in the assumptions of \autoref{L79.7} follows from \cite[Corollary~2.14]{Neeman17A}. Moreover, by \autoref{T701.9} and \autoref{R1.1}, the given
\tstrs on $\ct$ can be replaced by a compactly generated equivalent one
for which $\ct^{\geq0}$ is
closed under coproducts, achieving hypothesis (i) of \autoref{L79.7}. Hypothesis (ii) of the same lemma follows from \cite[Lemma~2.8]{Neeman17A}. Finally hypothesis (iii) of \autoref{L79.7} is the content of \autoref{L22.3}. 
Thus \autoref{L79.7} can be used and allows us to factor $f$ as
\[\xymatrix@C+20pt{
	F\ar[r]^-\ph &\ds\coprod_{\lambda\in\Lambda'}X_\lambda^{\leq-n_\lambda^{}}\ar[r] &
	\ds\coprod_{\lambda\in\Lambda'}X_\lambda\ar@{^{(}->}[r] &
	\ds\coprod_{\lambda\in\Lambda}X_\lambda.
}\]
Now consider the composite below
\[\xymatrix@C+20pt{
	F\ar[r]^-\ph &\ds\coprod_{\lambda\in\Lambda'}X_\lambda^{\leq-n_\lambda^{}}\ar[r]^-\sim &
	\ds\prod_{\lambda\in\Lambda'}X_\lambda^{\leq-n_\lambda^{}}.
}\]
Hypothesis (iii) of \autoref{L79.7}
tells us that the map from the coproduct to the product
is an isomorphism.
Now for each $\lambda\in\Lambda'$, \autoref{L79.8} permits us
to factor the morphism $F\la X_\lambda^{\leq-n_\lambda}$ as
$F\la F_\lambda\la X_\lambda^{\leq-n_\lambda^{}}$ with
$F_\lambda\in\ct^-_c\cap\ct^{\leq-n_\lambda^{}}$.
Since the coproduct and product agree for the collection
$\{F_\lambda\st\lambda\in\Lambda'\}$, this permits us to (uniquely)
fill in the dotted arrow in the diagram


\[\xymatrix@C+40pt{
	F\ar@{.>}[r]\ar@/-2pc/[dr]\ar@/^2pc/[rr]^\ph &\ds\coprod_{\lambda\in\Lambda'}F_\lambda\ar[r]\ar[d]^\wr &
	\ds\coprod_{\lambda\in\Lambda'}X_\lambda^{\leq-n_\lambda^{}}\ar[d]^\wr \\
	&\ds\prod_{\lambda\in\Lambda'}F_\lambda\ar[r] &
	\ds\prod_{\lambda\in\Lambda'}X_\lambda^{\leq-n_\lambda}.
}\]
The top row of the diagram provides the factorization $\ph$ which delivers
the result. 
\eprf

\subsection{Intrinsic descriptions of the subcategories of $\ct^-_c$ and $\ct^c$}
\label{S2}

We are ready to deal with the main task of this section: the intrinsic description of the subcategories of $\ct^-_c$. Let us start with a preliminary result.

\lem{L2.1}
Let $\ct$ be a weakly approximable triangulated category.  An object $F\in\ct^-_c$ belongs to the thick subcategory $\ct^c\subseteq\ct^-_c$
if and only if, for each object $Q\in\ct^-_c$, there exists an
integer $N=N(Q)>0$ such that $\Hom(F,\sh[n]{Q})=0$ for all $n>N$.
\elem

\prf
Choose an integer $A>0$, a compact generator $G\in\ct$, and
a \tstr\ $\tst\ct$ as in \autoref{dfn:wa}. If $F$ is compact then $F\in\genu G{ }{-B,B}$, for some
integer $B>0$ (see, for example, \cite[Remark~0.15 and Lemma~0.9(ii)]{Neeman17A}). Combining this
with the condition $\Hom(G,\ct^{\leq-A})=0$ of
\autoref{dfn:wa}(i), we have that $\Hom(F,\ct^{\leq-A-B})=0$. Let $Q\in\ct^-_c$ be arbitrary. As $Q\in\ct^-_c\subseteq\ct^-$, there must exist an integer $m>0$
with $Q\in\ct^{\leq m}$. Hence $\sh[n]{Q}\in\ct^{\leq-A-B}$ for all
$n>A+B+m$, and so $\Hom(F,\sh[n]{Q})=0$.

Conversely, suppose we are given an object 
$F\in\ct^-_c$ such that, for each
$Q\in\ct^-_c$, there exists an integer $N=N(Q)>0$ with $\Hom(F,\sh[n]{Q})=0$
for all $n>N$. We need to prove that $F$ is compact.
Because $F$ belongs to $\ct^-_c$ we may choose and fix, for any integer
$m>0$, a distinguished triangle $E_m\la F\la D_m$ with $E_m\in\ct^c$ and with
$D_m\in\ct^{\leq-m}$. Since in this triangle $E_m$ and $F$ both belong to $\ct^-_c$, also $D_m\in\ct^-_c$. Next form, in the category $\ct$, the object
\[
Q=\coprod_{m=1}^\infty\sh[-m]{D_{2m}}.
\]
By the above we have that $\sh[-m]{D_{2m}}\in\ct^{\leq-m}\cap\ct^-_c$,
and \autoref{L22.5} allows us to conclude that
$Q$ belongs to $\ct^-_c$.
By the hypothesis on $F$ there must exist an integer $M>0$ with
$\Hom(F,\sh[n]{Q})=0$ for all $n\geq M$. In particular $\Hom(F,\sh[M]{Q})=0$.
As $D_{2M}$ is a direct summand of
$\sh[M]{Q}=\coprod_{m=1}^\infty\sh[M-m]{D_{2m}}$, we have that $\Hom(F,D_{2M})=0$.
But then in the distinguished triangle $E_{2M}\la F\la D_{2M}$ the map
$F\la D_{2M}$ must vanish, making $F$ a direct summand of the compact
object $E_{2M}$. Therefore $F$ must be compact.
\eprf

We are now ready to prove the main result of this section.

\pro{P2.5}
Lat $\ct$ be a weakly approximable triangulated category. 
\be
\item
The full subcategory $\ct^c\subseteq\ct^-_c$ has for objects
all those $F\in\ct^-_c$ such that, for any object
$Q\in\ct^-_c$, there exists an integer $N>0$ with
$\Hom(F,\sh[n]{Q})=0$, for all $n>N$.
\item
For any classical
generator $G\in\ct^c$, 
the full subcategory $\ct^b_c\subseteq\ct^-_c$
is given by the formula
\[
\ct^{b}_c=\bigcup_{n=1}^\infty \big(G(-\infty,-n]\big)^\perp
  \]
where the perpendicular is taken in $\ct^-_c$. 
\item
For any classical
generator $G\in\ct^c$, 
the full subcategory $\ct^{c,b}\subseteq\ct^c$
is given by the same formula
\[
\ct^{c,b}=\bigcup_{n=1}^\infty \big(G(-\infty,-n]\big)^\perp
  \]
where this time the perpendicular is taken in $\ct^c$. 
\ee
\epro

\prf
The characterization in (i) of $\ct^c\subseteq\ct^-_c$ was proved
in \autoref{L2.1}. Thus we only need to prove (ii) and
(iii).

Recall
that a classical generator $G\in\cs=\ct^c$ is the
same as a compact generator for $\ct$.
If $\ca=G(-\infty,0]$, then in the category
$\ct$ we have the equalities
\[
\ca^\perp=\Coprod(\ca)^\perp=\ct^{\geq1}_G\ ,
\]  
where the first equality is by \autoref{R701.1}(ii)
and the second is by combining 
\autoref{R1.1} with \autoref{T701.9}.
Thus, still in the category $\ct$, we have
\[
\bigcup_{n=1}^\infty \big(G(-\infty,-n]\big)^\perp\eq
\bigcup_{n=1}^\infty\ct_G^{\geq-n+1}\eq
\ct^+\ ,
\]
and intersecting with $\ct^-_c$ gives (ii)
while intersecting with $\ct^c$ yields (iii).
\eprf

\section{The subcategories of $\ct^-$, and a recipe for
$\ct^b$ as a subcategory of $\ct^+$}
\label{S4}

This section is about the inclusions in
the diagram \eqref{eq:incl} of the subcategories of $\ct^-$,
and about describing $\ct^b$ intrinsically as
a subcategory of $\ct^+$. This will be elaborated in \autoref{subsec:inl2}.
But first, in \autoref{S3},
we discuss the intrinsic \tstrs in $\ct^?$, for $?=-,+,b$.

\subsection{The preferred equivalence class of \tstrs
	on $\ct^-$, $\ct^+$ and $\ct^b$}
\label{S3}

For the rest of this section, we assume $\ct$ to be a triangulated category with small coproducts and a single compact generator $G$. 
We know that $\ct$ is endowed with a preferred equivalence class of \tstr s. And, by the definition of
the ``preferred'' $\ct^-$, $\ct^+$ and $\ct^b$, any \tstrs on $\ct$, in the preferred equivalence class, restricts to a \tstrs on each of the three preferred subcategories $\ct^-$, $\ct^+$ and $\ct^b$. Furthermore, the restrictions of any two equivalent \tstr{s} are equivalent. This defines
for us preferred equivalence classes of \tstr{s} on each of full triangulated subcategories mentioned above. But this
definition is in terms of the embedding into $\ct$.

The aim of this section will be to show that, the preferred equivalence class of
\tstr s, on each of $\ct^-$, $\ct^+$ and $\ct^b$, has an intrinsic
description---by which we mean a description
which doesn't mention the embedding into $\ct$.

We begin with the following.

\dfn{D3.1}
Let $\cs$ be a triangulated category, and consider the class
$\Cp(\cs)$
of all full subcategories $P\subseteq\ct$ satisfying $\sh{P}\subseteq P$.
\be
\item
Two elements $P,Q\in\Cp(\cs)$ are \emph{equivalent}
if there exists an integer $A>0$ with $\sh[A]{P}\subseteq Q\subseteq\sh[-A]{P}$.
\item
Given two equivalence classes $[P]$ and $[Q]$ of elements of $\Cp(\cs)$,
we set $[P]\leq[Q]$ if, for a choice of representatives
$P\in[P]$ and $Q\in[Q]$, there exists an integer $A>0$ with 
$\sh[A]{P}\subseteq Q$.
\ee
\edfn

In the sequel, we will use the shorthand $\Cp:=\Cp(\cs)$.

\rmk{R3.3}
If we start with two \tstrs $\tstv\cs 1$ and $\tstv\cs 2$ on a triangulated category $\cs$ and we set $P:=\cs_1^{\leq 0}$ and $Q:=\cs_2^{\leq 0}$, then $P$ and $Q$ are equivalent as in \autoref{D3.1} if and only if the two \tstr s are equivalent in the sense discussed in \autoref{R1.1}. That is, there is a positive integer $N>0$ such that $\cs_1^{\leq -N}\subseteq\cs_2^{\leq 0}\subseteq\cs_2^{\leq N}$. The main point of \autoref{D3.1} is that we extend this equivalence
relation beyond \tstr s and, more importantly, we introduce
the partial order of \autoref{D3.1}(ii).
\ermk

Let us consider the bounded above case.

\pro{P3.5}
Let $\ct$ be a triangulated category with coproducts
and a single compact generator, and
put $\cs=\ct^-$, where $\ct^-$ is understood
to mean with respect to the preferred
equivalence class of \tstr{s} on $\ct$. Then a \tstr\ $\tst\cs$
belongs to the preferred equivalence class if and only if the following two conditions below are both satisfied:
\be
\item
$\cs=\bigcup_{m=1}^\infty\cs^{\leq m}$.
\item
The equivalence class $\left[\cs^{\leq0}\right]$ is the unique 
minimal one with respect to the partial odering $\leq$ for aisles of \tstr s satisfying (i).
\ee
\epro

\prf
If $\tst\ct$ is a \tstr\ on $\ct$ in the preferred equivalence
class, then by definition
\[
\cs=\ct^-=\cup_{m=1}^\infty\ct^{\leq m}.
\]
The restriction of $\tst\ct$ to a \tstr\ $\tst\cs$ on $\cs$ is defined
by the formulas $\cs^{\leq 0}=\cs\cap\ct^{\leq 0}$ and
$\cs^{\geq0}=\cs\cap\ct^{\geq0}$, which in the case of
$\cs=\ct^-$ yields $\cs^{\leq0}=\ct^{\leq0}$ and therefore 
$\cs^{\leq m}=\ct^{\leq m}$. Thus $\cs=\cup_{m=1}^\infty\cs^{\leq m}$ holds,
and the restricted \tstr\ on $\cs$
satisfies (i) in the statement.

What needs proof is that any \tstrs $\tst\cs$ on $\cs$, which satisfies
(i), must have that $\left[\cs^{\leq0}\right]$ is bigger than
$\left[\ct^{\leq0}\right]$
in the ordering $\leq$ of \autoref{D3.1}(ii). For this purpose, pick a compact
generator $G\in\ct$. As $G\in\ct_G^{\leq0}\subset\ct^-=\cs$,
and since
$\cs=\bigcup_{m=1}^\infty\cs^{\leq m}$, we must
have that $G$ belongs to $\cs^{\leq m}$ for some
integer $m>0$. But then $G(-\infty,0]\subset\cs^{\leq m}$,
and as the category $\cs^{\leq m}$ is the
aisle of a \tstrs it is closed under
extensions and coproducts which exist in $\ct^-$, giving
\[
\cs^{\leq m}\supseteq\Coprod\big(G(-\infty,0]\big)=\ct_G^{\leq0}.
\]
This proves the inequality
$\big[\ct_G^{\leq0}\big]\leq\big[\cs^{\leq0}\big]$ and, as
$\tstv\ct G$ is a representative of the preferred
equivalence class of \tstr{s}, the proof is complete.
\eprf

Let $\ct$ be a triangulated category with coproducts
and a single compact generator $G$.
Then the functor $\Hom_\ct^{}(G,-)$ is a homological
functor $\Hom_\ct^{}(G,-)\colon\ct\la\ab$ respecting coproducts. Because
$\qq/\zz$ is injective in $\ab$, the functor
$\fH\colon\ct\op\la\ab$ given by the formula
\[
\fH(-):=\Hom_{\ab}\big(\Hom_\ct^{}(G,-)\,,\,\qq/\zz\big)
\]
is homological and takes coproducts in $\ct$ to products
in $\ab$. Hence, by
\cite[Theorem~3.1]{Neeman96}, the functor $\fH$ is
representable, meaning that there exists an object $\bc(G)\in\ct$
and a natural isomorphism
$\fH(-)\cong\Hom_\ct\big(-,\bc(G)\big)$. The reason for the notation is that, in honor of  the paper \cite{Brown-Comenetz76} where the
construction was first used, the object $\bc(G)$ is called the \emph{Brown-Comenetz dual}
of $G$.

\lem{L3.9}
Let $\ct$ be a triangulated category with coproducts
and a single compact generator $G$, and assume
that $\Hom(G,\T^n G)=0$ for $n\gg0$. Then the object
$\bc(G)$ belongs to
$\ct^+$, and $\Hom\big(X,\bc(G)\big)=0$ if and only if
$\Hom(G,X)=0$.
\elem

\prf
For any $X\in\ct$ we have that $\Hom\big(X,\bc(G)\big)$
is canonically isomorphic to
\[
\Hom_{\ab}\big(\Hom_\ct^{}(G,X)\,,\,\qq/\zz\big).
\]
As $\qq/\zz$ is an injective cogenerator of the category of abelian groups, this vanishes if and only if $\Hom_\ct(G,X)=0$.
This proves the second assertion of the lemma.

To prove the first one, we apply the second assertion
to $X:=\sh[n]{G}$. By hypothesis there exists an
integer $B>0$ such that $\Hom(G,\sh[n]{G})$ vanishes for all $n>B$.
Hence $\Hom\big(\sh[n]{G},\bc(G)\big)$ vanishes for
all $n>B$, that is $^\perp\bc(G)$ contains $G(-\infty,-B-1]$.  But then
\[
{}^\perp\bc(G)\supseteq\Coprod\big( G(-\infty,-B-1]\big)=\ct_G^{\leq-B-1}.
\]
From this we deduce that $\bc(G)$ must lie in $\ct_G^{\geq -B}\subseteq\ct^+$.
\eprf

\ntn{Notforthis section}
Let $\cs$ be any triangulated category, and let
$\Cp=\Cp(\cs)$ be as in \autoref{D3.1}.
For any object $I\in\cs$ we can form 
$I[0,\infty)\subset\cs$ as in \autoref{N314.3}(i),
  and then let $P(I)\in\Cp(\cs)$ be given
  by the formula $P(I)={^\perp I[0,\infty)}$.
Spelling this out, we obtain
\[
P(I):={^\perp I[0,\infty)}=\bigcap_{m=0}^\infty{^\perp{\big(\sh[-m]{I}\big)}}\in\Cp.
\]  
We adopt this notation for the rest of this section.  
\entn

\exm{ex:bc}
Assume $\ct$ is a weakly approximable triangulated category
and 
$G\in\ct$ is a compact generator.
Let us work out what $P(I)$ is, in the case $\cs=\ct^+$ and
$I=\bc(G)$.
The fact that $I$ belongs to $\cs=\ct^+$ was proved in \autoref{L3.9} (note that $\Hom(G,\T^n G)=0$ for $n\gg0$ by (i) of \autoref{dfn:wa}).
On the other hand $\Hom(X,\sh[-m]{I})=0$ if and only if $\Hom(\sh[m]{X},I)=0$. By \autoref{L3.9}, the latter equality is equivalent to  $\Hom(G,\sh[m]{X})=0$. Thus with $I=\bc(G)$ we get the explicit description $P(I)=G[0,\infty)^\perp$. By \cite[Lemma~3.9(iv)]{Burke-Neeman-Pauwels18} we deduce that there exists
an integer $A>0$ with $P(I)\subseteq\cs\cap\ct^{\leq A}=\cs^{\leq A}$.
\eexm

This explicit description is useful in treating the bounded below case.

\pro{P3.11}
Let $\ct$ be a weakly approximable triangulated category, and put
$\cs=\ct^+$. Then, with $P(I)$ as
in \autoref{Notforthis section},
the collection of equivalence classes
in $\Cp(\cs)$ of the form $\big[P(I)\big]$
has a unique minimal member, in the partial order of \autoref{D3.1}. Moreover, a \tstrs $\tst\cs$ belongs to the preferred
equivalence class if and only if $\big[\cs^{\leq0}\big]$
is in the equivalence class of the minimal $\big[P(I)\big]$.
\epro

\prf
Let $\tst\cs$ be the intersection with $\cs$ of a \tstrs $\tst\ct$
in the preferred equivalence class. We need to show that
$\big[\cs^{\leq0}\big]\leq\big[P(J)\big]$ for every $J\in\cs$, and
that there exists an $I\in\cs$ with
$\big[P(I)\big]\leq\big[\cs^{\leq0}\big]$.

For any $J$ belonging to $\cs$ we may assume, after shifting,
that $J\in\cs^{\geq1}$. Therefore $J[0,\infty)\subset\cs^{\geq1}$,
  and $P(J)={^\perp J[0,\infty)}\supset\cs^{\leq0}$. Hence,
    in the partial order of \autoref{D3.1}, we have
\[
\big[\cs^{\leq0}\big]\leq\big[P(J)\big]
\]
for every $J\in\cs$.

It remains to observe that, if $G$ is a compact generator for $\ct$, then $I:=\bc(G)\in\cs$ and
\[
\big[P(I)\big]\leq\big[\cs^{\leq0}\big]\ .
\]
by the discussion in \autoref{ex:bc}.
\eprf

We are now ready to treat the bounded case.

\pro{P3.13}
Suppose $\ct$ is a weakly approximable triangulated category and put
$\cs=\ct^b$. Then, with $P(I)$ as
in \autoref{Notforthis section},
the collection of equivalence classes
in $\Cp(\cs)$ of the form $\big[P(I)\big]$
has a unique minimal member, in the partial order of \autoref{D3.1}.
Moreover, a \tstrs $\tst\cs$ belongs to the preferred
equivalence class if and only if $\big[\cs^{\leq0}\big]$
is in the equivalence class of the minimal $\big[P(I)\big]$.
\epro

\prf
Let $\tst\cs$ be the intersection with $\cs$ of a \tstrs $\tst\ct$
in the preferred equivalence class.
As in the previous proof, for any object $J$ belonging to $\cs$ we may
assume, after shifting,
that $J\in\cs^{\geq1}$. Therefore
$J[0,\infty)\subseteq\cs^{\geq1}$, hence
  $\cs^{\leq0}\subseteq{^\perp J[0,\infty)}= P(J)$,
    and this gives
\[
\big[\cs^{\leq0}\big]\leq\big[P(J)\big]
\]
for every $J$. Again, to conclude the proof, we need to exhibit an object $I\in\cs=\ct^b$ with
\[
\big[P(I)\big]\leq\big[\cs^{\leq0}\big]\ .
\]
Here we take
\[
I:=\bc(G)^{\leq0}\oplus\prod_{\ell=1}^\infty
\sh[\ell]{\left(\bc(G)^{\leq\ell}\right)^{\geq\ell}}\ .
\]

First we need to show that $I\in\ct^b$. As $\bc(G)$ belongs to $\ct^+$ it follows that
$\bc(G)^{\leq0}$ must belong to $\ct^b$.
And as all the other terms,
in the product defining $I$,
belong to $\ct^\heartsuit=\ct^{\leq0}\cap\ct^{\geq0}$,
\autoref{L22.2} guarantees that the product is in $\cs=\ct^b$.

Finally we claim that, with $A>0$ as in \autoref{ex:bc}, we have
that $P(I)\subseteq\cs^{\leq A}$. Indeed, suppose $X\in\cs=\ct^b$ lies outside $\cs^{\leq A}=\ct^b\cap\ct^{\leq A}$.
Then it also lies outside $\ct^+\cap\ct^{\leq A}$, and by setting
$J=\bc(G)$
and applying \autoref{ex:bc} we deduce that $X$ cannot belong to
$P(J)$. Therefore there
must exist in
$\ct^+$ a nonzero map $X\la\sh[-m]{J}$ for some $m\geq0$.
Since $X$ belongs to $\ct^b$, this map
must factor as $X\la \sh[-m]{\big(J^{\leq\ell}\big)}\la\sh[-m]{J}$,
for some $\ell\geq0$. As the composite is
nonzero, the map  $X\la \sh[-m]{\big(J^{\leq\ell}\big)}$ cannot
be trivial. Let us choose $\ell\geq0$ to be the smallest
$\ell$ for which such a factorization exists.

If $\ell=0$
this gives a nonzero map
$X\la\sh[-m]{J^{\leq0}}$. But $J^{\leq0}$ is a direct summand of $I$,
and hence we have a nonzero map $X\la\sh[-m]{I}$.

If $\ell>0$ consider the composite
$X\la \sh[-m]{\big(J^{\leq\ell}\big)}\la
\sh[-m]{\Big(\big(J^{\leq\ell}\big)^{\geq\ell}\Big)}$.
This composite cannot vanish by the minimality of $\ell$. Hence,
in this case, we get a nonzero map
$X\la \sh[-m]{\Big(\big(J^{\leq\ell}\big)^{\geq\ell}\Big)}$.
As $\sh[\ell]{\Big(\big(J^{\leq\ell}\big)^{\geq\ell}\Big)}$ is a
direct summand of $I$, this provides a nonzero map
$X\la\sh[-m-\ell]{I}$.

Either way, $X$ does not belong to $P(I)$. Since $X$ was an arbitrary
object outside $\cs^{\leq A}$, we deduce that $P(I)\subseteq\cs^{\leq A}$.
\eprf

\subsection{The intrinsic inclusions}\label{subsec:inl2}

We are now ready to give recipes for
the subcategories of $\ct^-$. Let us start with the case of compact objects.

\lem{L4.1}
Let $\ct$ be a weakly approximable triangulated category. Then $(\ct^-)^c=\ct^c$.
\elem

\prf
Since $\ct^c$ is contained in $\ct^-$ and since, by \autoref{L27.9},
the inclusion $\ct^-\hookrightarrow\ct$ preserves coproducts, we have the inclusion
$\ct^c\subset(\ct^-)^c$. We need to prove the reverse inclusion.

For this choose a \tstr\ $\tst\ct$ in the preferred equivalence class,
such that $\ct^{\geq0}$ is closed under coproducts and hence the truncation functor $(-)^{\leq0}$
commutes with coproducts. This can be achieved due to the discussion in \autoref{S701} (see, in particular, \autoref{rmk:propgen}).

Let $c\in(\ct^-)^c$ . Because $c$ belongs to $\ct^-$, we may (after shifting) assume
$c\in\ct^{\leq0}$. And now let $\{t_\lambda\st\lambda\in\Lambda\}$ be
an arbitrary collection of objects in $\ct$ and suppose we are given a
map
\[\xymatrix@C+30pt{
c\ar[r]^-{\psi}&\ds\coprod_{\lambda\in\Lambda}t_\lambda^{}\ .
}\]
As $c$ belongs to $\ct^{\leq0}$ the map $\psi$ factors through
$\left(\coprod_{\lambda\in\Lambda}t_\lambda\right)^{\leq0}$, but as
the functor $(-)^{\leq0}$ commutes with
coproducts we can write this factorization as
\[\xymatrix@C+30pt{
c\ar[r]^-\ph &\ds\coprod_{\lambda\in\Lambda}t_\lambda^{\leq0}
\ar[r] &\ds\coprod_{\lambda\in\Lambda}t_\lambda^{}
}\]
As $c\in(\ct^-)^c$, the map $\ph$ must factor
through a finite subcoproduct. But then the map
$\psi$ also factors through a finite subcoproduct.
\eprf

We conclude the task of this section, proving the following, additional parts of \autoref{thm:main1}.

\cor{cor:iota694}
If $\ct$ is a weakly approximable triangulated category, then the inclusions
in the diagram \eqref{eq:incl} of the subactegories
$\ct^b\subset\ct^-$, $\ct^b\subset\ct^+$ and $\ct^-_c\subset\ct^-$ are intrinsic.
\ecor

\prf
\autoref{L4.1} above provides a recipe for $\ct^c$
as a subcategory of $\ct^-$, while
in \autoref{P3.5} we gave a recipe for
the preferred equivalence class of \tstr s on $\ct^-$,
and in \autoref{P3.11} we gave a recipe for the
preferred equivalence class of \tstr{s} on $\ct^+$.
The proof of the current corollary combines the information.

In terms of the preferred equivalence class of \tstr s, we can describe $\ct^b$ as a subcategory of $\ct^-$ by
the recipe
\[
\ct^b=\ct^-\cap\ct^+=\bigcup_{m=1}^\infty(\ct^-)^{\geq-m}
\]
for any \tstr\ in the preferred equivalence class.
The term on the far right
has an intrinsic description in $\ct^-$ by \autoref{P3.5},
hence $\ct^b\subseteq\ct^-$ is intrinsic.  The formula
\[
\ct^b=\ct^-\cap\ct^+=\bigcup_{m=1}^\infty(\ct^+)^{\leq m}
\]
gives a
recipe for $\ct^b$ as a subcategory of $\ct^+$,
and in view of \autoref{P3.11} the term on the far right
is intrinsic in $\ct^+$.

The subcategory
$\ct^-_c\subseteq\ct^-$ is given by the formula
\[
\ct^-_c=\bigcap_{m=1}^\infty\ct^c\ast\ct^{\leq-m}=
\bigcap_{m=1}^\infty(\ct^-)^c\ast(\ct^-)^{\leq-m}\ ,
\]
where the second equality is by \autoref{L4.1}, where
$(\ct^-)^c$ is obviously intrinsic in $\ct^-$, and where
$(\ct^-)^{\leq-m}$ has (up to equivalence) an intrinsic
description
in $\ct^-$ by \autoref{P3.5}. Hence the term on
the far right gives 
an intrinsic recipe for $\ct^-_c\subset\ct^-$.
\eprf

\section{Generalities about pseudocompact and strongly pseudocompact objects}
\label{S5}

Using the discussion in the previous sections, in this one we introduce and study two new classes of objects: the pseudocompact and the strongly pseudocompact ones. Later in the paper, we will specialize to the case of weakly approximable triangulated categories.

\subsection{Pseudocompact objects}\label{subsec:pseudocomp}

Let us first introduce a slight generalization of the notion of compact object. We always assume that $\cs$ is a triangulated category and that $P$ is an element of $\Cp:=\Cp(\cs)$, with $\Cp(\cs)$ as in \autoref{D3.1}.

\dfn{D5.1}
\be
\item
An object $c\in\cs$ is \emph{$P$--pseudocompact}
if for any set
$\{s_\lambda\st\lambda\in\Lambda\}$ of objects,
which all belong to $P^\perp$ and whose coproduct
exists in $\cs$,  any morphism
\[\xymatrix@C+30pt{
c\ar[r] &\ds\coprod_{\lambda\in\Lambda}s_\lambda^{}
}\]
must factor through a finite subcoproduct.
\item
Given an equivalence class $[P]$ of elements in $\Cp$,
an object is \emph{$[P]$--pseudocompact} if
it is $Q$--pseudocompact for every $Q\in[P]$.
\ee
\edfn

The full subcategory of all $P$--pseudocompact objects will be denoted
$\cs^{pc}_P$. And the full subcategory of all $[P]$--pseudocompact objects
will be written $\cs^{pc}_{[P]}$.

\exm{R5.13.9}
There are some obvious examples of $P$-pseudocompact objects,
for special $P$s, and here we discuss a
particularly simple one.
Let $\tst\cs$ be a \tstr\ on $\cs$
such that $\cs^{\geq0}$ is closed in $\cs$ under coproducts. Then
every object in $\cs^{\leq 0}$ is $\cs^{\leq0}$--pseudocompact.
To see this let $c\in(\cs^{\leq 0})$ be any object. 
Assume
$\{s_\lambda\st\lambda\in\Lambda\}$ is a collection
of objects,
which all belong to $(\cs^{\leq 0})^\perp=\cs^{\geq1}$, 
and whose coproduct
exists in $\cs$.  By assumption $\cs^{\geq1}$ is closed
in $\cs$ under coproducts, hence a morphism
$c\la\ds\coprod_{\lambda\in\Lambda}s_\lambda$
is a map from $c\in\cs^{\leq0}$ to an object in $\cs^{\geq1}$ and
must vanish.
\eexm

The first result is easy.

\lem{L5.11}
Let $\cs$ be a triangulated category, and let $P\in\Cp$ be any
element. Then the subcategory
$\cs^{pc}_P$ is closed under direct summands and extensions. To
say it in symbols:
\[
\Smr(\cs^{pc}_P)=\cs^{pc}_P\qquad\text{ and }
\qquad\cs^{pc}_P*\cs^{pc}_P=\cs^{pc}_P\ .
\]
\elem

\prf
The assertion about direct summands is clear.
To prove the equality $\cs^{pc}_P*\cs^{pc}_P=\cs^{pc}_P$
let $A\la B\la C$ be a distinguished triangle in $\cs$, with
$A$ and $C$ both $P$--pseudocompact. Let
$\{s_\lambda\st\lambda\in\Lambda\}$ be a set of objects
in $P^\perp\subseteq\cs$, whose
coproduct exists in $\cs$.  Choose any morphism
\[\xymatrix@C+30pt{
B\ar[r]^-\ph &\ds\coprod_{\lambda\in\Lambda}s_\lambda^{}\ .
}\]
Because $A$ is $P$--pseudocompact the composite
\[\xymatrix@C+30pt{
A\ar[r] &B\ar[r]^-\ph &\ds\coprod_{\lambda\in\Lambda}s_\lambda^{}
}\]
must factor through a finite subcoproduct, which means that
the natural projection
to a subset $\Lambda'\subseteq\Lambda$ with a finite complement
gives a vanishing composite
\[\xymatrix@C+30pt{
A\ar[r] &B\ar[r]^-\ph &\ds\coprod_{\lambda\in\Lambda}s_\lambda^{}
\ar[r] &\ds\coprod_{\lambda\in\Lambda'}s_\lambda^{}\ .
}\]
Hence there is a factorization
\[\xymatrix@C+30pt{
B\ar[r] &C
\ar[r]^-\psi &\ds\coprod_{\lambda\in\Lambda'}s_\lambda^{}\ ,
}\]
and the $P$--pseudocompactness of $C$ guarantees that $\psi$
factors through a finite subcoproduct, meaning that by passing
to a subset $\Lambda''\subseteq\Lambda'$ with finite
complement we obtain the vanishing of
\[\xymatrix@C+30pt{
B\ar[r] &C
\ar[r] &\ds\coprod_{\lambda\in\Lambda'}s_\lambda^{}
\ar[r] &\ds\coprod_{\lambda\in\Lambda''}s_\lambda^{}\ .
}\]
Thus $\ph$ factors through the coproduct over $\Lambda\setminus\Lambda''$,
which is a finite set.
\eprf

\cor{C5.13}
Let $\cs$ be a triangulated category, and let $[P]\subseteq\Cp$
be an equivalence class of elements in $\Cp$. Then $\cs^{pc}_{[P]}$
is a thick subcategory of $\cs$.
\ecor

\prf
The closure under direct summands and extensions follows from \autoref{L5.11}.
The closure under suspension is because we
have an equality of equivalence classes $[P]=[\sh{P}]$.
\eprf

\exm{R5.13}
Let $\cs$ be a triangulated category, and let $\tst\cs$ be a \tstr\
such that $\cs^{\geq0}$ is closed in $\cs$ under coproducts.
By \autoref{R5.13.9}
every object of $\cs^{\leq 0}$ is $\cs^{\leq0}$--pseudocompact, and
clearly every object in $\cs^c$ is also $\cs^{\leq0}$--pseudocompact.
By \autoref{L5.11} every object of $\cs^c*\cs^{\leq0}$ is
$\cs^{\leq0}$--pseudocompact.
It follows that every object in
\[
\cs^-_{c,[\cs^{\leq0}]}:=\bigcap_{m=0}^\infty\cs^c*\cs^{\leq-m}
\]
is $\cs^{\leq-m}$--pseudocompact for every $m>0$, that is
$\cs^-_{c,[\cs^{\leq0}]}\subseteq\cs^{pc}_{[\cs^{\leq0}]}$ .

Now suppose that
$\ct$ is a
weakly approximable triangulated category, and $\tst\ct$
is a \tstrs in the preferred equivalence class. Then
we get the easy equality $\ct^-_{c,[\ct^{\leq0}]}=\ct^-_c$. Hence we get the inclusion
\[
\ct^-_c\subseteq\ct^{pc}_{[\ct^{\leq0}]}\ .
\]
We will return to this in \autoref{S6}.
\eexm

\lem{L5.14}
Let $\cs$ be a triangulated category, let $\tst\cs$ be a \tstr\ on
$\cs$, and assume $\cs^{\geq0}$ is closed in $\cs$ under coproducts.
If $C$ is a $\cs^{\leq0}$--pseudocompact object in $\cs$ then
so is $C^{\geq 2}$.
\elem

\prf
We have a distinguished triangle $C^{\leq1}\la C\la C^{\geq2}\la\sh{C^{\leq1}}$,
with $C$ in $\cs^{pc}_{\cs^{\leq0}}$ by hypothesis. Now
$\sh{C^{\leq1}}\in\cs^{\leq0}\subseteq\cs^{pc}_{\cs^{\leq0}}$,
where the inclusion in by \autoref{R5.13.9},
and the result now follows from \autoref{L5.11}.
\eprf

\subsection{Strongly pseudocompact objects}
\label{S67}

In this section we introduce the slightly
more restrictive notion of strongly pseudocompact objects.
The authors do not know whether all pseudocompact
objects are strongly pseudocompact as well.

Let us start from the definition and assume, as in the previous section, that $\cs$ is a triangulated category and that $\Cp=\Cp(\cs)$ is as in \autoref{D3.1}. 

\dfn{D67.1}
Given $P\in\Cp$, an object $c\in\cs$ is
\emph{$[P]$--strongly pseudocompact} if the following two conditions are satisfied:
\be
\item
$c$ is $[P]$--pseudocompact.
\item
Assume we are given any \tstr\ $\tst\cs$ on the category $\cs$,
such that $\big[\cs^{\leq0}\big]=[P]$. Suppose
that $\cs^{\geq0}$ is closed in $\cs$ under
coproducts and the heart $\cs^\heartsuit$ satisfies [AB5].
Then the functor $\Hom(c,-)$ commutes with filtered colimits in
$\cs^\heartsuit$.
\ee
\edfn

The collection of all $[P]$--strongly pseudocompact objects
will be written $\cs^{spc}_{[P]}$.

Let us discuss a slightly technical construction that will be used in a short while. Fix a \tstrs $\tst\cs$ on $\cs$ and assume that we are given in $\cs^\heartsuit$ a complex
\begin{equation}\label{eq:complex2}
\cf^\bullet:=\{\xymatrix@C+0pt{
 \cdots \ar[r] & \cf^{-4}\ar[r] & \cf^{-3}\ar[r] &
 \cf^{-2}\ar[r] &\cf^{-1}\ar[r] &\cf^{0}
}\}
\end{equation}
which is exact.

\exm{ex:compl}
We have already encountered an example of such a complex. Indeed, let $\ca$ be an abelian category satisfying [AB4] and let $I$ be a filtered category and let $\fF\colon I\la\ca$ be a functor. Let $\cf^\bullet$ be the standard cochain complex for computing
$\colimj F$, which we introduced in the proof of \autoref{L69.15}.
Let us briefly recall that $\cf^\bullet$ has the form
\[
	\xymatrix@C+0pt{
		\cdots \ar[r] & \cf^{-3}\ar[r] &
		\cf^{-2}\ar[r] &\cf^{-1}\ar[r] &\cf^{0}\ar[r] &0\ ,
	}
\]
the $(-j)\mth$ cohomology is $\colimj \fF$, and $\cf^{-n}$ is the
coproduct over sequences of composable morphisms in $I$
\[\xymatrix@C+0pt{
	i_0^{}\ar[r] & i_1^{}\ar[r] &\cdots\ar[r]&i_{n-1}^{}\ar[r] &i_n^{}
}\]
of $\fF(i_0^{})$. If $\ca$ satisfies [AB5], then
this cochain complex is a resolution of $\colim\, \fF$.
\eexm

Now let $\ck^n$ be the cokernel of the map $\cf^{n-1}\la\cf^n$.
Thus $\ck^0$ is the cokernel of $\cf^{-1}\la\cf^0$.
The acyclicity of the complex tells that $\ck^{-1}$ is the kernel
of $\cf^0\la\ck^0$ and, for all
$n<-1$, the object $\ck^n$ is
the kernel of the map $\cf^{n+1}\la\cf^{n+2}$. For $i\in\nn_{>0}$, we inductively construct the objects $C_i\in\cs$ with natural maps $C_i\la C_{i+1}\la\sh[i+1]{\cf^{-i-1}}$ as follows:
\begin{itemize}
\item Define $C_1$ and the morphisms $\cf^0\la C_1\la\sh{\cf^{-1}}$
by completing
$\cf^{-1}\la\cf^{0}$ to a distinguished triangle
$\cf^{-1}\la\cf^{0}\la C_1\la\sh{\cf^{-1}}$.
\end{itemize}
By construction we have that
$C_1^{\leq-1}=\sh{\ck^{-2}}$ and $C_1^{\geq0}=\ck^0$.
Now let us continue our inductive constructions as follows:
\begin{itemize}
\setcounter{enumi}{\value{enumiv}}
\item
Assume the object $C_n$ has been defined, and satisfies $C_n^{\geq0}=\ck^0$
and $C_n^{\leq-1}=\sh[n]{\ck^{-n-1}}$. Then we take the composite
\[\xymatrix{
\varphi\colon\sh[n]{\cf^{-n-1}}\ar[r] &\sh[n]{\ck^{-n-1}}\ar@{=}[r] & C_n^{\leq-1}
\ar[r] & C_n
}\]
And the morphisms $C_n\la C_{n+1}\la\sh[n+1]{\cf^{-n-1}}$ are
formed by completing the above composite to
a distinguished triangle
\begin{equation}\label{eq:mor1}
\xymatrix{
\sh[n]{\cf^{-n-1}}
\ar[r]^-{\varphi} & C_n
\ar[r] & C_{n+1}\ar[r] &\sh[n+1]{\cf^{-n-1}}
}
\end{equation}
\end{itemize}
By completing the composable morphisms 
$\sh[n]{\cf^{-n-1}}\la\sh[n]{\ck^{-n-1}}\la C_n$ to an octahedron the reader can easily
show that $C_{n+1}^{\geq0}=\ck^0$ and $C_{n+1}^{\leq-1}=\sh[n+1]{\ck^{-n-2}}$.

\lem{L67.7}
In the setting above, let $\cf^\bullet$ be a complex as in \eqref{eq:complex2} and let $X\in\cs$ be such that, for all $i\in\zz$,
$\Hom(\sh[i]{X},-)$ takes $\cf^\bullet$ to an exact sequence.
Then any morphism $X\la C_n$ must factor through the morphism
$\sh[n]{\cf^{-n-1}}\oplus\cf^0\la C_n$, where the morphism $\cf^0\la C_n$ is
the composite
\[\xymatrix@C+0pt{
 \cf^{0}\ar[r] & C_1\ar[r] & C_2\ar[r] &\cdots
 C_{n-1}\ar[r] & C_n\ ,
}\]
while the morphism $\sh[n]{\cf^{-n-1}}\la C_n$ is the composite
\[\xymatrix{
\sh[n]{\cf^{-n-1}}\ar[r] &\sh[n]{\ck^{-n-1}}\ar@{=}[r] & C_n^{\leq-1}
\ar[r] & C_n\ .
}\]
\elem

\prf
We prove the result by induction on $n>0$, starting with
$n=1$. Suppose we are given a morphism $f\colon X\la C_1$. The distinguished triangle
$C_1\la \sh{\cf^{-1}}\la\sh{\cf^0}$ tells us that the composite
$X\la C_1\la\sh{\cf^{-1}}$ gives a morphism $X\la\sh{\cf^{-1}}$ such that
the composite $X\la\sh{\cf^{-1}}\la\sh{\cf^0}$ vanishes.

The assumption in the statement implies that this morphism factors
as $X\la\sh{\cf^{-2}}\la\sh{\cf^{-1}}$. This gives us a morphism
$X\la\sh{\cf^{-2}}$, and we may now form the composite
\[
g\colon X\la\sh{\cf^{-2}}\la\sh{\ck^{-2}}=C_1^{\leq-1}\la C_1.
\]
The construction guarantees that
$(f-g)\colon X\la C_1$ is annihilated by the morphism $C_1\la\sh{\cf^{-1}}$,
and, by the distinguished triangle $\cf^0\la C_1\la\sh{\cf^{-1}}$, we get that
$(f-g)$ must factor as $X\la\cf^0\la C_1$. This proves
the case $n=1$.

Now assume the claim has been proved for $n>0$, and let $f\colon X\la C_{n+1}$
be a morphism. The distinguished triangle $C_{n+1}\la\sh[n+1]{\cf^{-n-1}}\la\sh{C_n}$ in \eqref{eq:mor1} gives the vanishing of the
composite
\[\xymatrix{
X\ar[r]^-f & C_{n+1}\ar[r] & \sh[n+1]{\cf^{-n-1}}\ar[r] &\sh{C_n}
}\]
and hence certainly also the vanishing of the longer composite
\[\xymatrix{
X\ar[r]^-f & C_{n+1}\ar[r] & \sh[n+1]{\cf^{-n-1}}\ar[r] &\sh{C_n}\ar[r] 
&  \sh[n+1]{\cf^{-n}},
}\]
where the morphism $C_n\la \sh[n]{\cf^{-n}}$ comes from the distinguished triangle
defining $C_n$. Thus we have
a vanishing composite $X\la\sh[n+1]{\cf^{-n-1}}\la\sh[n+1]{\cf^{-n}}$, and the
assumption in the statement allows us to factor the morphism
$X\la\sh[n+1]{\cf^{-n-1}}$
as $X\la\sh[n+1]{\cf^{-n-2}}\la \sh[n+1]{\cf^{-n-1}}$. But the morphism
$X\la \sh[n+1]{\cf^{-n-2}}$ may now be used to form the composite
\[
g\colon X\la\sh[n+1]{\cf^{-n-2}}\la\sh[n+1]{\ck^{-n-2}}=C_{n+1}^{\leq-1}\la C_{n+1}.
\]
The construction is such that the map $(f-g)\colon X\la C_{n+1}$ composes to
zero with the morphism $C_{n+1}\la\sh[n+1]{\cf^{-n-1}}$. This vanishing,
in combination with the distinguished triangle in \eqref{eq:mor1},
permits us to factor the morphism $(f-g)\colon X\la C_{n+1}$ as
$X\la C_n\la C_{n+1}$. By induction, the morphism
$X\la C_n$ must factor as
$X\la\sh[n]{\cf^{-n-1}}\oplus\cf^0\la C_n$. Now we observe that the
composite $\sh[n]{\cf^{-n-1}}\la C_n\la C_{n+1}$ vanishes by the distinguished
triangle \eqref{eq:mor1}, and hence the morphism
$f\colon X\la C_{n+1}$ does indeed factor as
$X\la\sh[n+1]{\cf^{-n-2}}\oplus\cf^0\la C_{n+1}$.
\eprf

\lem{L67.9}
In the setting above, let $\cf^\bullet$ be a complex as in \eqref{eq:complex2} and let $X\in\cs$ be such that, for all $i\in\zz$,
$\Hom(\sh[i]{X},-)$ takes
$\cf^\bullet$ to an exact sequence. Assume further that $\Hom(X,\cs^{\leq-n})=0$
for $n\gg0$.
Then the sequence
\[
\xymatrix@C+0pt{
\Hom(X,\cf^{-1})\ar[r] &\Hom(X,\cf^{0})\ar[r] &
\Hom(X,\ck^{0})\ar[r] & 0\ .
}
\]
is exact.
\elem

\prf
Choose an integer $n>0$ with $\Hom(X,\cs^{\leq-n})=0$, and consider the object
$C_n$ constructed above. Note that it is easy to see that the distinguished triangle
$C_n^{\leq-1}\la C_n\la C_n^{\geq0}\la \sh{C_n^{\leq-1}}$ identifies with
$C_n\la \ck^0\la\sh[n+1]{\ck^{-n-1}}$.

For any morphism
$X\la\ck^0$ we have that the composite $X\la \ck^0\la\sh[n+1]{\ck^{-n-1}}$
must vanish, and hence the morphism must factor as $X\la C_n\la\ck^0$.
\autoref{L67.7} permits us to further factor it through $X\la\sh[n]{\cf^{-n-1}}\oplus\cf^0\la C_n$. But by
	assumption any morphism $X\la\sh[n]{\cf^{-n-1}}$ vanishes, giving us the
	factorization of $X\la\ck^0$ through $X\la\cf^0\la\ck^0$.
	Thus the morphism $\Hom(X,\cf^0)\la\Hom(X,\ck^0)$ is surjective.

Now in $\cs^\heartsuit$ we have a short exact sequence
$0\la\ck^{-1}\la\cf^0\la\ck^0\la 0$, which means that in $\cs$ there is a
distinguished triangle $\ck^{-1}\la\cf^0\la\ck^0\la \sh{\ck^{-1}}$. If we are given a morphism
$X\la\cf^0$ such that the composite $X\la\cf^0\la\ck^0$ vanishes,
then it must factor through $\ck^{-1}$. But now we can apply the paragraph
above to the complex
\[
\xymatrix@C+0pt{
 \cdots \ar[r] & \cf^{-4}\ar[r]  & \cf^{-3}\ar[r] &
 \cf^{-2}\ar[r] &\cf^{-1}
}
\]
to deduce the surjectivity of $\Hom(X,\cf^{-1})\la\Hom(X,\ck^{-1})$.
Assembling this together gives the required exact sequence
\[\xymatrix@C+0pt{
\Hom(X,\cf^{-1})\ar[r] &\Hom(X,\cf^{0})\ar[r] &
\Hom(X,\ck^{0})\ar[r] & 0.
}\]
\eprf

The next result shows that, under mild assumptions, all compact
objects in $\cs$ are strongly pseudocompact.

\lem{L67.11}
Let $\cs$ be a triangulated category, and let
$P\in\Cp$. Assume that, for every
compact object $c\in\cs$, there exists an integer $n>0$
such that $\Hom(c,\sh[n]{P})=0$. Then $\cs^c\subseteq \cs^{spc}_{[P]}$.
\elem

\prf
It is obvious from \autoref{D5.1} that every compact
object is $[P]$--pesudocompact. What needs proof is
part (ii) of \autoref{D67.1}.

Assume therefore that we are given a \tstrs $\tst\cs$ on
the category $\cs$, as in \autoref{D67.1}. This
means that $\cs^{\leq0}$ must belong to $[P]$, the subcategory
$\cs^{\geq0}$ is closed in $\cs$ under coproducts, the heart
$\cs^\heartsuit$ satisfies [AB5], and the inclusion
$\cs^\heartsuit\la\cs$ respects coproducts.

Take any filtered category $I$ and any functor $\fF\colon I\la\cs^\heartsuit$,
and form the cochain complex 
$\cf^\bullet$ as in \autoref{ex:compl}.
In the category $\cs^\heartsuit$ the cohomology of the complex
$\cf^\bullet$ is
concentrated in degree $0$. If $c\in\cs$ is any compact
object and $i\in\zz$ is an integer,
then the functor $\Hom(\sh[i]{c},-)$ takes the complex $\cf^\bullet$
to a complex computing $\colimj$ for the functor
$\Hom(\sh[i]{c},\fF(-))\colon I\la\ab$. As the category of abelian groups
satisfies [AB5] these $\colimj$ vanish for every $j<0$,
and for $j=0$ we obtain $\colim\,\Hom\big(\sh[i]{c},\fF(-)\big)$.
The part about $j=0$ gives us 
the exact sequence
\[
\xymatrix@C+0pt{
\Hom(c,\cf^{-1})\ar[r] &\Hom(c,\cf^{0})\ar[r] &
\colim\,\Hom\big(c,\fF(-)\big)\ar[r] & 0\ ,
}
\]
while the assertion about $j<0$ tells us that,
for every integer $i\in\zz$, the functor
$\Hom(\sh[i]{c},-)$ takes
\[\xymatrix@C+0pt{
 \cdots \ar[r] & \cf^{-3}\ar[r] &
 \cf^{-2}\ar[r] &\cf^{-1}\ar[r] &\cf^{0}
}\]
to an exact sequence. Because $\cs^{\leq0}\in[P]$ and $\Hom(c,\sh[n]{P})=0$
for $n\gg0$, we have that $\Hom(c,\cs^{\leq-n})=0$ for
$n\gg0$. Thus by \autoref{L67.9} the sequence
\[
\xymatrix@C+0pt{
\Hom(c,\cf^{-1})\ar[r] &\Hom(c,\cf^{0})\ar[r] &
\Hom(c,\colim\,\fF)\ar[r] & 0\ .
}
\]
is exact. Combining the above we deduce that the natural map
\[
\xymatrix@C+30pt{
\colim\,\Hom\big(c,\fF(-)\big)\ar[r] &
\Hom(c,\colim\,\fF)
}
\]
is an isomorphism.
\eprf

\pro{P67.13}
Let $\cs$ be a triangulated category, and let
$P\in\Cp$ be such that $\cs$ admits a \tstrs as in part (ii) of \autoref{D67.1}. Assume moreover that, for every
compact object $c\in\cs$, there exists an integer $n>0$
such that $\Hom(c,\sh[n]{P})=0$. Then $\cs^-_{c,[P]}\subseteq \cs^{spc}_{[P]}$.
\epro

\prf
In \autoref{R5.13} we noted that $\cs^-_{c,[P]}\subseteq\cs^{pc}_{[P]}$, thus we need to verify \autoref{D67.1}(ii).
Assume therefore that we are given a \tstrs $\tst\cs$ on
the category $\cs$, satisfying the
hypotheses of  \autoref{D67.1}(ii).
Take any filtered category $I$ and any functor $\fF\colon I\la\cs^\heartsuit$.
If $x\in\cs^-_{c,[P]}$, then there exists a
distinguished triangle $b\la c\la x\la\sh{b}$ in $\cs$ with $c\in\cs^c$ and 
$b\in\cs^{\leq-1}$. But then, for every object
$H\in\cs^\heartsuit\subseteq\cs^{\geq0}$, we have that the map $c\la x$
induces an isomorphism $\Hom(x,H)\la\Hom(c,H)$. \autoref{L67.11}
tells us that $\Hom(c,-)$ commutes with filtered colimits
in $\cs^\heartsuit$, and hence so does $\Hom(x,-)$.
\eprf

\section{The subcategories of $\ct^b$}\label{sec:TbT-}

This section starts with a discussion of strongly pseudocompact objects in weakly approximable triangulated categories, see \autoref{S6}. This will be crucial in \autoref{S107}, where
we prove that the subcategories of $\ct^b$ in
the diagram \eqref{eq:incl} are
intrinsic. And, as the reader will see,
the recipe that works for these
subcategories in $\ct^b$ also describes them in the larger
$\ct^+$.

\subsection{Strongly pseudocompactness in weakly approximable triangulated categories}
\label{S6}

Let $\ct$ be a weakly approximable triangulated category.
In \autoref{S3} we observed that the subcategories $\ct^-$, $\ct^+$ and $\ct^b$ all have preferred equivalence classes of \tstr s with intrinsic
descriptions. Thus for $\cs$ being any of $\ct^b$, $\ct^+$, $\ct^-$
and $\ct$, it becomes interesting to study
the compact, $\big[\cs^{\leq0}\big]$--pseudocompact
and $\big[\cs^{\leq0}\big]$--strongly pseudocompact objects in
$\cs$, where $\tst\cs$
is a \tstrs in the preferred equivalence class.

To simplify the notation, when $\cs$ is a above, we say that an
object is \emph{pseudocompact,} without specifying with respect
to which class of objects in $\Cp(\cs)$, if it is
$\big[\cs^{\leq0}\big]$--pseudocompact, where $\tst\cs$
is a \tstr\ in the preferred equivalence class.
Similarly, an object is  \emph{strongly pseudocompact} if it is
$\big[\cs^{\leq0}\big]$--strongly pseudocompact.
The full subcategories of all pseudocompact (resp.\ strongly pseudocompact) objects in $\cs$ will
be denoted $\cs^{pc}$ (resp.\ $\cs^{spc}$).

\lem{L5.7}
Let $\ct$ be a weakly approximable triangulated category.
Then we have the equalities
\[
(\ct^b)^c=(\ct^b)^{pc}=\ct^b\cap(\ct^-)^{pc}=\ct^+\cap(\ct^-)^{pc}\qquad\text{and}\qquad(\ct^+)^c=(\ct^+)^{pc}=\ct^+\cap\ct^{pc}.
\]
\elem

\prf
First observe that, when $\cs=\ct^?$, where $?=b,+,-,\emptyset$, a set $\{X_\lambda\st\lambda\in\Lambda\}$ of objects of $\cs$ is contained in $P^\perp$ for some $P\in\big[\cs^{\leq0}\big]$ if and only if it is contained in $\cs^{\geq-n}$ for some integer $n>0$.

On the other hand, when $\cs=\ct^?$, where $?=b,+$, by
\autoref{L27.9} if the coproduct of a set 
$\{X_\lambda\st\lambda\in\Lambda\}$ of objects of $\cs$ exists in $\cs$, then again the set must be contained in $\cs^{\geq-n}$ for some integer $n>0$.

All the assertions in the statement follow easily from the two facts above.
\eprf

Not quite so trivial is the next Lemma.

\lem{L5.9}
Let $\ct$ be a weakly approximable triangulated category. Then
we have the equalities
\[
\ct^{pc}=(\ct^-)^{pc}\qquad\text{ and }\qquad
(\ct^+)^c=(\ct^b)^c=(\ct^+)^{pc}=(\ct^b)^{pc}=\ct^+\cap\ct^{pc}
\]
\elem

\prf
First of all, without loss of generality, we choose, in the preferred equivalence class, a \tstr\ $\tst\ct$
such that $\ct^{\geq0}$ is closed under coproducts. Suppose $t\in\ct$ is pseudocompact. Then the truncation morphisms
$t\la t^{\geq n}$ assemble to a single map
$t\la\prod_{n=0}^\infty t^{\geq n}$, but \autoref{L22.3} tells
us that the natural morphism
\[\xymatrix@C+30pt{
\ds\coprod_{n=0}^\infty t^{\geq n}\ar[r] &
\ds\prod_{n=0}^\infty t^{\geq n}
}\]
is an isomorphism. Thus the map $t\la\prod_{n=0}^\infty t^{\geq n}$
factors uniquely as
\[\xymatrix@C+30pt{
t\ar[r]^-\ph &\ds\coprod_{n=0}^\infty t^{\geq n}\ar[r] &
\ds\prod_{n=0}^\infty t^{\geq n}\ ,
}\]
and as $t$ is pseudocompact the map $\ph$ must factor through a
finite subcoproduct. Hence the truncation morphisms
$t\la t^{\geq n}$ must vanish for $n\gg0$, forcing $t$ to belong to
$\ct^-$. That is $\ct^{pc}=\ct^{pc}\cap\ct^-\subseteq(\ct^-)^{pc}$.

Now let $c$ be an object in $(\ct^-)^{pc}$, let $m$ be an integer,
let $\{t_\lambda\st\lambda\in\Lambda\}$ be a set of objects
in $\big(\ct^{\leq m}\big)^\perp=\ct^{\geq m+1}$, and let
\[
\xymatrix@C+30pt{
c\ar[r]^-\psi &\ds\coprod_{\lambda\in\Lambda}t_\lambda^{}
}
\]
be any morphism. Because $c$ belongs to $\ct^-$ we may choose
an integer $n>0$ with $c\in\ct^{\leq n}$, and the morphism
$\psi$ will have to factor through
$\left(\coprod_{\lambda\in\Lambda}t_\lambda^{}\right)^{\leq n}=
\coprod_{\lambda\in\Lambda}t_\lambda^{\leq n}$, where
the equality follows from the fact that, as we explained in \autoref{S701} (see, in particular, \autoref{rmk:propgen}), the truncation functors preserves coproducts. Hence $\psi$ must factor as
\[
\xymatrix@C+30pt{
c\ar[r]^-\rho &\ds\coprod_{\lambda\in\Lambda}t_\lambda^{\leq n}\ar[r] &
\ds\coprod_{\lambda\in\Lambda}t_\lambda^{}\ ,
}
\]
and the $\big[\ct^{\leq0}\big]$--pseudocompactness of
$c$ in the category $\ct^-$ forces $\rho$
to factor through a finite subcoproduct. Hence $\psi$ factors
through a finite subcoproduct, and we deduce that
$(\ct^-)^{pc}\subseteq\ct^{pc}$.

We have therefore proved that $(\ct^-)^{pc}=\ct^{pc}$. The remaining 
equalities in the statement come from
intersecting $(\ct^-)^{pc}=\ct^{pc}$ with $\ct^+$ and
combining with the
equalities of \autoref{L5.7}.
\eprf

For this article the important lemma will be the following.

\lem{L5.29}
Let $\ct$ be a weakly approximable triangulated category. Then
we have the equalities
\[
\ct^{spc}=(\ct^-)^{spc}\qquad\text{ and }\qquad
(\ct^+)^{spc}=(\ct^b)^{spc}=\ct^+\cap\ct^{spc}\ .
\]
\elem

\prf
For an object $c\in\cs$ the extra restriction imposed
by strong pseudocompactness, in addition
to the requirements of ordinary pseudocompactness, is an assertion
about the way $\Hom(c,-)$ behaves on
filtered colimits of objects in the hearts of certain \tstr s
in the preferred equivalence class. 
The point here is that the collection of such hearts is the same,
independent of which $\cs$ we choose in the collection
$\cs\in\{\ct^b,\ct^-,\ct^+,\ct\}$. Thus the equalities of the
current statement follow from the equalities 
\[
\ct^{pc}=(\ct^-)^{pc}\qquad\text{ and }\qquad
(\ct^+)^{pc}=(\ct^b)^{pc}=\ct^+\cap\ct^{pc}
\]
of \autoref{L5.9} by just adding the ``strongly'' restriction in each case.
\eprf

In addition, we can prove the following.

\lem{L5.30}
The subcategory $\ct^{spc}\subseteq\ct$ is thick.
\elem

\prf
In \autoref{C5.13} we saw that $\ct^{pc}\subseteq\ct$ is a thick
subcategory. The extra restriction, required of objects
in $\ct^{spc}\subseteq\ct^{pc}$, is clearly stable under shifts and direct
summands, and it remains to show that it is stable under
extensions.

Assume therefore that $a\la b\la c\la\sh{a}$ is
a distinguished triangle in $\ct$, with $a,c\in\ct^{spc}$. We need to prove
that $b\in\ct^{spc}$. Let $\tst\ct$ be a \tstr\ in the preferred equivalence class, with $\ct^{\ge0}$ closed under coproducts and $\ct^\heartsuit$ an abelian category satisfying [AB5].
We need to show that the functor $\Hom(b,-)$ commutes
with filtered colimits
in $\ct^\heartsuit$. Choose therefore
a filtered category $I$ and a functor $\fF\colon I\la\ct^\heartsuit$.
Now
the functor $\Hom(-,\colim\,\fF)$ is a
homological functor on $\ct$, as is the functor
$\colim\,\Hom\big(-,\fF(i)\big)$.
The natural transformation
\[
\xymatrix@C+40pt{
\colim\,\Hom\big(-,\fF(i)\big)\ar[r] &
\Hom(-,\colim\,\fF)
}
\]
is an isomorphism when evaluated on all suspensions of $a$ and $c$,
and the 5-Lemma tells us that it must also be an
isomorphism on $b$.
\eprf

We are now ready to prove the main result of this section.

\pro{P5.31}
Let $\ct$ be a weakly approximable triangulated category. Then
$\ct^{spc}=\ct^-_c$.
\epro

\prf
\autoref{P67.13} proved the inclusion $\ct^-_c\subseteq\ct^{spc}$,
in generality greater than we need. It remains to prove the
inclusion $\ct^{spc}\subseteq\ct^-_c$.

Without loss of generality, we may  choose and fix a compactly generated \tstrs
$\tstv\ct\ca$
in the preferred equivalence class. \autoref{L69.15} tells
us that the heart $\tsth\ct\ca$ of this \tstrs is an abelian
category satisfying [AB5].
The key will be to prove the following:

\smallskip

\noindent\emph{Claim.} With the notation as above, let $c$ be an object in
$\ct^{spc}\cap\ct^{\leq0}_\ca$. Then there exists a distinguished
triangle $b\la c\la d\la\sh{b}$ with $b\in\ct^c$ and with
$d\in\ct^{spc}\cap\ct^{\leq-1}_\ca$.

\smallskip

Assume that the claim holds
and choose an object $c\in\ct^{spc}=(\ct^-)^{spc}$, where the
equality is by \autoref{L5.29}. We wish to
show that $c\in\ct^-_c$. 
Because $c$ belongs to
$(\ct^-)^{spc}\subseteq\ct^-$, we have that $c\in\ct^{\leq n}_\ca$ for
some $n\geq0$. Shifting if necessary we may assume that
$c\in\ct^{\leq0}_\ca$. Now we apply the claim iteratively. Setting $c=d_0$ we produce
a sequence of morphisms $d_0\la d_1\la d_2\la\cdots$ in
the category $\ct^{spc}$, with $d_n\in\ct^{\leq-n}_\ca$. And by the claim this
can be done in such a way that in the distinguished triangles $b_n\la d_n\la d_{n+1}\la\sh{b_n}$
the object $b_n$ is compact. The octahedral axiom implies that the distinguished triangle $b\la c\la d_n\la\sh{b}$ has $b\in\ct^c$ and
$d_n\in\ct^{\leq-n}$, showing that $c\in\ct^-_c$.

Thus it remains to prove the claim. \autoref{L69.999} tells us that $c^{\geq0}\in\tsth\ct\ca$
must be
the filtered colimit of the finitely presented objects mapping to
it. So there exists a filtered category $I$, and a functor
$\fF\colon I\la\tsth\ct{\ca,c}$, such that $c^{\geq0}=\colim\, \fF$.
But because
$c$ is strongly pseudocompact the map $c\la c^{\geq0}=\colim\,\fF$
must factor through some $\fF(i)\in\tsth\ct{\ca,c}$. But the map
$c\la\fF(i)$ is a morphism from $c$ to
$\fF(i)\in\tsth\ct{\ca,c}\subseteq\ct^{\geq0}_\ca$, and must therefore
factor uniquely through $c^{\geq0}$. We deduce that the identity
map on $c^{\geq0}$ factors through $\fF(i)$, meaning that $c^{\geq0}$
is a direct summand of the object $\fF(i)\in\tsth\ct{\ca,c}$.
\autoref{C69.10} tells us that the category $\tsth\ct{\ca,c}$
is closed under retracts, and we deduce that
$c^{\geq0}\in\tsth\ct{\ca,c}$.
Hence there exists a compact object $b\in\ct^{\leq0}_{\ca,c}$ with
$b^{\geq0}\cong c^{\geq0}$.

Now consider the composite $f\colon b\la b^{\geq0}\cong c^{\geq0}$.
By construction $f^{\geq0}\colon b^{\geq0}\la c^{\geq0}$
is an isomorphism. \autoref{L69.7} yields a compact object $\wt b$, and morphisms $\ph\colon\wt b\la b$
and $g\colon\wt b\la c$ such that
\begin{itemize}
\item
The map $\ph^{\geq0}\colon\wt b^{\geq0}\la b^{\geq0}$ is an isomorphism.  
\item
The triangle below commutes
\[\xymatrix@C+40pt@R-20pt{
  & b^{\geq0}\ar[dd]^{f^{\geq0}}\\
\wt b^{\geq0} \ar[ru]^-{\ph^{\geq0}}\ar[rd]_{g^{\geq0}} &  \\
 & c^{\geq0} 
}\]
\end{itemize}
The commutativity of the triangle, coupled with the fact that
both $\ph^{\geq0}$ and $f^{\geq0}$ are isomorphisms, forces
$g^{\geq0}$ to be an isomorphism. Completing
$g$ to a distinguished triangle $\wt b\la c\la d\la\sh{\wt b}$
we deduce that $d$ must belong to $\ct^{\leq-1}$. Now $\wt b\in\ct^c\subseteq\ct^{psc}$ where the inclusion
is by \autoref{L67.11}. And the fact that in the distinguished triangle both $\wt b$ and $c$ belong
to $\ct^{spc}$, coupled with 
\autoref{L5.30}, tells us that $d\in\ct^{spc}$. 
This completes the proof of the claim and hence of the Proposition.
\eprf

\cor{C5.905}
Let $\ct$ be a weakly approximable triangulated category. Then
\[
\ct^{spc}=(\ct^-)^{spc}=\ct^-_c\qquad\text{ and }\qquad
(\ct^+)^{spc}=(\ct^b)^{spc}=\ct^b_c\ .
\]
\ecor

\prf
It follows directly from \autoref{P5.31} and \autoref{L5.29}.
\eprf

\subsection{The intrinsic subcategories of $\ct^b$ and $\ct^+$}
\label{S107}

The discussion in the previous section allows us to prove the following.

\lem{lemTbcT+b}
Let $\ct$ be a weakly approximable triangulated category and let $\cs=\ct^?$, with $?=+,b$. Then $\ct_c^b$ is an intrinsic subcategory of $\cs$
\elem

\prf
\autoref{C5.905} yields the equalities $\ct^b_c=(\ct^b)^{spc}=(\ct^+)^{spc}$. The claim is then obvious.
\eprf

The following result shows that also $\ct^{c,b}$ is an intrinsic subcategory both of $\ct^b$ and $\ct^+$.

\lem{L107.1}
Let $\ct$ be a weakly approximable triangulated category and let $\cs=\ct^?$, with $?=+,b$. Let
$\tst\cs$ be a \tstr\ in the preferred equivalence class. An object
$c\in\ct^b_c\subseteq\cs$
belongs to $\ct^{c,b}$
if and only if there exists an integer $n>0$ with
$\Hom(c,\cs^{\leq-n})=0$. 
\elem

\prf
If $c$ belongs to $\ct^c\cap\ct^b_c\subseteq\ct^c$ then there exists
an integer $n>0$ with $\Hom(c,\ct^{\leq-n})=0$, whence
$\Hom(c,\cs^{\leq-n})=0$.

Now assume that we have an object $c\in\ct^b_c\subseteq\cs$, and that
we are given an integer $n>0$ with $\Hom(c,\cs^{\leq-n})=0$.
If $t\in\ct^{\leq-n-1}$ is any object, then for each $\ell>n+1$ the
object $t^{\geq-\ell}$ belongs to $\cs^{\leq-n-1}$, whether $\cs=\ct^b$
or $\cs=\ct^+$. Then, for all $\ell>n+1$,
\[
\Hom\Big(c,t^{\geq-\ell}\Big)\eq0\eq\Hom\Big(c,\sh[-1]{\big(t^{\geq-\ell}\big)}\Big)\ .
\]
By \cite[Proposition~3.2]{Neeman18} the weakly approximable
triangulated category $\ct$ is left-complete, meaning that $t$
is isomorphic to $\holim t^{\geq-\ell}$. Hence there exists in $\ct$ a distinguished triangle
\[
\xymatrix@C+30pt{
\ds\prod_{\ell=n+2}^\infty\sh[-1]{\big(t^{\geq-\ell}\big)}
\ar[r] & t\ar[r] &
\ds\prod_{\ell=n+2}^\infty t^{\geq-\ell}\ar[r]& \ds\prod_{\ell=n+2}^\infty\big(t^{\geq-\ell}\big).
}
\]
Since $\Hom(c,-)$ is trivial when evaluated on the product terms in the triangle, we must have $\Hom(c,t)=0$.
As $t\in\ct^{\leq-n-1}$ is arbitrary we conclude that
$\Hom(c,\ct^{\leq-n-1})=0$.

But now $c$ belongs to $\ct^b_c\subseteq\ct^-_c$, and \autoref{L2.1}
together with the vanishing of $\Hom(c,\ct^{\leq-n-1})$ implies
that $c\in\ct^c$.
\eprf

\section{The case $\ct^c\subseteq\ct^b_c$: a procedure
  to recover $\ct^c$}
\label{S112}

Let $\ct$ be a weakly approximable triangulated category. As we said in the introduction, we do not in general understand how to recognise in $\ct^b_c$ which objects belong to the subcategory $\ct^{c,b}=\ct^c\cap\ct^b_c$.
There are, however, two special cases in which we can say something useful:
\be
\item[(a)]
If $\ct$ is \emph{coherent}.
\item[(b)]
If $\ct^c\subseteq\ct^b_c$ and furthermore $^\perp(\ct^b_c)\cap\ct^-_c=\{0\}$.
\ee
For an exposition of the algorithm in
(a), in full generality, the reader is referred to  
\cite{Neeman18A}.

In \autoref{subsect:gencase} we will concern ourselves
with the precise statement and  proof of (b). The recipe
for $\ct^{c,b}=\ct^c$ is in
\autoref{P112.791}. An example of a $\ct$, to which (b) applies, is $\ct=\Dqcs Z(X)$, with $X$ a
quasi-compact, quasi-separated scheme and with $Z\subseteq X$
a closed subset with quasi-compact complement (see \autoref{P1000.1}). Furthermore, the hypothesis $^\perp(\ct^b_c)\cap\ct^-_c=\{0\}$ also holds for any coherent weakly approximable triangulated category. Therefore, if we specialize to the case where
$\ct^c\subseteq\ct^b_c$, then (b) gives an approach to
(a) different form \cite[Propositions~5.6 and 6.5]{Neeman18A}.
See \autoref{subsec:coherent} for further discussion.

\subsection{Classical generators and special equivalence classes}\label{subsec:generaltech}

In this section we develop general results which will help to recognize the subcategory $\ct^{c.b}$ of $\ct^b_c$.
We start from a definition which is very similar, in spirit, to \autoref{D3.1}.

\dfn{D112.1}
Let $\cs$ be a triangulated category, and consider the class
$\cq=\cq(\cs)$
of all full subcategories $P\subseteq\cs$ satisfying $P\subseteq\sh{P}$.
\be
\item
Two elements $P,Q\in\cq$ are \emph{equivalent}
if there exists an integer $A>0$ with $\sh[-A]{P}\subseteq Q\subseteq \sh[A]{P}$.
\item
Given two equivalence classes $[P]$ and $[Q]$ of elements of $\cq$,
then $[P]\leq[Q]$ if, for a choice of representatives
$P\in[P]$ and $Q\in[Q]$, there exists an integer $A>0$ with 
$\sh[-A]{P}\subseteq Q$.
\ee
\edfn

The difference between the $\Cp(\cs)$ of \autoref{D3.1}
and the $\cq(\cs)$ of \autoref{D112.1}
is that the $P\in\Cp(\cs)$ are assumed to
satisfy $\sh{P}\subset P$, while the $P\in\cq(\cs)$ must have
$P\subseteq\sh{P}$. Next we focus on the $P\in\Cp(\cs)$ and $Q\in\cq(\cs)$
that will interest us in this
section.

Given $H\in\cs$, we define the following two subcategories of $\cs$:
\begin{equation}\label{eq:specialsub}
P_H(\cs):=H[0,\infty)^\perp\qquad Q_H(\cs):=H(-\infty,0]^\perp.
\end{equation}
When there is no confusion about the triangulated category $\cs$, we use the shorthands $P_H=P_H(\cs)$ and $Q_H=Q_H(\cs)$.

\rmk{R112.4}
For any object $H\in\cs$ we have the inclusions and equalities
\[
\sh{P_H}\subseteq P_H,\quad P_H*P_H=P_H,\quad
\add(P_H)=P_H,\quad\Smr(P_H)=P_H
\]
as well as
\[
Q_H\subseteq\sh{Q_H},\quad Q_H*Q_H=Q_H,\quad
\add(Q_H)=Q_H,\quad\Smr(Q_H)=Q_H.
\]
If $X$ is an object of $P_H$, it follows from the
first set of inclusions that $\genul X{}0\subseteq P_H$.
If $X$ is an object of $Q_H$, the second set of
inclusions gives that $\genuf X{}0\subseteq Q_H$.
\ermk

For an object $H$ of a triangulated category $\cs$ we may ask the following hypothesis to be satisfied:

\hyp{H112.11}
There exists an integer $A>0$ such that:
\be
\item  
For every object $X\in\cs$ there exists an integer $B>0$, with $X\in\sh[-B]{P_H}\cap\sh[B]{Q_H}$.
\item
$\Hom(\sh[A]{P_H},Q_H)=0$.
\item
For every object $F\in P_H$ and every integer $m>0$,
there exists a distinguished triangle $E\la F\la D$ in $\cs$ with
$E\in\genu H{}{1-m-A,A}$ and with
$D\in\sh[m]{P_H}$.
\ee
\ehyp

\rmk{R112.98}
We take a bit of time to elaborate on \autoref{H112.11}(iii). As stated, this item starts with an object $F\in P_H$ and an
integer $m>0$. We claim that, in fact, the general case can be deduced from the special
case where $m=1$, as follows:
\begin{itemize}
\item
Take $F\in P_H$. By the case $m=1$ of \autoref{H112.11}(iii), we may
construct a distinguished triangle $E_1\la F\la D_1$ with $E_1\in\genu H{}{-A,A}$ and with
$D_1\in\sh{P_H}$.
\item
Now we proceed by induction. Suppose we have constructed the sequence
\[\xymatrix{
	F\ar@{=}[r]&D_0\ar[r] & D_1\ar[r] &D_2\ar[r] &\cdots\ar[r] &D_{n-1}\ar[r] & D_n,
}\]
with $D_i\in\sh[i]{P_H}$ and such that, in the distinguished triangle $\wt E_i\la D_i\la D_{i+1}$,
we have $\wt E_i\in \genu H{}{-i-A,-i+A}$. Then the previous step allows us to 
continue this a step further.
\end{itemize}
It is easy to see that, if we complete the composite $F\la D_m$ to a distinguished triangle
$E_m\la F\la D_m$, then $E_m\in\genu  H{}{1-m-A,A}$.
Moreover, for any
$i$ in the interval $0<i<m$, we have that the $E_m$ we constructed
satisfies
\[
E_m\in\genu H{}{1-i-A,A}*\genu H{}{1-m-A,-i+A}.
\]
More precisely, in the distinguished triangle $E_i\la E_m\la\wt D$ we have
$E_i\in\genu H{}{1-i-A,A}$ and $\wt D\in\genu H{}{1-m-A,-i+A}$.
\ermk

The hypothesis above looks technical at first sight. The first aim of this section is to study which objects do satisfy it.

\lem{L112.5}
Let $\ct$ be a weakly approximable triangulated category,
and let $\tst\ct$ be a \tstr\ on $\ct$ in the preferred
equivalence class and assume that
$\ct^c\subseteq\ct^b_c$.

If $G\in\ct^c\subseteq\ct^b_c$ is such that $\ct^c=\langle G\rangle$, then
\[
[P_G(\ct^b_c)]=\big[\ct^b_c\cap\ct^{\leq0}\big]
  \qquad\text{ and }\qquad
[Q_G(\ct^b_c)]=\big[\ct^b_c\cap\ct^{\geq0}\big]\ .
\]
\elem

\nin
Let us remind that an object $G\in\ct^c$ such that $\ct^c=\langle G\rangle$ is called a \emph{classical generator} for
$\ct^c$. And it is automatic that a classical generator
$G\in\ct^c$ compactly
generates the larger category $\ct$.

\prf
Since we only want to compute the equivalence classes of
$P_G=P_G(\ct^b_c)$ and $Q_G=Q_G(\ct^b_c)$, we are free to replace
the \tstr\ $\tst\ct$ by the equivalent \tstrs $\tstv\ct G$.
In this case
$\ct_G^{\leq0}=\Coprod\big(G(-\infty,0]\big)$,  providing the second equality in
\[
\ct_G^{\geq1}=\Big(\ct_G^{\leq0}\Big)^\perp=\Coprod\big(G(-\infty,0]\big)^\perp=G(-\infty,0]^\perp,
\]
where the orthogonals are taken in $\ct$ and the last equality is due to \autoref{R701.1}(ii). Hence
the equality $Q_G=\ct^b_c\cap\ct_G^{\geq1}$ comes by intersecting with $\ct^b_c$.

Now we turn to computing the equivalence
class of $P_G$. By \autoref{dfn:wa}(i) there exists
an integer $A>0$ such that
$G\in{^\perp\ct^{\leq-A}}$. Hence
$G[0,\infty)\subset{^\perp\ct^{\leq-A}}$, or equivalently
\[
\ct^{\leq0}\subseteq G[A,\infty)^\perp
  \]
where the orthogonals are taken in $\ct$.
The inclusion $\ct^b_c\cap\ct^{\leq0}\subseteq\sh[-A]{P_G}$
comes from intersecting with $\ct^b_c$.
On the other hand,
\cite[Lemma~3.9(iv)]{Burke-Neeman-Pauwels18}
permits us to choose an integer $B>0$ such that
\[
G[-B,\infty)^\perp\subseteq\ct^{\leq0} ,
\]
and intersecting with $\ct^b_c$ gives the inclusion
$\sh[B]{P_G}\subseteq\ct^b_c\cap\ct^{\leq0}$.
\eprf

\cor{C112.7}
Let $\ct$ be a weakly approximable triangulated category and assume that
$\ct^c\subseteq\ct^b_c$. If $G\in\ct^c\subseteq\ct^b_c$ is such that $\ct^c=\langle G\rangle$, then $\Hom(\sh[A]{P_G},Q_G)=0$ for $A\gg0$.
\ecor

\prf
Let $\tst\ct$ be a \tstr\ on $\ct$ in the preferred
equivalence class.
By \autoref{L112.5} there exists an integer $B>0$ such that
\[
\sh[B]{P_G}\subseteq\ct^b_c\cap\ct^{\leq0}\qquad\text{ and }\qquad
\sh[-B]{Q_G}\subseteq\ct^b\cap\ct^{\geq0}.
\]
But as $\Hom\big(\sh[n]{(\ct^{\leq0})},\ct^{\geq0}\big)=0$ for $n>0$, we deduce
that $\Hom(\sh[A]{P_G},Q_G)=0$ for $A>2B$.
\eprf

Putting all these results together, we get that classical generators of $\ct^c$ satisfy our technical assumption.

\pro{L112.9}
Let $\ct$ be a weakly approximable triangulated category and assume that
$\ct^c\subseteq\ct^b_c$.
If $G\in\ct^c\subseteq\ct^b_c$ is such that $\ct^c=\langle G\rangle$, then $G$ satisfies \autoref{H112.11} (with $\cs=\ct^b_c$).
\epro

\prf
We begin with (i) in  \autoref{H112.11}. Let $\tst\ct$ be a \tstrs
in the preferred equivalence class,
and let $X\in\ct^b_c$ be an object.
Because $\ct^b_c$ is contained in $\ct^b$, the object
$X\in\ct^b_c$ must lie in $\ct^{\leq n}\cap\ct^{\geq-n}$ for some $n\ge0$.
By \autoref{L112.5} there exists an integer $B>0$ with
\[
\ct^b_c\cap\ct^{\leq n}\subseteq \sh[-n-B]{P_G}
\qquad\text{ and }\qquad
\ct^b_c\cap\ct^{\geq-n}\subseteq \sh[n+B]{Q_G}
\]
and hence
\[
X\in\ct^b_c\cap\ct^{\leq n}\cap\ct^{\geq-n}
\subseteq
\sh[-n-B]{P_G}\cap\sh[n+B]{Q_G} .
\]

Item (ii) of  \autoref{H112.11} directly follows from \autoref{C112.7}.
So let us deal with item (iii) in the hypothesis. Replacing
the \tstrs $\tst\ct$ by an equivalent one if necessary, we may
choose an integer $A$ such that
$\ct^b_c\cap\ct^{\leq-A}\subseteq P_G\subseteq\ct^b_c\cap\ct^{\leq0}$.
And now \cite[Corollary~2.14]{Neeman17A} allows us to choose
an integer $A'>0$ such that, for
any integer $m>0$, any object
$F\in P_G\subseteq\ct^b_c\cap\ct^{\leq0}$ admits a
distinguished triangle $E\la F\la D$, with
$D\in\ct^b_c\cap\ct^{\leq-A-m}\subseteq \sh[m]{P_G}$ and with
$E\in\genu G{}{1-m-A-A',A'}$.
\eprf

\subsection{The general criterion}\label{subsect:gencase}

Let us put ourselves in the following context. Let $\ct$ be a weakly approximable triangulated category, assume that $\ct^c$
is contained in $\ct^b_c$, and let $H\in\cs=\ct^b_c$ be an object satisfying 
\autoref{H112.11}.
Choose a compactly generated \tstr\ $\tst\ct$ in the preferred equivalence
class, and let $\ch\colon\ct\la\tsth\ct{}$ be the induced homological functor.

Let $F\in P_H$ be any object and, for $m>0$, let
\begin{equation}\label{eq:triarmk}
E_m\la F\la D_m\la\sh{E_m}
\end{equation}
be the distinguished triangles in \autoref{H112.11}(iii). We assume we constructed them according to \autoref{R112.98}. In particular, they form sequences which allow us to define
\[
E:=\hoco E_m\qquad\text{ and }\qquad D:=\hoco D_m.
\]
Again, for each $m>0$, we are given the natural morphisms
$E_m\la E$ and $D_m\la D$ which we can complete in $\ct$ to the distinguished triangles
\[
E_m\la E\la \wt E_m\qquad\text{ and }\qquad D_m\la D\la \wt D_m.
\]

\lem{L112.75}
In the setting above, the following hold:
\be
\item
There exists an integer $B>0$ such that, for all
$m>0$, both $\wt E_m$ and $\wt D_m$
lie in $\ct^{\leq-m+B}$.
\item
The objects
$E$ and $D$
belong to $\ct^-_c$.
\item
There is in $\ct^-_c$
a distinguished triangle $E\la F\la D\la\sh{E}$.
\ee
\elem

\prf
We are assuming that $H$ belongs to $\ct^b_c\subseteq\ct^-$, and hence we
may assume, up to shift, that $H\in\ct^{\leq0}$. Let $F\in P_H$.

Take the triangle \eqref{eq:triarmk} constructed in
$\ct^b_c$ in \autoref{R112.98}. By the same remark, for any pair of integers $m<n$, the morphism $E_m\la E_n$ may be completed to
a distinguished triangle $E_m\la E_n\la\wt D_{m,n}$ with
$\wt D_{m,n}\in\genu H{}{1-n-A,-m+A}\subseteq\ct^{\leq-m+A}$.
Fix the integer $m>0$;
the long exact sequence in cohomology associated to such a triangle tells us that
$\ch^i(E_m)\la\ch^i(E_n)$ is an isomorphism for all $i>-m+A+1$. Now
form the morphism of triangles
\[\xymatrix{
E_m\ar[r]\ar[d] & F\ar@{=}[d]\ar[r] & D_m\ar[r]\ar[d] &\sh{E_m}\ar[d] \\
E_n\ar[r] & F\ar[r] & D_n\ar[r] &\sh{E_n}
}\]
and, by applying $\ch$ to it,
we learn that the induced morphism $\ch^i(D_m)\la\ch^i(D_n)$ is also an isomorphism in the
range $i>-m+A+1$. Taking the colimit as $n\la\infty$ and
combining  with
\cite[Lemma~1.4(iv)]{Neeman17A} we deduce that, for $i>-m+A+1$,
the functor $\ch^i$
takes the natural maps $E_m\la E$ and $D_m\la D$ to isomorphisms.
The long exact sequence in cohomology, applied to the distinguished triangles
$E_m\la E\la \wt E_m$ and $D_m\la D\la \wt D_m$, tells us that
$\wt E_m$ and $\wt D_m$ both belong to $\ct^{\leq-m+A+1}$,
proving (i).

By construction, the objects
$E_m$ and $D_m$ both
belong to $\ct^b_c\subseteq\ct^-_c$, and from
the distinguished triangles $E_m\la E\la \wt E_m$ and
$D_m\la D\la \wt D_m$ we now deduce that
$E$ and $D$ belong to $\ct^-_c*\ct^{\leq-m+A+1}$. On the other hand
$\ct^-_c=\bigcap_{n=1}^\infty\big(\ct^c*\ct^{\leq-n}\big)$,
and hence $E$ and $D$ both belong to
\[
\ct^-_c*\ct^{\leq-m+A+1}
\subseteq \big(\ct^c*\ct^{\leq-m+A+1}\big)*\ct^{\leq-m+A+1}
=\ct^c*\big(\ct^{\leq-m+A+1}*\ct^{\leq-m+A+1}\big)
=\ct^c*\ct^{\leq-m+A+1}.
\]
As $m>0$ is arbitrary, we deduce that $D,E\in\ct^-_c$, proving (ii).

Now consider the commutative square
\[\xymatrix@C+50pt{
\ds\coprod_{n=1}^\infty D_n\ar[d]\ar[r]^-{1-\shi} &\ds\coprod_{n=1}^\infty D_n\ar[d] \\
\ds\coprod_{n=1}^\infty \sh{E_n}\ar[r]_-{1-\shi}
& \ds\coprod_{n=1}^\infty \sh{E_n}
}\]
and complete to a morphism of triangles
\[\xymatrix@C+50pt{
\ds\coprod_{n=1}^\infty D_n\ar[d]\ar[r]^-{1-\shi} &\ds\coprod_{n=1}^\infty D_n\ar[d] 
\ar[r] & D\ar[d]^\ph\\
\ds\coprod_{n=1}^\infty \sh{E_n}\ar[r]_-{1-\shi}
& \ds\coprod_{n=1}^\infty \sh{E_n}\ar[r] & \sh{E}.
}\]
Next we may complete $\ph$ to a distinguished triangle $E\la \wt F\la D\stackrel\ph\la\sh{E}$. For every integer $m>0$, we have the commutative diagram
\[
\xymatrix@C+30pt{
E_m\ar[r] & F\ar[r] & D_m\ar[r]\ar[d] & \sh{E_m}\ar[d] \\
E  \ar[r] &\wt F\ar[r] & D\ar[r]^-\ph & \sh{E}
}\]
which we can complete to a morphism of triangles
\[\xymatrix@C+30pt{
E_m\ar[r]\ar[d] & F\ar[r]\ar[d]^-{\lambda_m} & D_m\ar[r]\ar[d] & \sh{E_m}\ar[d] \\
E  \ar[r] &\wt F\ar[r] & D\ar[r]^-\ph & \sh{E}.
}
\]
Applying the functor $\ch$ gives a map of long exact sequences, and as
$\ch^i$ takes $E_m\la E$ and $D_m\la D$ to isomorphisms for all $i>-m+A+1$,
we deduce that $\ch^i$ takes the map $\lambda_m$ to an isomorphism
whenever $i>-m+A+2$. Since $F$ is assumed to belong to $\ct^b_c\subseteq\ct^b$,
there exists an integer $B>0$ with $\ch^i(F)=0$, for all $i<-B$.
For any
integer $i<-B$, choose an integer $m>0$ with $-m+A+2<i$, and the isomorphism
$\ch^i(\lambda_m)\colon\ch^i(F)\to\ch^i(\wt F)$ establishes that $\ch^i(\wt F)\cong\ch^i(F)=0$. Thus $\ch^i(\wt F)$ also vanishes for all $i<-B$.

Therefore for any integer $m>B+A+2$ we have that, for every $i\in\zz$,
the functor $\ch^i$ takes the map $\lambda_m\colon F\la\wt F$ to an isomorphism.
For $i\geq -B$ this is because $m$ is chosen big enough, and for $i<-B$ it
is because $\ch^i(F)\cong0\cong\ch^i(\wt F)$.
By \cite[Lemma~2.7]{Burke-Neeman-Pauwels18} the
\tstr\ $\tst\ct$ is non-degenerate (meaning that the intersections of all the $\ct^{\leq n}$ and of all the $\ct^{\geq n}$ are trivial),
and hence the morphism $\lambda_m\colon F\la\wt F$ must be
an isomorphism in $\ct$ (see \cite[Proposition 1.3.7]{BeiBerDel82}). 
Replacing 
the triangle $E\la\wt F\la D\stackrel\ph\la\sh{E}$ by an isomorphic
triangle
$E\la F\la D\stackrel\ph\la\sh{E}$ completes the proof of (iii).
\eprf


In conclusion, for every
$F\in\ct^b_c$, we obtain a distinguished triangle $E\la F\la D\la\sh{E}$ in $\ct^-_c$. Since we need to remember the dependence on $H$ and
on $F$, we will write $E_{(F,H)}$ and $D_{(F,H)}$ for the
objects $E$ and $D$.

\lem{L112.975}
In the setting at the beginning of this section, if $F\in P_H$, then $\ct^b_c\subseteq\big(D_{(F,H)}^{}\big)^\perp$.
\elem

\prf
Let $X\in\ct^b_c$ be any object, and we need to prove the
vanishing of $\Hom(D,X)$, where $D=D_{(F,H)}$. Because $X$ belongs to $\ct^b_c$ and $H$ satisfies
\autoref{H112.11}(i), there exists an integer $B>0$
with $X\in \sh[B]{Q_H}$. Because
\autoref{H112.11}(ii) holds for $H$, we may (after increasing
the integer $B>0$) assume that $\Hom(\sh[B]{P_H},X)=0$.
The construction of the triangle \eqref{eq:triarmk} is such that
$D_m\in\sh[m]{P_H}$, and hence for all $m>B$ we have
$\Hom(D_m,X)=0$.

Complete the morphism $D_m\la D$ to a distinguished triangle
\begin{equation}\label{eq:t1}
D_m\la D\la\wt D_m\la\sh{D_m}.
\end{equation}
By \autoref{L112.75}(i) we may, after increasing the
integer $B$, assume that $\wt D_m\in\ct^{\leq-m+B}$.
But $X$ belongs to $\ct^b_c\subseteq\ct^+$, and hence
$X\in\ct^{\geq-N}$, for some $N>0$. Hence, for
all $m>B+N$, we have $\wt D_m\in\ct^{<-B-N+B}=\ct^{<-N}$
and so $\Hom(\wt D_m,X)=0$.

Choose therefore any $m>B+N$. In the distinguished triangle \eqref{eq:t1} we have
$\Hom(D_m,X)=0=\Hom(\wt D_m,X)$. Hence
$\Hom(D,X)=0$.
\eprf

We can now prove the main result of this section which gives an intrinsic description of $\ct^c$ in $\ct_c^b$ under a technical assumption.

\pro{P112.791}
Let $\ct$ be a weakly approximable triangulated category
such that $\ct^c\subseteq\ct^b_c$
and suppose that $^\perp(\ct^b_c)\cap\ct^-_c=\{0\}$.
Then $[P_G]$ is maximal with respect to the partial order $\leq$, among the classes $[P_H]$ with
$H\in\cs=\ct^b_c$ satisfying \autoref{H112.11}.

Furthermore, an object $X\in\ct^b_c$ belongs to
the subcategory $\ct^c\subseteq\ct^b_c$ if and
only if, for some object $H\in\ct^b_c$
satisfying \autoref{H112.11} and with $[P_H]$ maximal,
we have $\Hom(X,P_H)=0$.
\epro

\prf
Let $H\in\ct^b_c$ be an object satisfying \autoref{H112.11}. Up to shift, we can assume $G\in P_H$, thanks to \autoref{H112.11}(i). So we can form in
$\ct^-_c$ the distinguished triangle
\[
E_{(G,H)}^{}\la G\la D_{(G,H)}^{}\la\sh{E_{(G,H)}}
\]
of \autoref{L112.75}(iii). By
\autoref{L112.975} we have that $D_{(G,H)}^{}\in\ct^-_c$
is an object with
$\big(D_{(G,H)}\big)^\perp$ containing $\ct^b_c$. By hypothesis
$^\perp(\ct^b_c)\cap\ct^-_c=\{0\}$, and hence $D_{(G,H)}=0$.
Therefore the map from $E=E_{(G,H)}\la G$ is an
isomorphism. Thus we can factor the identity of $G$ as the composite $G\la E\la G$.

On the other hand,
by assumption, we have that $E=\hoco E_m$,
with each $E_m$ in some
$\genu H{}{1-m-A,A}$. The map
from the compact object $G$ to $E=\hoco E_m$ factors
through some $E_m\la E$, allowing
us to write the identity on $G$ as a composite
$G\la E_m\la G$. Therefore $G$ is a direct summand
of an object  $E_m$, and hence belongs to $\genu H{}{-B,B}$ for some integer $B>0$. We deduce that $G[0,\infty)\subseteq\genuf H{}{-B}$.
Taking orthogonals in $\ct^b_c$ now gives the inclusion
\[
\sh[B]{P_H}=H[-B,\infty)^\perp=
\left(\genuf H{}{-B}\right)^\perp\subseteq G[0,\infty)^\perp=P_G.
\]
This shows $[P_H]\leq[P_G]$. The first part of the statement is then proven as \autoref{L112.9} shows that $G$ satisfies \autoref{H112.11}.

To prove the second part of the statement, suppose that $X\in\ct^b_c$
is such that $\Hom(X,P_H)=0$ for some $H$ with $[P_H]$
maximal. Combining the first part of the statement with \autoref{L112.5}
we have that $X^\perp$ contains $\ct^b_c\cap\ct^{\leq-B}$
for some $B>0$. Since $G\in\ct^c$ is a
classical generator, \cite[Corollary~2.14]{Neeman17A}
gives the existence of a distinguished triangle
\begin{equation}\label{eq:t2}
E\la X\la D\la\sh{E},
\end{equation}
with
$E\in\gen G{}$ and $D\in\ct^{\leq-B}$. As we are assuming
that $\ct^c\subseteq\ct^b_c$ we have that $E\in\ct^c$ must
belong to $\ct^b_c$, and the distinguished triangle \eqref{eq:t2}, coupled
with the fact that both $E$ and $X$ belong to $\ct^b_c$,
tells us that $D\in\ct^b_c$ as well. Hence
$D$ belongs to $\ct^b_c\cap\ct^{\leq-B}$, and so the map $X\la D$
in \eqref{eq:t2} must vanish. Thus
$X$ is a direct summand of the object $E\in\ct^c$ and $X$ must belong to $\ct^c$.

To complete the proof note that, if $X\in\ct^c$ and $G\in\ct^c\subseteq\ct^b_c$ is a classical generator, then
$X\in\genu G{}{-B,B}$ for some $B>0$. Now, by definition, $^\perp P_{\sh[B]{G}}$ contains
$\genu G{}{-B,B}$, and hence $X^\perp$ contains
$P_{\sh[B]{G}}$, with $[P_{\sh[B]{G}}]=[P_G]$ maximal by the first part.
\eprf

\section{Examples and applications}\label{sec:examples}

The first part of this section is about two geometric examples where the assumption $^\perp(\ct^b_c)\cap\ct^-_c=\{0\}$ is automatically verified: $\ct=\D(\Mod{R})$ for $R$ a coherent ring or, more generally, $\ct$ a coherent weakly approximable triangulated category, and $\ct=\Dqcs Z(X)$, with $X$ quasi-compact and quasi-separated with $Z\subseteq X$ a closed subset with quasi-compact complement. This is done in \autoref{subsec:coherent} and \autoref{S1000}, respectively.

\subsection{The coherent case}\label{subsec:coherent}

The first situation where \autoref{P112.791} applies is provided by weakly approximable triangulated categories of the following special type.

\dfn{dfn:coherent}
A weakly approximable triangulated category $\ct$ is \emph{coherent} if, for any
\tstrs $\tst\ct$ in the preferred equivalence class, there
exists an integer $N>0$ such that every
object $Y\in\ct^-_c$ admits a distinguished triangle
$X\la Y\la Z$ with
$X\in\ct^-_c\cap\ct^{\leq N}$ and with $Z\in\ct^b_c\cap\ct^{\geq0}$.
\edfn

\exm{ex:coherent}
Let $\ct$ be a weakly approximable triangulated category.
Suppose there exists, in the preferred equivalence class, a \tstr\
$\tst\ct$ such that $\big(\ct^-_c\cap\ct^{\leq0},\ct^-_c\cap\ct^{\geq0}\big)$
is a \tstr\ on $\ct^-_c$. Then it is automatic that $\ct$ is coherent.
After all: for every object $Y\in\ct^-_c$, the truncation triangles
$Y^{\leq0}\la Y\la Y^{\geq1}$ will satisfy $Y^{\leq0}\in\ct^-_c\cap\ct^{\leq0}$ and
$Y^{\geq1}\in\ct^-_c\cap\ct^{\geq1}=\ct^b_c\cap\ct^{\geq1}$

Concrete examples abound. If there exists a left coherent ring $R$
and a triangulated equivalence of categories $\ct\cong\D(\Mod R)$, then $\ct$
is weakly approximable, the
standard \tstr\ on $\D(\Mod R)$ is in the preferred
equivalence class, and it restricts to a \tstr\ on
$\ct^-_c=\D^-(\mmod R)$. Similarly: if there
exists a coherent and quasi-separated scheme $X$, with
$Z\subseteq X$ a closed subscheme such that $X\setminus Z$ is quasi-compact,
and a triangulated equivalence of
categories $\ct\cong\Dqcs{Z}(X)$,
then $\ct$ is a weakly approximable, the
standard \tstr\ on $\Dqcs{Z}(X)$ is in the
preferred equivalence class,
and it restricts to a \tstr\ on the subcategory
$\ct^-_c=\dcohs Z(X)$.
\eexm

\pro{prop:coherent2}
If $\ct$ is a weakly approximable and coherent triangulated category such that $\ct^c\subseteq\ct^b_c$, then the
condition $^\perp(\ct^b_c)\cap\ct^-_c=\{0\}$ of
\autoref{P112.791} is satisfied.
\epro

\prf
Let  $\tst\ct$ be a \tstrs
in the preferred equivalence class and assume that $D\in\ct^-_c$ is an object with
$D^\perp\supseteq\ct^b_c$.
By the coherence assumption, we get distinguished triangles
\[
X_m\la\sh[m]{D}\la Z_m\la\sh{X_m},
\]
with $X_m\in\ct^-_c\cap\ct^{\leq N}$ and $Z_m\in\ct^-_c\cap\ct^{\geq0}\subseteq\ct^b_c$. But then the morphisms $\sh[m]{D}\la Z_m$ vanish, for all $m\in\zz$ and $\sh[m]{D}$ is
a direct summand of $X_m\in\ct^{\leq N}$.
Therefore $D$ must belong to
$\bigcap_{m=0}^\infty\ct^{\leq-m}=\{0\}$, where
the vanishing is by
\cite[Lemma~2.7]{Burke-Neeman-Pauwels18}. Thus $D=0$ and we are done.
\eprf

\subsection{The case of $\Dqcs Z(X)$}
\label{S1000}

Let $X$ be a quasi-compact, quasi-separated scheme, and let $Z\subseteq X$
be a closed subset with quasi-compact complement. The category $\Dqcs Z(X)$
is weakly approximable by \autoref{ex:wa}(ii). On the other hand, \cite[Theorem~3.2(i)]{Neeman22A}
says that  the compact objects in $\Dqcs Z(X)$ are all perfect complexes. In particular they have bounded cohomology. Thus we are in the situation where $\ct^c\subseteq\ct^b_c$

We can now prove the following key result.

\pro{P1000.1}
The assumption $^\perp(\ct^b_c)\cap\ct^-_c=\{0\}$ holds for
$\ct=\Dqcs Z(X)$.
\epro

\prf
Let $D\in\ct^-_c$ be an object such that $D^\perp$ contains
$\ct^b_c$, and hence $\ct^c$. In order to prove that $D=0$, it is enough to show that, for all open
immersions $i\colon W\hookrightarrow X$ with $W\subseteq X$ an affine open subset,
we have $i^*D=0$.

For  $i\colon W\hookrightarrow X$ as above, \cite[Lemma~7.5(i)]{Neeman22A}
shows that, for any
compact object $F\in\Dqcs{Z\cap W}(W)$ and any compact
generator $G\in\Dqcs Z(X)$, there exists an integer
$M>0$ with
$i_*F\in\ogenu G{}{-M,M}$. But then
any morphism
$D\la i_*F$ is
a map from $D\in\ct^-_c$ to
$i_*F\in\ogenu G{}{-M,M}$, and \cite[Lemma~2.7]{Neeman18} tells
us that it
must factor as $D\la E\la i_*F$ for
some object $E\in\genu G{}{-M,M}$.
On the other hand, the map $D\la E$ must vanish because $D^\perp$ contains
$\genu G{}{-M,M}\subseteq\ct^c$. Hence, for every
compact object $F\in\Dqcs{Z\cap W}(W)$, we have
\[
\Hom(i^*D,F)\cong\Hom(D,i_*F)\cong0.
\]
The output is that we can reduce to the case where $X$ is an affine scheme.

Let $j\colon X\setminus Z\la X=\spec R$ be the open immersion,
and complete
the unit of adjunction $\co_X^{}\la j_*j^*\co_X^{}$ to
a distinguished
triangle
\begin{equation}\label{eq:tr3}
L'\la \co_X^{}\la j_*j^*\co_X^{}\la\sh{L'}.
\end{equation}
Choose any compact generator $H\in\Dqcs{Z}(X)$. By \cite[Lemma~7.4]{Neeman22A}
there exists an integer $r>0$ such that
$L'\in\ogenu H{}{-r,r}\subseteq\Dqcs{Z}(X)$.
Suppose $D\la \sh[n]{L'}$ is some morphism. Since it is a morphism
from an object in $\ct^-_c$ to
an object in $\ogenu H{}{-r-n,r-n}$,
\cite[Lemma~2.7]{Neeman18} allows us to factor it through an object $F\in\genu H{}{-r-n,r-n}\subseteq\ct^c$.
As $D^\perp$ contains $\ct^c$, the map $D\la F$ and therefore the composite $D\la F\la\sh[n]{L'}$ must vanish.
As a consequence, the functor $\Hom(D,-)$ applied to $\sh[n]{L'}$ is trivial, for all $n\in\zz$. On the other hand, 
\[\Hom(D,\sh[n]{j_*j^*\co_X})\cong\Hom(j^*D,\sh[n]{j^*\co_X})=0
\]
because $D\in\Dqcs Z(X)$ and hence $j^*D\in j^*\Dqcs Z(X)=\{0\}$. From
\eqref{eq:tr3}, we deduce that $\Hom(D,-)$ applied to $\sh[n]{\co_X}$ is trivial, 
for all $n\in\zz$. Under the
equivalence $\ct\subseteq\Dqc(X)\cong\D(R)$, the complex
$D\in\ct^-_c$ is a bounded above cochain complex
of finitely generated projective $R$-modules
and we have just
proved that $\Hom(D,\sh[n]{R})=0$ for
all $n\in\zz$. Applying \autoref{thm:appendix-main}
to the complex $P^\bullet:=\Hom(D,R)$,
we deduce
that $D=0$.
\eprf




\appendix

\section{Appendix by Christian Haesemeyer}
\label{A1}

In this appendix, we prove the following:

\begin{thmApp}\label{thm:appendix-main}
Let $R$ be a commutative ring, and $P^\bullet$ a bounded below cochain complex of finitely generated projective $R$-modules. If $P^\bullet$ is acyclic, then it is chain contractible.
\end{thmApp}

\begin{remarkApp}\label{RA1.1}
The statement also holds for (not necessarily commutative) rings with bounded finitistic dimension. (Noetherian local commutative rings fall into this class, by the Auslander - Buchsbaum formula.) For more on this, see \cite{Shaul23}. However, the statement is false, in general, over non-commutative rings (see \cite{Positselski23}).
\end{remarkApp}

In order to prove the theorem, it suffices to show that the cocycle modules $Z^n(P^\bullet)$ are projective for all $n$. Since they are finitely presented by hypothesis, we can and will assume that $R$ is local  and we are dealing with a complex of free modules. We employ the following definition, due to Hochster \cite{Hochster74}; we will use the version of the definition given by Northcott in \cite[Section 5.5]{Northcott76}. We remind the reader that the term ``grade" is often called ``depth" in recent literature; we will stick with grade to align with the references we use. 

\begin{definApp}\label{defn:grade}
Let $I$ be an ideal in $R$, and $M$ an $R$-module. The {\em true} or {\em polynomial} grade of $I$ on $M$, denoted $\mathrm{Gr}(I;M)$, is the limit (possibly $\infty$) of the increasing sequence
\[
\mathrm{gr} (I R[x_1,\dots,x_n]; R[x_1,\dots,x_n]\otimes_R M),
\]
where $\mathrm{gr} (I R[x_1,\dots,x_n]; R[x_1,\dots,x_n]\otimes_R M)$ denotes the (classical) grade, that is the upper bound on the length of $R[x_1,\dots,x_n]\otimes_R M$-regular sequences in $I R[x_1,\dots,x_n]$ (again, possibly infinite). Following usual practice, we write $\mathrm{Gr}(I)$ for $\mathrm{Gr}(I;R)$.
\end{definApp}
		
Northcott proves (see \cite[Chapter 5, Theorems 9 and 13]{Northcott76}):

\begin{lemmaApp}\label{lem:grade_vs_generators}
Let $I\subseteq R$ be a proper ideal that can be generated by $n$ elements. Then $\mathrm{Gr}(I)\leq n$. 
\end{lemmaApp}

Given an $R$-homomorphism $\phi\colon E\to F$ of finite rank free $R$-modules, write $\mathfrak{A}_\nu(\phi)$ for the determinantal ideal of $R$ generated by the $\nu\times\nu$-minors of $\phi$ (this ideal does not depend on a choice of basis of the modules $E$ and $F$). The rank of $\phi$ is defined as the maximal $\nu$ such that $\mathfrak{A}_\nu(\phi)\neq 0$ (see \cite[Equation (3.2.1)]{Northcott76}); we denote the determinantal ideal with $\nu$ the rank of $\phi$ by $\mathfrak{A}(\phi)$.

Now suppose we are given a (chain) complex $0\to F_n\xrightarrow{d_n} F_{n-1}\xrightarrow{d_{n-1}} \dots \xrightarrow{d_1} F_0$ of finite rank free $R$-modules which is exact except, possibly, at $F_0$. 
Northcott proves the following result (see \cite[Chapter 6, Theorem 14]{Northcott76}); the noetherian version of this theorem is due to Buchsbaum and Eisenbud  \cite{Buchsbaum-Eisenbud73}.

\begin{thmApp}\label{thm:BE}
For $0< i\leq n$, we have the inequality $\mathrm{Gr}(\mathfrak{A}(d_i))\geq i$.
\end{thmApp}

With this in hand, we are ready to prove \autoref{thm:appendix-main}.

\prf (of \autoref{thm:appendix-main})
As discussed above, we may assume that $R$ is local, and the complex $P^\bullet$ is a complex of free modules of finite rank; it then suffices to show show that $Z^k(P^\bullet)$ is free for all integers $k$. Consider the differential $d^{k-2}\colon P^{k-2}\to P^{k-1}$. Its determinantal ideal $\mathfrak{A}(d^{k-2})$ is finitely generated, say, by $s$ generators. Applying \autoref{thm:BE} to the complex $P^\bullet$, brutally truncated in degree $k + s$, we conclude that $\mathfrak{A}(d^{k-2})$ must have true grade greater than $s$, and hence, cannot be a proper ideal by \autoref{lem:grade_vs_generators}. By \cite[Chapter 6, Exercise 11]{Northcott76}, the cokernel of $d^{k-2}$, that is, the cocycle module $Z^k(P^\bullet)$, is free. 	
\eprf

\bibliographystyle{abbrv}

\end{document}